\documentclass{article}

\usepackage{arxiv}

\usepackage[utf8]{inputenc} 
\usepackage[T1]{fontenc}    
\usepackage{hyperref}       
\usepackage{url}            
\usepackage{amsfonts}       
\usepackage{nicefrac}       
\usepackage{microtype}      
\usepackage{lipsum}

\usepackage{amsmath}
\usepackage{amssymb}
\usepackage{amsfonts}
\usepackage{amstext}

\usepackage[utf8]{inputenc}
\usepackage{chngcntr}
\usepackage{appendix}
\usepackage{amsthm}

\usepackage{xcolor, etoolbox}

\usepackage{hyperref}

\hypersetup{
   colorlinks=true,
   urlbordercolor={1 1 1},
   citecolor={blue}
}

\usepackage[T1]{fontenc}    
\usepackage{hyperref}       
\usepackage[left]{lineno}

\usepackage{url}            
\makeatletter
\g@addto@macro{\UrlBreaks}{\UrlOrds}
\makeatother

\usepackage{booktabs}       
\usepackage{amsfonts}       
\usepackage{nicefrac}       
\usepackage{microtype}      
\usepackage{lipsum}

\usepackage{amsfonts}
\usepackage{relsize}
\usepackage{comment}
\usepackage{listings}

\usepackage{multirow}
\usepackage{tabulary,capt-of}
\usepackage{adjustbox}
\usepackage{diagbox}

\usepackage{algpseudocode}
\usepackage{multicol}
\usepackage{caption}

\usepackage{tabularx}
\usepackage{pdfpages}
\usepackage{graphicx}
\usepackage{caption, subcaption}
\usepackage{booktabs}

\usepackage{mathtools}
\newtagform{blue}{\color{red}(}{)}
\newtagform{redandblue}[\textcolor{red}]{\color{red}(}{)}

\usetagform{blue}

\newtheorem{theorem}{Theorem}

\newtheorem{corollary}[theorem]{Corollary}

\newtheorem{example}[theorem]{Example}

\newtheorem{lemma}[theorem]{Lemma}

\newtheorem{proposition}[theorem]{Proposition}
\newtheorem{remark}[theorem]{Remark}

\newcounter{unnumber}

\usepackage{multirow}
\usepackage{colortbl}
\usepackage{multicol}

\usepackage{mathpazo}
\usepackage{physics}

\usepackage[ruled,vlined]{algorithm2e}

\title{A theoretical and empirical study of new adaptive algorithms with additional momentum steps and shifted updates for stochastic non-convex optimization}

\author{
 Cristian Daniel Alecsa \\
  Romanian Institute of Science and Technology \\
  Technical University of Cluj-Napoca\\
  Cluj-Napoca, Romania \\
  \texttt{alecsa.cd@gmail.com} \\ \\
  \textbf{ORCID }: $0000$-$0003$-$2324$-$7711$
}

\begin{document}
\nocite{*}

\maketitle
\begin{abstract}
It is known that adaptive optimization algorithms represent the key pillar behind the rise of the Machine Learning field. In the Optimization literature numerous studies have been devoted to accelerated gradient methods but only recently adaptive iterative techniques were analyzed from a theoretical point of view. In the present paper we introduce new adaptive algorithms endowed with momentum terms for stochastic non-convex optimization problems. Our purpose is to show a deep connection between accelerated methods endowed with different inertial steps and \emph{AMSGrad}-type momentum methods. Our methodology is based on the framework of stochastic and possibly non-convex objective mappings, along with some assumptions that are often used in the investigation of adaptive algorithms. In addition to discussing the finite-time horizon analysis in relation to a certain final iteration and the almost sure convergence to stationary points, we shall also look at the worst-case iteration complexity. This will be followed by an estimate for the expectation of the squared Euclidean norm of the gradient. Various computational simulations for the training of neural networks are being used to support the theoretical analysis. For future research we emphasize that there are multiple possible extensions to our work, from which we mention the investigation regarding non-smooth objective functions and the theoretical analysis of a more general formulation that encompass our adaptive optimizers in a stochastic framework. 
\end{abstract}

\section{Preliminaries}\label{Section_Preliminaries}
We will consider in the present paper the following unconstrained expectation minimization problem
\begin{align}\label{OptPb}\tag{OptPb}
\min\limits_{\vb* \lambda \in \mathbb{R}^d} F(\vb* \lambda) = \mathbb{E}[f(\vb* \lambda; \xi)],    
\end{align}
where the mapping $f : \mathbb{R}^d \times \mathbb{R}^m \to \mathbb{R}$ is continuously Fr\' echet differentiable and possibly non-convex. The expectation $\mathbb{E}[\cdot]$ is taken with respect to $\xi$, where $\xi$ denotes a random vector with an unknown distribution $\mathbb{P}$. We assume that the gradients of the objective function $F$ can only be accessed through noisy information. The so-called expected risk with respect to the minimization problem \eqref{OptPb} is defined as $R(\vb* \lambda) = \mathbb{E}[f(\vb* \lambda; \xi)]$. On the other hand, the empirical risk which is often used as an approximation for the expected risk, takes the form
\begin{align}\label{EmpiricalRisk}\tag{EmpRisk}
R_N(\vb* \lambda) = \dfrac{1}{N} \sum\limits_{i = 1}^{N} f_i(\vb* \lambda),    
\end{align}
where the mapping $f_i : \mathbb{R}^d \to \mathbb{R}$ is the loss function taken with respect to the $i^{th}$ sample, and $N$ represents the number of the sample data. By considering a sample set $\lbrace (\vb x_i, \vb t_i) \rbrace_{i=1}^{N}$, the functions $f_i$ can be written as $f_i(\vb* \lambda) = l(h(\vb x_i, \vb* \lambda); \vb* t_i)$, where $l: \mathbb{R}^{d_t} \times \mathbb{R}^{d_t} \to \mathbb{R}$ is a general loss function and $h : \mathbb{R}^{d_x} \times \mathbb{R}^d \to \mathbb{R}^{d_t}$ is a prediction function. Here, $\mathbb{R}^{d_x}$ represents the input space, while $\mathbb{R}^{d_t}$ is the output space. The general objective of Machine Learning problems is to find a prediction function $h \in \mathcal{H}$ (where $\mathcal{H}$ represents the family of prediction functions) that minimizes the loss function $l$, which relates to inaccurate predictions between the true output and the expected output.

\subsection{Related works}\label{Subsection_RelatedWorks}
\textbf{Stochastic gradient descent-type optimizers} \\
It is well known that optimization techniques can be divided into two major groups from a practical standpoint: stochastic and batch methods. In the fundamental work of Robbins and Monro \cite{RobbinsMonro}, the classical Stochastic Gradient Descent method (briefly, \emph{SGD}) was introduced as a Markov chain. The examination carried out at \cite{Bottou}, in which the conventional asymptotic analysis of \emph{SGD} was done, is an exhaustive overview in the context of learning from i.i.d. observations. Non-asymptotic properties, in the wider context of Hilbert spaces, have been discussed in \cite{BachMoulines} as part of an extended theoretical analysis of SGD-type algorithms. However, in \cite{GhadimiLanFirstOrder}, Ghadimi and Lan have presented a finite-time analysis of a randomized first-order algorithm in a non-convex setting where a given random variable with a certain probability distribution plays the role of a final iteration. For preconditioned stochastic first-order algorithms, a continuation of this investigation was conducted in \cite{StochasticQuasiNewton}, with some additional findings related to worst-case iteration complexity. According to \cite{LiOrabona}, an adaptive version of SGD was used, for which an almost sure convergence of the gradients to zero was shown. Sebbouh et. al. in \cite{DefazioHeavyBall} have shown the almost sure convergence rates for the minimum squared gradient norm along the trajectories of the \emph{SGD} and for the heavy ball method for non-convex objective functions. They also provided convergence rates for the iterates values in the objective function. \\ \\
\textbf{Accelerated first-order algorithms} \\ 
We are going to focus now on accelerated first-order stochastic methods. From the work of Polyak \cite{Polyak}, who introduced the heavy-ball method, to the pioneering work of Nesterov \cite{Nesterov}, which introduced the first accelerated inertial method, the first-order algorithms with an accelerated rate of convergence in the iterates of the objective function have spawned an entire research field concerning optimization algorithms. Yan et. al. provided a unified framework for momentum methods and they have developed the convergence analysis in the stochastic scenario in \cite{UnifiedMomentum}. We also like to point out that in the article \cite{GhadimiHeavyBall} the investigation of the convergence analysis of the heavy ball method was studied, but only for convex objective functions. In order to deal with non-convex mappings, Ochs et. al. \cite{Ochs} have designed the \emph{IPiano} algorithm, which is a generalization of the heavy ball algorithm in the form of a forward-backward splitting method. Although Nesterov's accelerated gradient method (\emph{NAG}) was originally proposed for solving convex problems, in \cite{LSCNesterov} L\' aszl\' o has considered the convergence of the aforementioned algorithm in a deterministic non-convex context. Also, for non-convex and stochastic optimization problems, the convergence properties of generalized approaches related to \emph{NAG} have been independently investigated in \cite{GhadimiLanAccelerated} by Ghadimi and Lan. \\ \\
\textbf{First-order adaptive methods} \\
The emergence of Machine Learning has led to the popularity of momentum-type methods and recent studies like \cite{Sutskever} have shown that inertial algorithms can be used to tackle large scale non-convex problems. Notwithstanding these inertial algorithms, the true rejuvenation behind the study of qualitative properties of neural networks is often related to the case of adaptive optimizers. Just to name a few of them, first-order adaptive algorithms like \emph{RMSProp} \cite{RMSProp}, \emph{Adadelta} \cite{Adadelta}, \emph{Adagrad} \cite{Adagrad} and \emph{Adam} \cite{Adam} were originally proposed for solving unconstrained optimization problems, and are currently the powerhouse behind the empirical performances of neural networks. \\
The most well-known of these adaptive iterative approaches is the \emph{Adam} algorithm, which constitutes a turning point moment in the history of numerical optimization algorithms.
The unscaled version of this algorithm was recently studied in \cite{BarakatBianchi} concerning non-convex deterministic and stochastic optimization problems. In addition, the more complex case of the scaled \emph{Adam} method endowed with bias correction terms has just been studied from a dynamical system point of view. More precisely, in \cite{DaSilva}, the \emph{Adam} method was shown to be the forward Euler discretization of a dynamical system which was thoroughly studied in the convex setting (see also the work \cite{BarakatBianchiDynSys} for the non-convex setting). 
Despite the fact that this adaptive dynamical system includes the one for which the heavy ball is a numerical discretization, the natural dynamical system that arises as a continuous version of \emph{NAG} was only recently explored in the framework of convex objective functions in \cite{ALP}. \\
Quite recently, in \cite{ReddiAdamConvergence}, \emph{Adam} was shown that it does not convergence even on simple convex examples, and the authors have proposed \emph{AMSGrad} to solve the issues linked to its convergence. Zaheer et. al. \cite{Zaheer} have showed a $\mathcal{O}(1/n)$ convergence rate for increasing mini-batch sizes, but the proof was merely presented for the \emph{RMSProp} algorithm. Through techniques related to convergence rates with high probability bounds, in \cite{ZhouAMSGrad} Zhou et. al. have considered a theoretical study underlying the capabilities of adaptive methods and it improved the dependency of the convergence results with respect to the work of Chen et. al. \cite{ChenAdam}. Also, in \cite{ZouAdam} the authors have investigated the convergence of \emph{Adam} and \emph{RMSProp} but, as in many research papers, the theoretical bounds related to hyper-parameters depend on the dimension of the space. On the other hand, Li and Orabona \cite{LiOrabona} proved the convergence rate for the \emph{AdaGrad-Norm} algorithm assuming the classical assumption concerning Lipschitz continuity of the gradient and the boundedness of the variance. This differs from the work \cite{BottouAdagrad}, where the authors did not required knowledge of the Lipschitz constant, but they assume that the gradient is uniformly bounded by some finite value. It is worth emphasizing the recent work \cite{DefossezAdam} of D\' efossez et. al. in which they have proved the convergence of some variants of \emph{Adagrad} and \emph{Adam} algorithms for smooth and non-convex optimization problems under the assumption of an almost surely uniform bound of the norm of the gradients. 
There are other various recent research works devoted to the analysis of adaptive-type methods. As an example, Baggioni and Tarrago in \cite{BaggioniAdaptive} considered a theoretical study regarding the non-asymptotic analysis of adaptive gradient methods for strongly convex objective functions, by focusing on applications related to linear regression and regularized generalized linear models for \emph{Adagrad} and \emph{stochastic Newton} iterative processes, respectively. A different line of work is that of Xiao et al. \cite{XiaoAdam} which proved some convergence guarantees for \emph{Adam}-type methods in the setting of non-convex and possibly non-smooth objective functions, and they provided applications on the training of the \emph{ResNet-50} neural network model for image classification tasks on the benchmark datasets \emph{CIFAR10} and \emph{CIFAR100}, respectively. 
Finally, we note that Gadat and Gavra recently published in \cite{GadatGavra} a rigorous analysis regarding the theoretical convergence of single-scaled \emph{RMSProp}-type algorithms in a general non-convex stochastic setting.

\setlength{\tabcolsep}{3pt}
\renewcommand{\arraystretch}{0.4}
\begin{table}[hb!]
\centering
\begin{tabular}{p{45mm} p{30mm} p{30mm} p{30mm} p{20mm} p{30mm}}
\bottomrule \bottomrule \\ [0.5ex] 
\texttt{Articles} & \texttt{Convexity} & \texttt{Smoothness} & \texttt{Setting} & \texttt{Method} \\ [0.5ex]
\hline \\ [0.5ex]
Ghadimi $\&$ Lan (2016) \cite{GhadimiLanAccelerated} & non-convex & smooth & stochastic & algorithm \\ \\ [0.35ex]
Wang et al. (2017) \cite{StochasticQuasiNewton} & non-convex & smooth & stochastic & algorithm \\ \\ [0.35ex]
Chen et al. (2018) \cite{ChenAdam} & non-convex & smooth & stochastic & algorithm \\ \\ [0.35ex]
Alecsa et al. (2019) \cite{ALV} & non-convex & smooth & deterministic & algorithm \\ \\ [0.35ex]
Alecsa et al. (2020) \cite{ALP} & convex & smooth & deterministic & dyn. sys. \\ \\ [0.35ex] 
Barakat $\&$ Bianchi (2020) \cite{BarakatBianchi} & non-convex & smooth & deterministic \newline stochastic & algorithm \\ \\ [0.35ex] 
Da Silva $\&$ Gazeau (2020) \cite{DaSilva} & convex \newline non-convex & smooth & deterministic & dyn. sys. \\ \\ [0.35ex]
Barakat $\&$ Bianchi (2021) \cite{BarakatBianchiDynSys} & non-convex & smooth & stochastic & dyn. sys. \newline algorithm \\ \\ [0.35ex]
Gadat $\&$ Gavra (2022) \cite{GadatGavra} & non-convex & smooth & stochastic & algorithm \\ \\[0.35ex] 
Xiao et al. (2023) \cite{XiaoAdam} & non-convex & non-smooth & stochastic & dyn.sys \newline algorithm \\ \\ [0.5ex]
\bottomrule \bottomrule
\end{tabular}
\vspace{+0.5em}
\caption{Recent research contributions}\label{table_recent_research}
\end{table}

We end the present subsection by describing the investigation that was done in different papers which are in connection to our theoretical results. For this we refer to table \eqref{table_recent_research} where we have introduced some reference articles along with the following properties that were studied: \emph{Convexity} (if the objective function is convex or non-convex), \emph{Smoothness} (if the objective function has a differentiability property like $C^1$ or $C^2$), \emph{Setting} (which consists of the deterministic or stochastic framework) and \emph{Method} (which shows if dynamical systems or algorithms were studied).

\subsection{Notations}\label{Subsection_Notations}
Through the present paper, scalars (along with functions and operators) are denoted by lower and upper case letters, and (random) vectors are denoted by lower $\&$ upper case bold letters. On the other hand, even though the random vector $\xi$ (or any elements which are related to this) with respect to the stochastic optimization problem \eqref{OptPb} are random vectors, we will denote them with usual lower case letters, in order to distinguish them from the random vectors which will represent the iterates of our algorithms. Further, we consider the following entrywise Hadamard notations that are similar to the ones used in \cite{DaSilva}. For a vector $\vb u \in \mathbb{R}^d$, the $i^{th}$ component is written as $u_{[i]}$. Hence, for a sequence of vectors $(\vb u_n)_{n \in \mathbb{N}} \subset \mathbb{R}^d$, we denote by $u_{j, [i]} \in \mathbb{R}$ the $i^{th}$ component of $\vb u_j$. Now, given two vectors $\vb u = (u_{[1]}, \ldots, u_{[d]})$ and $\vb v = (v_{[1]}, \ldots, v_{[d]})$ from $\mathbb{R}^d$ and $\delta, \rho \in \mathbb{R}$, we consider the following vectors in $\mathbb{R}^d$:
$\vb u + \rho = (u_{[1]} + \rho, \ldots, u_{[d]} + \rho)$,
$\vb u \rho := \vb u \cdot \rho = (u_{[1]} \cdot \rho, \ldots, u_{[d]} \cdot \rho)$,
$\vb u \odot \vb v = (u_{[1]} \cdot v_{[1]}, \ldots, u_{[d]} \cdot v_{[d]})$, $\dfrac{\vb u}{\vb v} = \left( \dfrac{u_{[1]}}{v_{[1]}}, \ldots, \dfrac{{u}_{[d]}}{v_{[d]}} \right)$,
$[\vb u]^\delta = \left( u_{[1]}^\delta, \ldots, u_{[d]}^\delta \right)$,
$|\vb u| = (|u_{[1]}|, \ldots, |u_{[d]}|) \text{ and } \sqrt{\vb{u}} = (\sqrt{u_{[1]}}, \ldots, \sqrt{u_{[d]}})$. At the same time, the vector-type inequality (evidently, this also stands also for equality identities) $\vb u \leq \rho$ means that $u_{[i]} \leq \rho$ for every $i \in \lbrace 1, \ldots, d \rbrace$. Evidently, this notations is valid also for the reverse inequality, i.e. $\vb u \geq \rho$ (and also for non-strict vector-type inequalities).
Furthermore, a vector $\vb u \in \mathbb{R}^d$ is positive if and only if $\vb u > 0$, while $\vb u \geq 0$ represents that $\vb u \in \mathbb{R}^d$ is non-negative (these notations also stand for scalar variables). \\   
Moreover, throughout our article, a constant vector $\textbf{c} \in \mathbb{R}^d$ is viewed as the vector with equal components, i.e. $\textbf{c} = ( \underbrace{c, \ldots, c}_{\text{d times}})$, where $c \in \mathbb{R}$. Also, $\| \cdot \|$ will always represent the Euclidean norm $\| \cdot \|_2$ in the finite dimensional space $\mathbb{R}^d$. When discussing the characteristics of random variables or random vectors, we shall use the notation \emph{w.p.1} to denote \emph{'with probability one'} or, alternatively, \emph{'almost surely'}. Consequently, random variables/vectors with values in the finite dimensional space $\mathbb{R}^d$ are subject to all the concepts and notations for vectors in $\mathbb{R}^d$.

\subsection{Organization of the paper}
In section \eqref{Section_AdaptiveHeavyBall} we introduce different formulations of our adaptive methods. More precisely, the optimization algorithm \eqref{AcceleratedGradientMethod} can be perceived as an extension of the algorithm presented in \cite{GhadimiLanAccelerated}, while \eqref{AdaptiveHeayBallMethod} (which is equivalent with the previous method \eqref{AcceleratedGradientMethod}) represents an adaptive version of accelerated gradient methods. On the other hand, our last iterative formulation is the adaptive momentum type method \eqref{AdaptiveAcceleratedMomentumMethodShiftedUpdates} which is inspired by the recent work \cite{BarakatBianchi}. In section \eqref{Section_AssumptionSetting} we present the assumptions used by us in all of the remaining sections: boundedness of the objective function, Lipschitz continuity of the gradient of the objective function, the property of the stochastic gradient to be an unbiased estimate of the true gradient of the loss (in the shifted updates), bounded variance in connection with the stochastic gradient, measurability of certain random vectors and bounded stochastic gradients (used often in the case of adaptive methods), respectively. In section \eqref{Section_ConvergenceResults} we present our theoretical results from which we remind the descent type property presented in Lemma \eqref{Lemma_AlmostSureConvergence}, the convergence rate analysis with respect to a chosen final iteration along with an important asymptotic almost sure inequality considered in Theorem \eqref{Theorem_AsymptoticBehavior}, and the long time convergence of the gradient of the objective function	presented in Remark \eqref{Remark_GradientConvergence}. In section \eqref{Section_ExamplesAdaptiveMethods} we introduce the following: in example \eqref{Example_1} we consider a new type of adaptive methods using the idea of effective stepsize \eqref{a_n} and accumulated squared gradients \eqref{v_n} (along with the underlying theoretical rate of convergence), and in proposition \eqref{Proposition_ComplexityResult} we investigate a complexity result for a particular case of our adaptive momentum algorithm with respect to a given final time iteration. In section \eqref{EquivalentForms} we show in \eqref{AdaptiveSutskever} a practical reformulation of our adaptive methods, and in Algorithm \eqref{AlgDescription_Sutskever_AAMMSU} entitled \emph{Sutskever formulation of AAMMSU} we consider our associated \texttt{PyTorch} algorithm. On the other hand, in section \eqref{Section_NeuralNetworks} we present our simulations concerning classification tasks involving different Machine Learning models, while in section \eqref{Conclusions} we outline some perspectives and limitations regarding our work. Finally, the Appendix section \eqref{Section_Appendix} contains the proofs of our theoretical results.

\subsection{The proposed adaptive algorithm}

In this subsection we consider some preliminary results regarding the comparison of our method presented in Algorithm \eqref{AlgDescription_Sutskever_AAMMSU} entitled \emph{Sutskever formulation of AAMMSU}, along with various adaptive optimizers. Also, we will track down the construction of our adaptive algorithm in a step by step fashion.

\begin{figure}[h!]
\centering
    \begin{subfigure}[b]{\textwidth}
        \includegraphics[width=\textwidth, height=0.28\textheight]{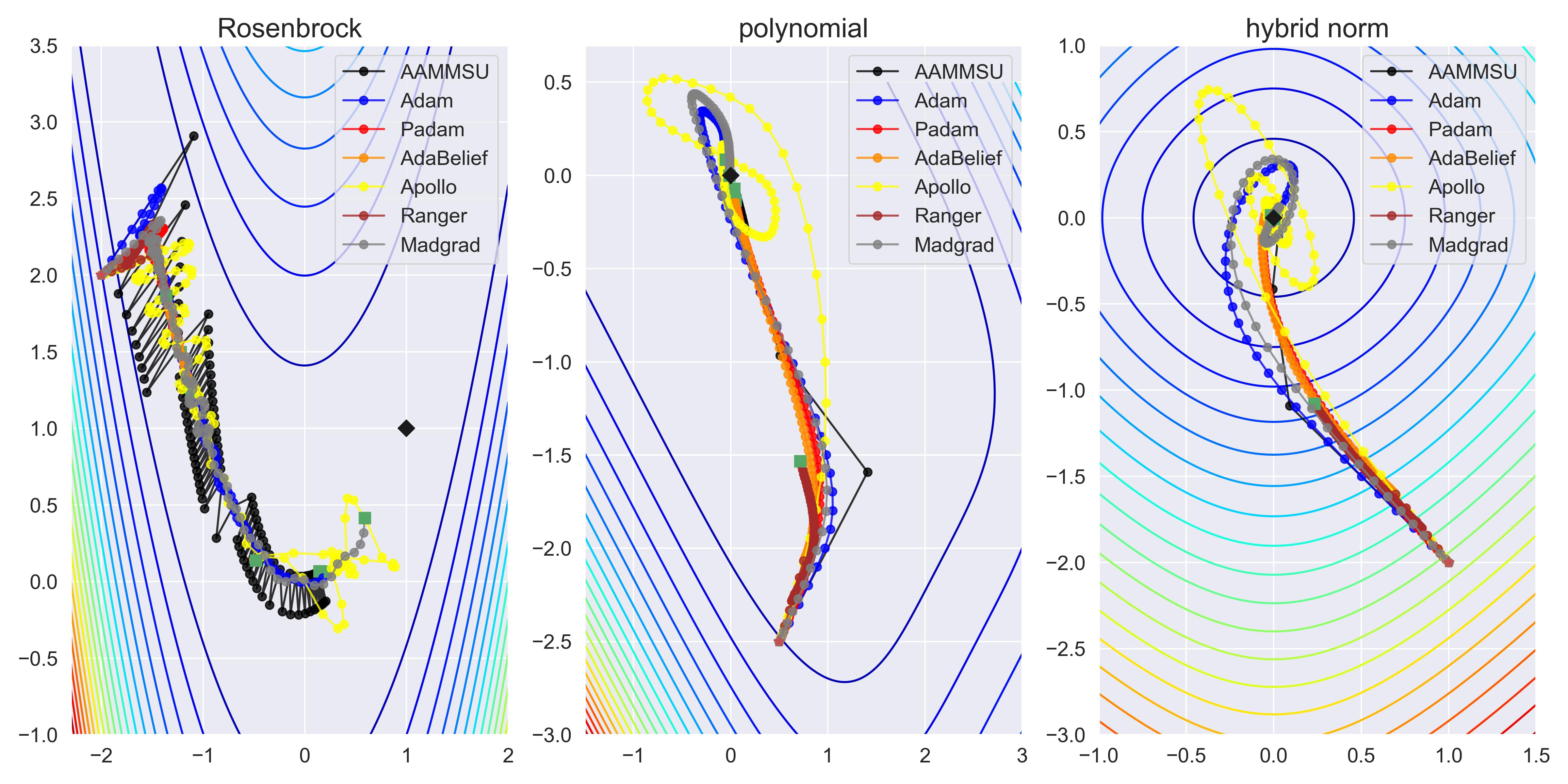}
    \end{subfigure}
    \begin{subfigure}[b]{\textwidth}
        \includegraphics[width=\textwidth, height=0.28\textheight]{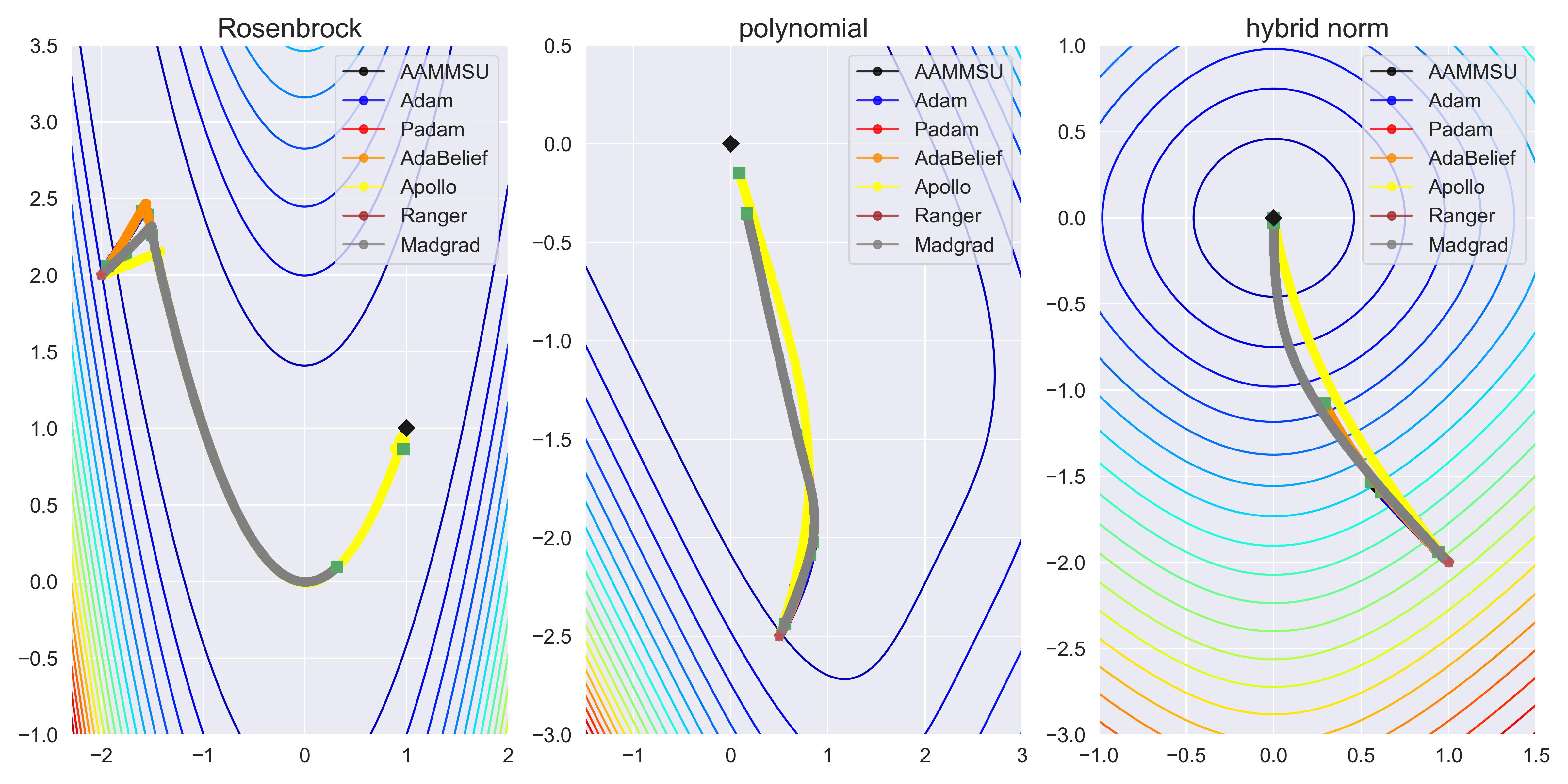}
    \end{subfigure}
    \begin{subfigure}[b]{\textwidth}
        \includegraphics[width=\textwidth, height=0.28\textheight]{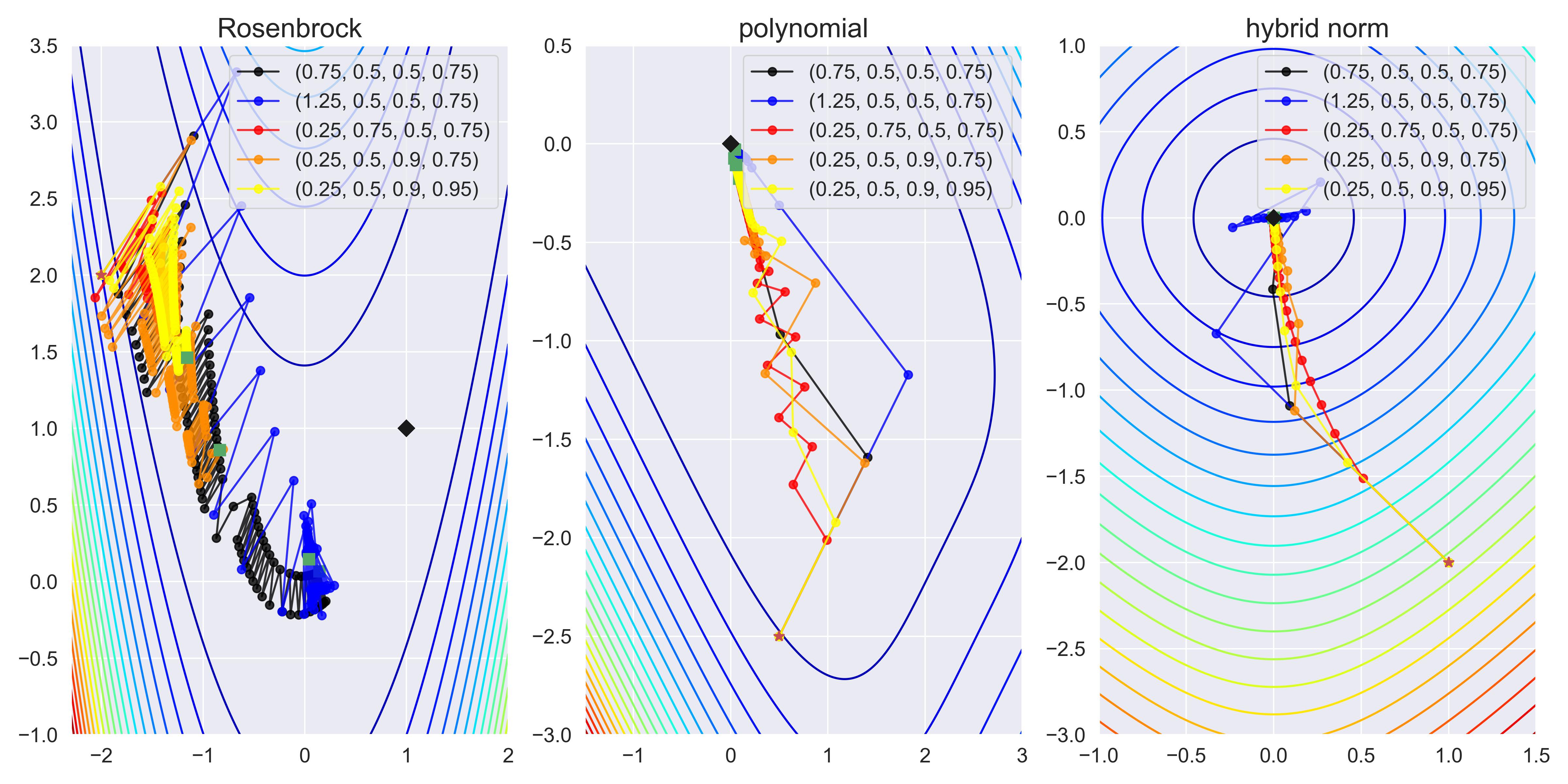}
    \end{subfigure}
\caption{Comparison of optimizers}
\label{fig:comparison_various_optimizers}
\end{figure}

In figure \eqref{fig:comparison_various_optimizers} we have chosen as benchmark objective functions the following mappings: the non-convex \emph{Rosenbrock} function endowed with a single global minimum at $(1, 1)^T$, the \emph{hybrid norm} function used in \cite{AlecsaOptimization}, defined as $f(x, y) = \sqrt{1 + x^2} + \sqrt{1 + y^2}$ (which has the minimum point $(0, 0)^T$ and is a convex function), and the \emph{polynomial} function employed in \cite{DaSilva} defined as $f(x, y) = (x + y)^4 + \left( \dfrac{x}{2} - \dfrac{y}{2} \right)^4$ which has the minimum point $(0, 0)^T$. Furthermore, for each of the aforementioned functions we have chosen the initial states $(-2, 2)^T$, $(1, -2)^T$ and $(0.5, -2.5)^T$, respectively. In our first two experiments shown in the first and second rows of figure \eqref{fig:comparison_various_optimizers} we have compared our \eqref{AdaptiveAcceleratedMomentumMethodShiftedUpdates} optimizer with the following adaptive optimization algorithms: \emph{Adam} from \cite{Adam}, \emph{Padam} from \cite{Padam}, \emph{AdaBelief} from \cite{AdaBelief}, \emph{Apollo} presented in \cite{Apollo}, \emph{Ranger} from \cite{Ranger} and Madgrad introduced in \cite{Madgrad}, respectively. For the experiment presented in the first and second rows we have chosen the number of iterations to be equal to $100$ and $400$, while the learning rates were chosen as $1\text{e-}1$ and $1\text{e-}3$. For all the adaptive methods, we chose $\varepsilon = 1\text{e-}8$ and $(\beta_1, \beta_2) = (0.9, 0.999)$ along with their default parameters. However, the only parameters which were selected by us in a different manner were: $(M, \mu, \nu, \tilde{\gamma}) = (0.75, 0.5, 0.5, 0.75)$ for \eqref{AdaptiveAcceleratedMomentumMethodShiftedUpdates}, \emph{momentum} term equal to $0.9$ for \emph{Madgrad} and \emph{amsgrad} option equal to \emph{True} for both \emph{Padam} and \emph{AdaBelief}. From the first row of figure \eqref{fig:comparison_various_optimizers} we observe that \eqref{AdaptiveAcceleratedMomentumMethodShiftedUpdates} presents an oscillating convergent behavior with respect to the non-convex landscape of the \emph{Rosenbrock} function, which is influenced by the high value of the learning rate. On the other hand, for the convex mappings presented in the last two columns the power of our algorithm can be seen due to the fact that \eqref{AdaptiveAcceleratedMomentumMethodShiftedUpdates} converges faster and in larger steps than the other chosen optimizers. This is in contrast especially with the \emph{Apollo} method which presents a circular-type trajectory near the minimum point. For the experiment presented in the first row of figure \eqref{fig:comparison_various_optimizers} we have set a low number of iterations and a high learning rate in order to observe the trajectory behavior of the optimizers under consideration. On the other hand, in the second row of figure \eqref{fig:comparison_various_optimizers} we have considered a larger number of iterations and a lower learning rate in order to see the full behavior of the adaptive methods. It can be easily observed that our \eqref{AdaptiveAcceleratedMomentumMethodShiftedUpdates} algorithm is competitive with the other adaptive optimizers (the only algorithm which stands out is the \emph{Madgrad} method for which its trajectories seem to be different that the ones of the other optimizers, in the case of convex functions). In our last row of figure \eqref{fig:comparison_various_optimizers} we consider comparing the variation of \eqref{AdaptiveAcceleratedMomentumMethodShiftedUpdates} algorithm with respect to the different choices of the parameters, where in the legend of the plots the underlying parameters are chosen in the order $(M, \mu, \nu, \tilde{\gamma})$. We highlight the following empirical observations of the consequences of different choices of the coefficients: when $M = 0.25$ or $\nu = 0.9$ the convergence behavior is very different for convex and non-convex mappings (this can also be the effect of the high value of the $\tilde{\gamma}$ parameter), while a big value of $M$, i.e. $1.25$ makes an oscillating behavior in the case of non-convex functions along with a very clear trajectory for convex functions. Moreover, although the behavior is oscillating when $(\mu, \nu, \tilde{\gamma}) = (0.5, 0.5, 0.75)$ the convergence is faster in the non-convex case of the \emph{Rosenbrock} function. \\
\vskip+0.15cm
We end this subsection with a step-by-step guide regarding how our algorithm is actually constructed:
\begin{itemize}
\item We extend the algorithms presented in \cite{GhadimiLanAccelerated} to the method \eqref{AcceleratedGradientMethod} endowed with two inertial steps.
\item We show that \eqref{AcceleratedGradientMethod} is equivalent to the \eqref{AdaptiveHeayBallMethod} which represents a generalization of accelerated gradient methods.
\item Inspired by \cite{BarakatBianchi} we introduce our method \eqref{AdaptiveAcceleratedMomentumMethodShiftedUpdates} which is equivalent to the aforementioned \eqref{AcceleratedGradientMethod} through the changes of variables given in \eqref{ChangesVariablesAlgorithms}.
\item Under the assumptions stated in section \eqref{Section_AssumptionSetting} we present theoretical convergence results for \eqref{AcceleratedGradientMethod} type methods.
\item In section \eqref{Section_ExamplesAdaptiveMethods} we construct an optimizer equivalent with and adaptive version of the \eqref{AdaptiveAcceleratedMomentumMethodShiftedUpdates} method, which can be regarded as a generalized version of \emph{AMSGrad}. 
\item By taking into account the equivalence between \eqref{AcceleratedGradientMethod}, \eqref{AdaptiveHeayBallMethod} and \eqref{AdaptiveAcceleratedMomentumMethodShiftedUpdates}, we present a practical formulation of our algorithm in the section \eqref{EquivalentForms}.
\item The algorithmic description of our adaptive method is presented in Algorithm \eqref{AlgDescription_Sutskever_AAMMSU} under the name \emph{Sutskever formulation of AAMMSU}. 
\item The associated \texttt{PyTorch} implementation can be found at \url{https://github.com/CDAlecsa/AAMMSU}.
\end{itemize}

\section{Adaptive algorithms seen as accelerated heavy ball type methods}\label{Section_AdaptiveHeavyBall}
In the present section, inspired by the results from \cite{GhadimiLanAccelerated} regarding accelerated first-order methods, the equivalence between adaptive Adam-like optimizers and accelerated Nesterov-type algorithms is given. For this, we firstly examine the following accelerated method, briefly called \emph{2 Steps Accelerated Gradient Method}, which represents an extension of the ones in \cite{GhadimiLanAccelerated} to two inertial steps ($\mu_n \in \mathbb{R}$ and $\tilde{\vb* \gamma}_n$ from $\mathbb{R}$ or possibly $\mathbb{R}^d$), namely for each $n \geq 1$ it takes the form
\begin{align}\label{AcceleratedGradientMethod}\tag{2SAGM}
\begin{cases}
\vb z_n = (1 - \tilde{\vb* \gamma}_n) \odot \vb* \theta_n + \tilde{\vb* \gamma}_n \odot \vb w_n \\ \\
\vb y_n = (1 - \mu_n) \vb* \theta_n + \mu_n \vb w_n \\ \\
\vb w_{n+1} = \vb w_n - \vb* \lambda_n \odot \vb g_n \\ \\
\vb* \theta_{n+1} = \vb y_n - \vb* \alpha_n \odot \vb g_n,
\end{cases}
\end{align}
where for every $n \geq 1$, $\tilde{\vb* \gamma}_n \in \mathbb{R}^d$ (or from $\mathbb{R}$ in the simplest case) and $\mu_n \in \mathbb{R}$ represent the positive coefficients of the inertial steps, while $\vb* \lambda_n$ and $\vb* \alpha_n$ represent two (possible non-constant) stepsizes which are also positive. It is important to note that, for a particular iteration $n$, the notations $\vb y_n$ and $\vb z_n$ underline the fact that these variables were computed previously, while $\vb* \theta_{n+1}$ and $\vb w_{n+1}$ actually represent that these are computed at the current iteration. Hence, for $n = 1$, we can let $\vb y_1$ and $\vb z_1$ to be initial values. But, in order to simplify the analysis and to not have three initial data $(\vb y_1, \vb z_1, \vb w_1)$ be given, we will take, in practical applications, only $(\vb* \theta_1, \vb w_1)$ to be the initial data, while $\vb y_1$ and $\vb z_1$ to be defined through the formula of \eqref{AcceleratedGradientMethod} (clearly, the initial adaptive stepsize values and the inertial coefficient values go together with all of these). Then, for $n = 1$, instead of $(\vb y_1, \vb z_1, \vb w_1, \vb* \alpha_1, \vb* \lambda_1)$, we have the given initial data $(\vb* \theta_1, \vb w_1, \vb* \alpha_1, \vb* \lambda_1)$ and also $(\tilde{\vb* \gamma}_1, \mu_1)$. At the same time, as we will see in Example \eqref{Example_1}, we take $\mu_1 = 1$ and $0 \leq \mu_n < 1$ for each $n \geq 2$, along with the assumption that $\dfrac{\| 1 - \tilde{\vb* \gamma}_n \|^2}{(1-\mu_n)}$ when $\tilde{\vb* \gamma}_n \in \mathbb{R}^d$ or $\dfrac{(1 - \tilde{\vb* \gamma}_n)^2}{(1-\mu_n)}$ when $\tilde{\vb* \gamma}_n \in \mathbb{R}$, are bounded (almost surely if these are adaptive, namely when they depend on some random vectors from \eqref{AcceleratedGradientMethod}). \\
Further, the random vector $\tilde{\vb* \gamma}_n$ will be specified later, but we mention that it must satisfy $\tilde{\vb* \gamma}_1 = 1$. On the other hand, we are considering $\vb g_n$ to be a noisy stochastic evaluation of the true gradient, defined at the shifted position $\vb z_n$. We will demonstrate how \eqref{AcceleratedGradientMethod} can be viewed as a straightforward accelerated gradient method with two inertial steps, similar to the algorithms shown in \cite{ALP}. To do this, we follow the steps presented below. From the fact that $\vb w_{n+1} = \vb w_n - \vb* \lambda_n \odot \vb g_n$ and that $\vb* \lambda_n \neq \textbf{0}$, one has, for every $n \geq 1$, that
\begin{align*}
- \vb g_n = \dfrac{\textbf{1}}{\vb* \lambda_n} \odot \left( \vb w_{n+1} - \vb w_n \right).    
\end{align*}
In light of the fact that $\vb* \theta_{n+1} = \vb y_n - \vb* \alpha_n \odot \vb g_n$, we obtain, for every $n \geq 1$, that
\begin{align*}
\vb* \theta_{n+1} = \vb y_n + \dfrac{\vb* \alpha_n}{\vb* \lambda_n} \odot \left( \vb w_{n+1} - \vb w_n \right).
\end{align*}
Now, because $\vb y_n = (1 - \mu_n) \vb* \theta_n + \mu_n \vb w_n$, for each $n \geq 1$ we get that
\begin{align*}
\vb* \theta_{n+1} = (1-\mu_n) \vb* \theta_n + \dfrac{\vb* \alpha_n}{\vb* \lambda_n} \odot \vb w_{n+1} + \left( \mu_n - \dfrac{\vb* \alpha_n}{\vb* \lambda_n} \right) \odot \vb w_n.   
\end{align*}
For $n \geq 1$, multiplying the above with $\dfrac{\vb* \lambda_n}{\vb* \alpha_n} \neq \textbf{0}$, we have
\begin{align*}
\vb w_{n+1} = \dfrac{\vb* \lambda_n}{\vb* \alpha_n} \odot \vb* \theta_{n+1} - (1 - \mu_n) \dfrac{\vb* \lambda_n}{\vb* \alpha_n} \odot \vb* \theta_n - \left( \mu_n \dfrac{\vb* \lambda_n}{\vb* \alpha_n} - 1 \right) \odot \vb w_n,   
\end{align*}
hence for every $n \geq 2$, it follows that
\begin{align}\label{eq:1}
\vb w_n = \dfrac{\vb* \lambda_{n-1}}{\vb* \alpha_{n-1}} \odot \vb* \theta_n - (1 - \mu_{n-1}) \dfrac{\vb* \lambda_{n-1}}{\vb* \alpha_{n-1}} \odot \vb* \theta_{n-1} - \left( \mu_{n-1} \dfrac{\vb* \lambda_{n-1}}{\vb* \alpha_{n-1}} - 1 \right) \odot \vb w_{n-1}.   
\end{align}
Inserting the above expression for $\vb w_n$ into $\vb y_n = (1 - \mu_n) \vb* \theta_n + \mu_n \vb w_n$, for $n \geq 2$ it follows that 
\begin{align}\label{eq:2}
\vb y_n = \vb* \theta_n + \mu_n \dfrac{\vb* \lambda_{n-1}}{\vb* \alpha_{n-1}} \odot \left[ \left( 1 - \dfrac{\vb* \alpha_{n-1}}{\vb* \lambda_{n-1}} \right) \odot \vb* \theta_n - (1 - \mu_{n-1}) \vb* \theta_{n-1} \right] - \mu_n \left( \mu_{n-1} \dfrac{\vb* \lambda_{n-1}}{\vb* \alpha_{n-1}} - 1 \right) \odot \vb w_{n-1}.    
\end{align}
Since $\vb w_n = \vb w_{n+1} + \vb* \lambda_n \odot \vb g_n$ for $n \geq 1$, hence $\vb w_{n-1} = \vb w_n + \vb* \lambda_{n-1} \odot \vb g_{n-1}$ for each $n \geq 2$, then by combining $\vb w_{n-1} = \vb w_n + \vb* \lambda_{n-1} \odot \vb g_{n-1}$ for $n \geq 2$ with the fact that $\vb w_n = \dfrac{1}{\mu_n} \vb y_n - \dfrac{1-\mu_n}{\mu_n} \vb* \theta_n$ for $n \geq 1$ (where $\mu_n \neq 0$) then, for each $n \geq 2$, we get that
\begin{align}\label{eq:3}
\vb w_{n-1} = \dfrac{1}{\mu_n} \vb y_n - \dfrac{1-\mu_n}{\mu_n} \vb* \theta_n + \vb* \lambda_{n-1} \odot \vb g_{n-1}.    
\end{align}
Then, \eqref{eq:2} and \eqref{eq:3} imply that for all $n \geq 2$
\begin{align}\label{eq:4}
\vb y_n &= \vb* \theta_n + \dfrac{\vb* \alpha_{n-1}}{\mu_{n-1} \vb* \lambda_{n-1}} \odot \left( - \dfrac{\mu_{n-1} \vb* \lambda_{n-1}}{\vb* \alpha_{n-1}} + 1 + \mu_n \dfrac{\vb* \lambda_{n-1}}{\vb* \alpha_{n-1}} \odot \left( 1 - \dfrac{\vb* \alpha_{n-1}}{\vb* \lambda_{n-1}} \right) + (1 - \mu_n) \dfrac{\vb* \lambda_{n-1}}{\vb* \alpha_{n-1}} \odot \left( \mu_{n-1} - \dfrac{\vb* \alpha_{n-1}}{\vb* \lambda_{n-1}} \right)  \right) \odot \vb* \theta_n \nonumber \\
&- \dfrac{\mu_n}{\mu_{n-1}} (1 - \mu_{n-1}) \vb* \theta_{n-1} - \mu_n \left( \vb* \lambda_{n-1} - \dfrac{\vb* \alpha_{n-1}}{\mu_{n-1}} \right) \odot \vb* g_{n-1}.   
\end{align}
However, it is clear from some straightforward calculations that
\begin{align*}
- \dfrac{\mu_{n-1} \vb* \lambda_{n-1}}{\vb* \alpha_{n-1}} + 1 + \mu_n \dfrac{\vb* \lambda_{n-1}}{\vb* \alpha_{n-1}} \odot \left( 1 - \dfrac{\vb* \alpha_{n-1}}{\vb* \lambda_{n-1}} \right) + (1 - \mu_n) \dfrac{\vb* \lambda_{n-1}}{\vb* \alpha_{n-1}} \odot \left( \mu_{n-1} - \dfrac{\vb* \alpha_{n-1}}{\vb* \lambda_{n-1}} \right) = \mu_n (1 - \mu_{n-1}) \dfrac{\vb* \lambda_{n-1}}{\vb* \alpha_{n-1}},
\end{align*}
therefore, for every $n \geq 2$, \eqref{eq:4} becomes
\begin{align*}
\vb y_n = \vb* \theta_n + \dfrac{\mu_n}{\mu_{n-1}} (1 - \mu_{n-1}) (\vb* \theta_n - \vb* \theta_{n-1}) - \dfrac{\mu_n}{\mu_{n-1}} (\mu_{n-1} \vb* \lambda_{n-1} - \vb* \alpha_{n-1}) \odot \vb* g_{n-1},    
\end{align*}
where the coefficient of $\vb g_{n-1}$ can be considered as the non-\emph{NAG} term. \\
We now follow the same procedure for the second inertial term $\vb z_n$. From the fact that $\vb z_n = (1 - \tilde{\vb* \gamma}_n) \odot \vb* \theta_n + \tilde{\vb* \gamma}_n \odot \vb w_n$ for $n \geq 1$ and from \eqref{eq:1} which is stated only from $n \geq 2$, we obtain the following computations for every $n \geq 2$:
\begin{align*}
\vb z_n = \vb* \theta_n + \tilde{\vb* \gamma}_n \odot \dfrac{\vb* \lambda_{n-1}}{\vb* \alpha_{n-1}} \odot \left( \left( 1 - \dfrac{\vb* \alpha_{n-1}}{\vb* \lambda_{n-1}} \right) \odot \vb* \theta_n - (1 - \mu_{n-1}) \vb* \theta_{n-1} \right) - \tilde{\vb* \gamma}_n \odot \left( \dfrac{\mu_{n-1} \vb* \lambda_{n-1}}{\vb* \alpha_{n-1}} - 1 \right) \odot \vb w_{n-1}.    
\end{align*}
For each $n \geq 2$, utilizing that $\vb w_{n-1} = \dfrac{1}{\tilde{\vb* \gamma}_n} \odot \vb z_n - \dfrac{1 - \tilde{\vb* \gamma}_n}{\tilde{\vb* \gamma}_n} \odot \vb* \theta_n + \vb* \lambda_{n-1} \odot \vb g_{n-1}$, it follows that
\begin{align}\label{eq:5}
\vb z_n &= \vb* \theta_n + \dfrac{\vb* \alpha_{n-1}}{\mu_{n-1} \vb* \lambda_{n-1}} \odot \left( 1 - \dfrac{\mu_{n-1} \vb* \lambda_{n-1}}{\vb* \alpha_{n-1}} + \tilde{\vb* \gamma}_n \odot \dfrac{\vb* \lambda_{n-1} - \vb* \alpha_{n-1}}{\vb* \alpha_{n-1}} + (1 - \tilde{\vb* \gamma}_n) \odot \dfrac{\mu_{n-1} \vb* \lambda_{n-1} - \vb* \alpha_{n-1}}{\vb* \alpha_{n-1}} \right) \odot \vb* \theta_n \nonumber \\
&- \dfrac{\tilde{\vb* \gamma}_n}{\mu_{n-1}} \odot (1 - \mu_{n-1}) \vb* \theta_{n-1} - \dfrac{\tilde{\vb* \gamma}_n}{\mu_{n-1}} \odot (\mu_{n-1} \vb* \lambda_{n-1} - \vb* \alpha_{n-1}) \odot \vb g_{n-1}.
\end{align}
For the coefficient of $\vb* \theta_n$ from the right hand side, we obtain that
\begin{align}\label{eq:6}
\dfrac{\vb* \alpha_{n-1}}{\mu_{n-1} \vb* \lambda_{n-1}} \odot \left( - \dfrac{\mu_{n-1} \vb* \lambda_{n-1}}{\vb* \alpha_{n-1}} + 1 + \tilde{\vb* \gamma}_n \odot \dfrac{\vb* \lambda_{n-1} - \vb* \alpha_{n-1}}{\vb* \alpha_{n-1}} + (1 - \tilde{\vb* \gamma}_n) \odot \dfrac{\mu_{n-1} \vb* \lambda_{n-1} - \vb* \alpha_{n-1}}{\vb* \alpha_{n-1}} \right) = \dfrac{\tilde{\vb* \gamma}_n}{\mu_{n-1}} (1 - \mu_{n-1}).
\end{align}
When \eqref{eq:5} and \eqref{eq:6} are combined, it means that for every $n \geq 2$
\begin{align*}
\vb z_n = \vb* \theta_n + \dfrac{\tilde{\vb* \gamma}_n}{\mu_{n-1}} \odot (1 - \mu_{n-1}) (\vb* \theta_n - \vb* \theta_{n-1}) - \dfrac{\tilde{\vb* \gamma}_n}{\mu_{n-1}} \odot (\mu_{n-1} \vb* \lambda_{n-1} - \vb* \alpha_{n-1}) \odot \vb g_{n-1}.
\end{align*}
All the calculations that have been made so far confirm our conclusion that \eqref{AcceleratedGradientMethod} is actually equivalent to the following \emph{Adaptive Heavy Ball Method} from $n \geq 2$, namely
\begin{align}\label{AdaptiveHeayBallMethod}\tag{AHBM}
\begin{cases}
\vb y_n = \vb* \theta_n + \beta_n (\vb* \theta_n - \vb* \theta_{n-1}) - \mu_n \vb* \omega_n \odot \vb* g_{n-1} \\ \\
\vb z_n = \vb* \theta_n + \vb* \gamma_n \odot (\vb* \theta_n - \vb* \theta_{n-1}) - \tilde{\vb* \gamma}_n \odot \vb* \omega_n \odot \vb* g_{n-1} \\ \\
\vb* \theta_{n+1} = \vb y_n - \vb* \alpha_n \odot \vb* g_n,
\end{cases}
\end{align}
where $\beta_n := \dfrac{\mu_n}{\mu_{n-1}}(1 - \mu_{n-1}) \in \mathbb{R}$, $\vb* \gamma_n := \dfrac{\tilde{\vb* \gamma}_n}{\mu_{n-1}}(1 - \mu_{n-1}) \in \mathbb{R}^d$ (or $\mathbb{R}$ in the simplest case which we will employ in the last section) and 
\begin{align}\label{omega_n}
    \vb* \omega_n := \vb* \lambda_{n-1} - \dfrac{\vb* \alpha_{n-1}}{\mu_{n-1}} \in \mathbb{R}^d. 
\end{align}
In \eqref{AdaptiveHeayBallMethod}, the initial data are stated for $n = 2$, namely $(\vb y_2, \vb z_2, \vb* \alpha_2)$, but as in the case of \eqref{AcceleratedGradientMethod} we compute $\vb y_2$ and $\vb z_2$ through the above formulas, hence $(\vb* \theta_1, \vb* \theta_2, \vb z_1, \vb* \alpha_1, \vb* \alpha_2, \vb* \lambda_1, \mu_1, \mu_2, \tilde{\vb* \gamma}_2)$ are given, under the assumption that $\mu_1 = \tilde{\vb* \gamma}_1 = 1$. Evidently, $(\mu_1, \vb* \alpha_1, \vb* \lambda_1, \vb* \theta_1)$ are also initial conditions for \eqref{AcceleratedGradientMethod}, while $(\vb* \theta_2, \vb z_1)$ can be computed explicitly through the iterative process given by \eqref{AcceleratedGradientMethod}. 
Similar to \eqref{AcceleratedGradientMethod}, $\vb y_n$ and $\vb z_n$ are considered computed in the iteration before $\vb* \theta_{n+1}$ (this is exactly the reason why we have denoted the former ones with a shifted index). \\
Now we shall focus on Adam-type algorithms. Inspired by \cite{BarakatBianchi}, we consider the following momentum-like method and we will eventually show that our adaptive algorithm is actually equivalent to \eqref{AdaptiveHeayBallMethod}. Our iterative process defined for $n \geq 2$, which will be named \emph{Adaptive Accelerated Momentum Method with Shifted Updates} reads as follows.
\begin{align}\label{AdaptiveAcceleratedMomentumMethodShiftedUpdates}\tag{AAMMSU}
\begin{cases}
\vb p_{n+1} = \vb* \tau_n \odot \vb p_n + \vb* \nu_n \odot \vb g_n - \vb k_n \odot \vb g_{n-1} \\
\vb* \theta_{n+1} = \vb* \theta_n - \vb r_n \odot \vb p_{n+1},
\end{cases}    
\end{align}
with $\vb g_n$ being the stochastic gradient estimate in the shifted iterate $\vb z_n$ (which will be defined later). At this moment, we are not interested in the initial data for \eqref{AdaptiveAcceleratedMomentumMethodShiftedUpdates}, since we will actually show that \eqref{AdaptiveAcceleratedMomentumMethodShiftedUpdates} is equivalent to \eqref{AdaptiveHeayBallMethod} for some choices of coefficients starting from $n \geq 3$. Now, from the algorithm's definition, for every $n \geq 2$, we have that
\begin{align*}
\vb* \theta_{n+1} = \left( \vb* \theta_n - \vb r_n \odot \vb* \tau_n \odot \vb p_n + \vb r_n \odot \vb k_n \odot \vb g_{n-1} \right) - \vb r_n \odot \vb* \nu_n \odot \vb g_n.    
\end{align*}
Using that the stepsize is $\vb* \alpha_n$, we firstly deduce that $\vb* \alpha_n = \vb r_n \odot \vb* \nu_n$. Second of all, from the definition of $\vb y_n$ from \eqref{AdaptiveHeayBallMethod}, we infer, for every $n \geq 2$, that
\begin{align*}
\vb* \theta_n + \beta_n (\vb* \theta_n - \vb* \theta_{n-1}) - \mu_n \vb* \omega_n \odot \vb g_{n-1} = \vb* \theta_n - \vb r_n \odot \vb* \tau_n \odot \vb p_n + \vb r_n \odot \vb k_n \odot \vb g_{n-1}.
\end{align*}
For $n \geq 2$, we easily obtain that 
\begin{align}\label{eq:7}
\vb k_n = - \mu_n \dfrac{\vb* \omega_n}{\vb r_n}.    
\end{align}
Moreover, we discover that the momentum term takes, for $n \geq 2$, the form 
\begin{align}\label{eq:8}
\vb p_n = - \dfrac{\beta_n}{\vb r_n \odot \vb* \tau_n} \odot (\vb* \theta_n - \vb* \theta_{n-1}).
\end{align}
Since, for $n \geq 2$, $\vb* \theta_{n+1} - \vb* \theta_n = - \vb r_n \odot \vb p_{n+1}$, then, for every $n \geq 3$, we get that 
\begin{align}\label{eq:9}
\vb p_n = - \dfrac{1}{\vb r_{n-1}} \odot (\vb* \theta_n - \vb* \theta_{n-1}).    
\end{align}
For $n \geq 3$, from \eqref{eq:8} and \eqref{eq:9}, it follows that
\begin{align}\label{eq:tau_n}
\vb* \tau_n = \beta_n \dfrac{\vb r_{n-1}}{\vb r_n}.    
\end{align}
We will now address the shifted update $\vb z_n$ as follows. Since $\vb z_n = \vb* \theta_n + \vb* \gamma_n \odot (\vb* \theta_n - \vb* \theta_{n-1}) - \tilde{\vb* \gamma}_n \odot \vb* \omega_n \odot \vb g_{n-1}$ for $n \geq 2$, and using the above computations, one gets, for all $n \geq 3$, that
\begin{align}\label{zn_AAMMSU}
\vb z_n = \vb* \theta_n - \vb* \gamma_n \odot \vb r_{n-1} \odot \vb p_n - \tilde{\vb* \gamma}_n \odot \vb* \omega_n \odot \vb g_{n-1}.    
\end{align}
The final computation we have to deal with is to find the coefficients $\tilde{\gamma}_n$, $\mu_n$, $\vb* \lambda_n$ and $\vb* \alpha_n$ in terms of the given data $\vb r_n$, $\vb* \tau_n$, $\vb* \nu_n$, $\vb k_n$ and $\vb* \gamma_n$, respectively. This is of utmost importance, since our convergence results are proved for the algorithms which take the form \eqref{AcceleratedGradientMethod} as in \cite{GhadimiLanAccelerated}, and only after that we will find the corresponding adaptive Adam-type methods which are given in \eqref{AdaptiveAcceleratedMomentumMethodShiftedUpdates}. In order to do this, we will continue using the following strategy. \\
In what follows, we will take into account that, in Example \eqref{Example_1}, $\beta_2 = \vb* \gamma_2 = 0$, hence we avoid $\beta_n$ and $\vb* \gamma_n$ to be at the denominator in the computations from below. On the other hand, from the fact that $\vb* \gamma_n = \dfrac{\tilde{\vb* \gamma}_n}{\mu_{n-1}}(1-\mu_{n-1})$ and $\beta_n = \dfrac{\mu_n}{\mu_{n-1}}(1-\mu_{n-1})$ for $n \geq 2$, one gets that $1 - \mu_{n-1} = \dfrac{\vb* \gamma_n}{\tilde{\vb* \gamma}_n} \mu_{n-1}$ and $1 - \mu_{n-1} = \dfrac{\beta_n}{\mu_n} \mu_{n-1}$, where we have used that $\tilde{\vb* \gamma}_n, \, \mu_n \neq 0$ for each $n \geq 2$, respectively. Hence, since $\mu_{n-1} \neq 0$, it follows that
\begin{align}\label{eq:10}
\beta_n = \dfrac{\mu_n \vb* \gamma_n}{\tilde{\vb* \gamma}_n},
\end{align}
where $\mu_n$ must be viewed as $\mu_n \cdot \textbf{1}$. Also, since $\vb* \omega_n = \vb* \lambda_{n-1} - \dfrac{\vb* \alpha_{n-1}}{\mu_{n-1}}$, for $n \geq 2$, we have that
\begin{align*}
\vb* \omega_n = \vb* \lambda_{n-1} - \dfrac{\vb* \gamma_{n-1}}{\tilde{\vb* \gamma}_{n-1}} \odot \dfrac{\vb* \alpha_{n-1}}{\beta_{n-1}}.
\end{align*}
Because $\vb* \tau_n = \beta_n \dfrac{\vb r_{n-1}}{\vb r_n}$ from \eqref{eq:tau_n}, we find that $\beta_n = \vb* \tau_n \odot \dfrac{\vb r_n}{\vb r_{n-1}}$ for $n \geq 3$, where $\beta_n$ must be interpreted as $\beta_n \cdot \textbf{1}$. When combined this with \eqref{eq:10}, it follows, also for each $n \geq 3$, that
\begin{align}\label{eq:11}
\mu_n \dfrac{\vb* \gamma_n}{\tilde{\vb* \gamma}_n} = \dfrac{\vb r_n}{\vb r_{n-1}} \odot \vb* \tau_n.
\end{align}
Nevertheless, in light of \eqref{eq:11} and that $\dfrac{\vb* \gamma_n}{\tilde{\vb* \gamma}_n} = \dfrac{1 - \mu_{n-1}}{\mu_{n-1}}$ for every $n \geq 2$, it implies that the sequence with the general term $\mu_n$ satisfies the implicit recursion for every $n \geq 3$, namely
\begin{align*}
\dfrac{1 - \mu_{n-1}}{\mu_{n-1}} \mu_n = \dfrac{\vb r_n}{\vb r_{n-1}} \odot \vb* \tau_n.
\end{align*}
Combining \eqref{eq:7} with the fact that $\vb* \omega_n = \vb* \lambda_{n-1} - \dfrac{\vb* \alpha_{n-1}}{\mu_{n-1}}$ for $n \geq 2$, it leads to
\begin{align*}
\vb* \lambda_{n-1} = - \dfrac{\vb k_n \odot \vb r_n}{\mu_n} + \dfrac{\vb* \alpha_{n-1}}{\mu_{n-1}}.  
\end{align*}
Following all of these computations, we find the following changes in variables that link the adaptive momentum algorithm \eqref{AdaptiveAcceleratedMomentumMethodShiftedUpdates} with the accelerated gradient method that contains two inertial steps \eqref{AcceleratedGradientMethod}, namely:
\begin{align}\label{ChangesVariablesAlgorithms}\tag{VarChanges}
\begin{cases}
\vb* \alpha_n = \vb r_n \odot \vb* \nu_n \text{ for every } n \geq 2 \\ \\
\beta_n = \dfrac{\mu_n \vb* \gamma_n}{\tilde{\vb* \gamma}_n} \text{ for every } n \geq 2 \\ \\
\vb* \gamma_n = \dfrac{\tilde{\vb* \gamma}_n}{\mu_{n-1}} (1-\mu_{n-1}) \text{ for every } n \geq 2 \\ \\
\vb* \lambda_{n-1} = - \dfrac{\vb k_n \odot \vb r_n}{\mu_n} + \dfrac{\vb* \alpha_{n-1}}{\mu_{n-1}} \text{ for every } n \geq 2 \\ \\
\vb* \tau_n = \beta_n \dfrac{\vb r_{n-1}}{\vb r_n} \text{ for every } n \geq 3.
\end{cases}
\end{align}

\section{The assumptions setting for our algorithms}\label{Section_AssumptionSetting}
In this section, we will look at the key assumptions that we will apply in the following sequel to get our main results. We shall stick to the hypotheses that are frequently employed in Machine Learning (for these, we especially refer to \cite{ChenAdam}, \cite{DefossezAdam} and \cite{GadatGavra}). Furthermore, the setting which we employ focuses on the section \eqref{Section_Preliminaries}, concerning the notation with respect to stochastic gradients. Our first assumption refers to the boundedness of the objective function $F$, similar to the work \cite{ALV}, namely
\begin{center}
\fbox{$(\vb{A_1}) \quad F : \mathbb{R}^d \to \mathbb{R}$ is bounded from below.}    
\end{center}
Our second assumption refers to the classical globally Lipschitz continuity of the gradient of the objective function. This represents a standard tool in the Optimization community and it was used in \cite{BarakatBianchi}, \cite{Bottou}, \cite{DefossezAdam}, \cite{GadatGavra}, \cite{ZouAdam} and references therein.
\begin{center}
\fbox{$(\vb{A_2}) \quad \exists L > 0, \, \text{ such that } \forall (x, y) \in \mathbb{R}^d \times \mathbb{R}^d: \, \| \nabla F(x) - \nabla F(y) \| \leq L \| x - y \|$.}    
\end{center}
Following \cite{ALV}, \cite{GadatGavra} and \cite{GhadimiLanAccelerated}, it is worth mentioning that the Lipschitz assumption implies the so-called descent inequality which is the backbone behind our main result, i.e.
\begin{align}\label{DL_Inequality}\tag{DL}
\forall (x, y) \in \mathbb{R}^d \times \mathbb{R}^d: F(y) \leq F(x) + \langle \nabla F(x), y - x \rangle + \dfrac{L}{2} \| y - x \|^2.
\end{align}
In what follows, we consider the noisy evaluation of the true gradient, namely $\vb g_n = \vb* \delta_n + \nabla F(\vb z_n)$, where $\vb* \delta_n = \sigma \zeta_n$ being the noisy term, where $\sigma > 0$. Furthermore, we consider $(\mathcal{F}_n)_{n \geq 1}$ to be an increasing sequence of $\sigma-$fields. \\
Similarly to the stochastic setting from \cite{GadatGavra}, we impose also as a hypothesis the property of the stochastic gradient to be an unbiased estimate of the true gradient of the loss $F$ at the shifted position $\vb z_n$ (and not in $\vb* \theta_n$ as in the most of the related works), namely
\begin{center}
\fbox{$(\vb{A_3})$  For $n \geq 1$, we have $\mathbb{E}[\vb g_n \, | \, \mathcal{F}_n] = \nabla F(\vb z_n)$.} 
\end{center}
The next assumption refers to classical bounded variance assumption with respect to the stochastic gradient, i.e.
\begin{center}
\fbox{$(\vb{A_4})$ There exists $\sigma > 0$, such that for any $n \geq 1$, we have $\mathbb{E}[\| \vb* \delta_n \|^2 \, | \, \mathcal{F}_n] \leq \sigma^2, \, (w.p.1)$.} 
\end{center}
If we suppose the assumption (as in \cite{GadatGavra}) that $\mathbb{E}[\| \zeta_n \|^2 \, | \, \mathcal{F}_n] \leq 1 , \, (w.p.1)$, then we easily find that 
\begin{align*}
\mathbb{E}[\| \vb* \delta_n \|^2 \, | \, \mathcal{F}_n]  = \sigma^2 \mathbb{E}[\| \zeta_n \|^2 \, | \, \mathcal{F}_n] \leq \sigma^2 , \, (w.p.1),
\end{align*}
which means that $\sigma$ from \emph{($A_4$)} is the same as the one from the noisy gradient error $\vb* \delta_n = \sigma \zeta_n$. Now, for every $n \geq 1$, we will consider the following assumption regarding the measurability of some random vectors, namely
\begin{center}
\fbox{$(\vb{A_5})$ For any $n \geq 1$, the random variables $\vb w_n$, $\vb z_n$ $\vb* \lambda_n$, $\vb* \alpha_n$ are $\mathcal{F}_n$ measurable.}
\end{center}
For the adaptive cases (see our Example \eqref{Example_1}), $\vb* \lambda_n$ and $\vb* \alpha_n$ may depend on $\vb z_n$, which implies the measurability of the former random vectors, due to the fact that $\vb z_n$ is $\mathcal{F}_n$ measurable. But, for non-adaptive cases and from the above assumption, if $\vb* \theta_n \in \mathcal{F}_n$, we observe that $\vb z_n$ and $\vb y_n$ are also $\mathcal{F}_n$ measurable random vectors, since in \eqref{AcceleratedGradientMethod} $\vb z_n$ and $\vb y_n$ are convex combinations of $\vb* \theta_n$ and $\vb w_n$. So, since $\vb z_n$ is a convex combination of $\vb* \theta_n$ and $\vb w_n$, then we can require instead in \emph{($A_5$)} that $\vb* \theta_n \in \mathcal{F}_n$. Conversely, for $\vb w_n, \vb z_n \in \mathcal{F}_n$, we obtain that $\vb* \theta_n \in \mathcal{F}_n$. Also, we mention that in Example \eqref{Example_1} which we will be presenting, $\vb* \alpha_n = \nu \dfrac{\eta_n}{\varepsilon + \sqrt{\vb v_n}}$ for some $\eta_n$ and $\nu$ from $\mathbb{R}$, and where $\vb v_n$, which depends on $[\vb g_{n}]^2$, is the moving average of past squared gradients. On the other hand, in the last section of numerical experiments \eqref{Section_NeuralNetworks} we will take $\vb* \lambda_n = M \vb* \alpha_n$ for some positive term $M$ and $\vb* \gamma_n \in \mathbb{R}$ a deterministic constant, then assumption \emph{$(A_5)$} becomes $\vb w_n$, $\vb z_n \in \mathcal{F}_n$ for every $n \geq 1$. This makes clear that $\mathcal{F}_n$ represents the history of the algorithm \eqref{AcceleratedGradientMethod} up to time $n$, just before the stochastic gradient $\vb g_n$ is generated. Additionally, we mention that the assumption \emph{$(A_5)$} reveals that $\vb w_n, \vb z_n$, $\vb* \lambda_n$ and $\vb* \alpha_n$ are constructed in the previous iteration, before $\vb w_{n+1}$ and $\vb* \theta_{n+1}$ are defined, as in the method \eqref{AcceleratedGradientMethod}. Likewise, conditioning on $\mathcal{F}_n$ is equivalent to conditioning on the random vectors $\vb w_1, \ldots, \vb w_n, \vb* \theta_1, \ldots, \vb* \theta_n, \vb z_1, \ldots, \vb z_n, \vb y_1, \ldots, \vb y_n, \vb* \alpha_1, \ldots, \vb* \alpha_n$, and also $\vb* \lambda_1, \ldots, \vb* \lambda_n$ in the more general case (and also $\tilde{\vb* \gamma}_1, \ldots, \tilde{\vb* \gamma}_n$ when these are non-deterministic and from $\mathbb{R}^d$). Since, in \eqref{AcceleratedGradientMethod}, $\vb* \theta_{n+1}$ and $\vb w_{n+1}$ are generated at the current iteration, while $\vb z_n$ and $\vb y_n$ are generated in the previous iteration, then $\mathcal{F}_n$ could be considered the canonical filtration associated to the sequences of the iterates of \eqref{AcceleratedGradientMethod} (as in \cite{GadatGavra}), namely $\mathcal{F}_n = \sigma \left( (\vb* \theta_k)_{1 \leq k \leq n}, (\vb w_k)_{1 \leq k \leq n}, (\vb z_k)_{1 \leq k \leq n}, (\vb y_k)_{1 \leq k \leq n}, (\tilde{\vb* \gamma}_k)_{1 \leq k \leq n}, (\vb* \alpha_k)_{1 \leq k \leq n}, (\vb* \lambda_k)_{1 \leq k \leq n} \right)$.
Evidently, in the canonical filtration $\mathcal{F}_n$ we can get rid of $(\vb* \alpha_k)_{1 \leq k \leq n}$ and $(\vb* \lambda_k)_{1 \leq k \leq n}$ if we take into consideration the definition of these adaptive stepsizes as in Example \eqref{Example_1} and section \eqref{Section_NeuralNetworks}. \\
The final concept we discuss is that assumptions \emph{$(A_3)$}, \emph{$(A_4)$} and \emph{$(A_5)$} can be also applied to the i.i.d. setting as in \cite{StochasticQuasiNewton}. By considering $\xi_n$ to be independent samples, where each $\xi_n$ is independent of $(\vb z_j)_{1 \leq j \leq n}$, then \emph{$(A_3)$} can be replaced with $\mathbb{E}_{\xi_n} \left[ g(\vb z_n, \xi_n) \right] = \nabla F(\vb z_n)$, hence $\mathbb{E}[\vb g_n \, | \, \vb z_n] = \nabla F(\vb z_n)$. Moreover, assumption \emph{$(A_4)$} can be replaced with $\mathbb{E}_{\xi_n} \left[ \| g(\vb z_n, \xi_n) - \nabla F(\vb z_n) \|^2 \right] \leq \sigma^2$. Here, the notation $\mathbb{E}_{\xi_n}[\cdot]$ represents the expectation taken with respect to the distribution of $\xi_n$ which represents the random sampling in the $n$ iteration. Furthermore, proceeding similar to \cite{StochasticQuasiNewton}, we can take $\vb* \alpha_n$ and $\vb* \lambda_n$ to be dependent only on the random samplings from the first $n$ iterations, which is a hypothesis that is checked by Example \eqref{Example_1} and in our last section \eqref{Section_NeuralNetworks}.

Regarding the scenarios where the concept of bounded stochastic gradients is present (such as in \cite{ChenAdam}, \cite{DefossezAdam}, \cite{ReddiAdamConvergence}, \cite{ZhouAMSGrad} and \cite{ZouAdam}) we end the present section by considering our last assumption which will be necessary for the case of our \emph{AMSGrad}-like methods, i.e.
\begin{center}
\fbox{$(\vb{A_6})$  There exists $\mathcal{K} > 0$, such that for every $n \geq 1$, one has that $\| \vb g_n \| \leq \mathcal{K}, \, (w.p.1).$} 
\end{center}

\section{Convergence guarantees}\label{Section_ConvergenceResults}
In this section we present our convergence results regarding our accelerated gradient algorithm \eqref{AcceleratedGradientMethod}. 
We also analyze the worst-case iteration complexity of our adaptive optimizers. Furthermore, our qualitative results for \eqref{AcceleratedGradientMethod} are inspired by \cite{GhadimiLanAccelerated} and \cite{StochasticQuasiNewton}, while the convergence results for \eqref{AdaptiveAcceleratedMomentumMethodShiftedUpdates} is inferred through the equivalence of both algorithms as done in \cite{BarakatBianchi}. We mention that the proofs of all our results that we present are postponed to the Appendix section \eqref{Section_Appendix}. \\
The first result we show is a technical lemma regarding the almost sure inequality between the norm of a squared vector (in the componentwise sense) and the squared norm of the given vector.
\begin{proposition}\label{Proposition_SquaredVector}
Let $\vb u \in\mathbb{R}^d$ be an arbitrarily given vector. Then, it follows that
\begin{align*}
\| [\vb u]^2 \| \leq \| \vb u \|^2.    
\end{align*}
\end{proposition}

Our second result refers to a inequality involving the entrywise Hadamard product between two vectors.
\begin{proposition}\label{Proposition_HadamardProduct}
Let $\vb u, \vb v \in \mathbb{R}^d$ be two vectors chosen arbitrarily. Then, we have that
\begin{align*}
\| \vb u \odot \vb v \| \leq \| \vb u \| \cdot \| \vb v \|.    
\end{align*}
\end{proposition}

For completeness, we will establish the following simple algebraic procedure for a finite summation, as used in \cite{GhadimiLanAccelerated}.
\begin{proposition}\label{Proposition_SwitchSummation}
Let $N \geq n$, where $n \geq 1$ is a given natural number. Also, consider $e_1, \ldots, e_N$ and $f_1, \ldots, f_N$ to be real numbers. Then
\begin{align*}
\sum\limits_{n = 1}^{N} e_n \sum\limits_{j = 1}^{n} f_j = \sum\limits_{n = 1}^{N} \left( \sum\limits_{j = n}^{N} e_j \right) f_n.    
\end{align*}
\end{proposition}

We now provide the central pillar of our theoretical analysis, which is an almost surely descent-type property in relation to the underlying random variables.
\begin{lemma}\label{Lemma_AlmostSureConvergence}
We consider the accelerated gradient method \eqref{AcceleratedGradientMethod} given by $\vb* \theta_{n+1}$ with the shifted updates $\vb z_n$ starting from $n \geq 1$, where the auxiliary sequence $(\mu_n)_{n \geq 1}$ is chosen such that $\mu_1 = 1$ and $\mu_n \in (0, 1)$ for each $n \geq 2$. Further, suppose that the adaptive sequence $(\tilde{\vb* \gamma}_n)_{n\geq 1}$ defined by $\tilde{\vb* \gamma}_n$ is bounded, and therefore denote $B$ as a major bound for $\dfrac{ \|1 - \tilde{\vb* \gamma}_n \|^2}{(1 - \mu_n)}$ for $n \geq 2$. Also, suppose that $\tilde{\vb* \gamma}_1 = \textbf{1}$ and that the assumptions from the previous section \emph{$(A_1)$}, \emph{$(A_2)$}, \emph{$(A_3)$}, \emph{$(A_4)$} and \emph{$(A_5)$} are satisfied. By choosing $N \geq n$ (where $n \geq 1$) and taking into consideration $\Gamma_n$ defined in \eqref{eq:Gamma}, we denote 
\begin{align*}
\begin{cases}
C_{n, N} := \sum\limits_{j = n}^{N} \Gamma_j, \\
\vb Q_n := \langle \vb* \delta_n, \vb* \lambda_n \odot \left( \nabla F(\vb w_n) - L \vb* \lambda_n \odot \nabla F(\vb z_n) \right) - LB \dfrac{C_{n,N}}{\mu_n \Gamma_n} [\vb* \lambda_n - \vb* \alpha_n]^2 \odot \nabla F(\vb z_n) \rangle \\
\vb D_n := \vb* \lambda_n - L [\vb* \lambda_n]^2 - \dfrac{LB}{2} \dfrac{C_{n,N}}{\mu_n \Gamma_n} [\vb* \lambda_n - \vb* \alpha_n]^2.
\end{cases}
\end{align*}
In addition, assume that for all $n \geq 1$, there exists $m_n > 0$, such that $\vb D_n \geq m_n$ almost surely and in the Hadamard sense. 
Then, the following almost sure inequality holds with respect to the algorithm's random variables:
\footnotesize{
\begin{align}\label{eq:Results_AlmostSureInequality}
F(\vb w_{N+1}) - F(\vb w_1) + \sum\limits_{n = 1}^{N} m_n \cdot \| \nabla F(\vb z_n) \|^2 &\leq - \sum\limits_{n = 1}^{N} \vb Q_n + \dfrac{L}{2} \sum\limits_{n = 1}^{N} \| \vb* \delta_n \|^2 \cdot \left( \| \vb* \lambda_n \|^2 + B \dfrac{C_{n,N}}{\mu_n \Gamma_n} \| \vb* \lambda_n - \vb* \alpha_n \|^2 \right) , \, (w.p.1).
\end{align}
}
\end{lemma}

The following remark we present is related to the term $\tilde{\vb* \gamma}_n$.
\begin{remark}\label{Remark_NormBound}
It is easy to observe that, in the proof of Lemma \eqref{Lemma_AlmostSureConvergence}, if $\vb* \gamma_n$ is from $\mathbb{R}$, then the proof can be adjusted by taking the absolute value instead of the Euclidean norm of the random vector. More precisely, when $\tilde{\vb* \gamma}_n \in \mathbb{R}$ (possibly deterministic), we can utilize that $\dfrac{1}{2} \| (1 - \tilde{\vb* \gamma}_n)(\vb* \theta_n - \vb w_n) \|^2 = \dfrac{(1 - \tilde{\vb* \gamma}_n)^2}{2} \cdot \| \vb* \theta_n - \vb w_n \|^2$ in \eqref{eq:16}, hence the assumption from Lemma \eqref{Lemma_AlmostSureConvergence} must be $\dfrac{(1 - \gamma_n)^2}{(1-\mu_n)} \leq B$. When $\tilde{\vb* \gamma}_n \in \mathbb{R}$, this is in contrast with $\| 1 - \tilde{\vb* \gamma}_n \|^2$ which becomes $d (1 - \tilde{\vb* \gamma}_n)^2$ if we follow step by step the usage of the Euclidean norm from the aforementioned proof.
\end{remark}

The first conclusion of the result from below is related to the convergence rate with respect to a final iteration $N$,while the second conclusion is concerned with an asymptotic almost sure inequality. 
\begin{theorem}\label{Theorem_AsymptoticBehavior}
Consider the hypotheses along with the notations from the statement of Lemma \eqref{Lemma_AlmostSureConvergence}. Furthermore, suppose that $(\vb* \lambda_n)_{n \geq 1}$ and $(\vb* \alpha_n)_{n \geq 1}$ are almost surely bounded. Also, assume that the real valued sequence $(\Gamma_n)_{n \geq 1}$ is bounded, and that a (strictly positive) major bound can be found for $C_{n,N}$ (where $n$ and $N$ were defined already in Lemma \eqref{Lemma_AlmostSureConvergence}). Thus, for every $n \geq 1$ such that $n \leq N$, there exists $R_n >0$ that may depend on $N$ (but we did not emphasized the use of index $N$ in order to not have a cumbersome notation, but it can depend on $N$ hence also $C_{n,N}$ can have a major bound depending on $N$), such that $\| \vb* \lambda_n \|^2 + B \dfrac{C_{n,N}}{\mu_n \Gamma_n} \| \vb* \lambda_n - \vb* \alpha_n \|^2 \leq R_n$ almost surely, with respect to the entrywise Hadamard notations. Then, for the final iteration $N$, we have the convergence rate in the finite time horizon
\begin{align}\label{MinimumGradientRate}
\mathbb{E} \left[ \min_{n = 1, \ldots, N} \| \nabla F(\vb z_n) \|^2 \right] &\leq \dfrac{\overline{M} + \dfrac{L}{2} \sigma^2 \sum\limits_{n = 1}^{N} R_n}{\sum\limits_{n = 1}^{N} m_n}.
\end{align}
Let us suppose, on the other hand,
\begin{align}\label{Limit_Infinity_Gradient}
\sigma^2 \sum\limits_{n = 1}^{+\infty} R_n < +\infty,
\end{align}
where $\sum\limits_{n = 1}^{+\infty} R_n$ is a shorthand notation for $\lim\limits_{N \to +\infty} \sum\limits_{n = 1}^{N} R_n = \lim\limits_{N \to +\infty} \sum\limits_{n = 1}^{N} R_{n,N}$ where, as we have said before, $R_{n, N} := R_n$ can actually depend on N. Then, we obtain the following asymptotic behavior:
\begin{align}\label{AsymptoticConvergence}
\sum\limits_{n = 1}^{+\infty} m_n \cdot \mathbb{E} [\| \nabla F(\vb z_n) \|^2] &\leq \overline{M} + \dfrac{L}{2} \sigma^2 \sum\limits_{n = 1}^{+\infty} R_n < +\infty.
\end{align}
\end{theorem}

Our next result is an observation linked to the long time convergence of the gradient of the objective function, and it follows closely \cite{BertsekasGradient} and \cite{StochasticQuasiNewton}.
\begin{remark}\label{Remark_GradientConvergence}
Suppose that the assumptions from Theorem \eqref{Theorem_AsymptoticBehavior} are satisfied, with respect to the optimization algorithm \eqref{AcceleratedGradientMethod}. In addition, assume that the following property holds:
\begin{align}\label{Assumption_Stepsize}
\sum\limits_{n = 1}^{+\infty} m_n = +\infty,
\end{align}
where $m_n$ was defined in Lemma \eqref{Lemma_AlmostSureConvergence}. Due to \eqref{Limit_Infinity_Gradient}, \eqref{AsymptoticConvergence} and the fact that all of the random variables which are involved in these computations have non-negative values (which leads to the interchange between expected values and infinite summations), one has that
\begin{align*}
\mathbb{E} \left[ \sum\limits_{n = 1}^{+\infty} m_n \cdot \| \nabla F(\vb z_n) \|^2 \right] = \sum\limits_{n = 1}^{+\infty} m_n \cdot \mathbb{E} \left[ \| \nabla F(\vb z_n) \|^2 \right] < +\infty,
\end{align*}
hence
\begin{align*}
\sum\limits_{n = 1}^{+\infty} m_n \| \nabla F(\vb z_n) \|^2 < +\infty, \, (w.p.1).
\end{align*}
Taking into account \eqref{Assumption_Stepsize}, it follows that
\begin{align*}
\liminf\limits_{n \to +\infty} \| \nabla F(\vb z_n) \| = 0, \, (w.p.1).    
\end{align*}
\end{remark}
The second to last conclusion we provide is the convergence rates with regard to the squared norm of the gradient, and it follows the results of \cite{ZouAdam}.
\begin{proposition}\label{ChebyshevProbabilityGradient}
Let $\delta^\prime > 0$ be arbitrarily chosen. Additionally, consider to be true the assumptions used in Theorem \eqref{Theorem_AsymptoticBehavior}. Consider $R_n$ to be an almost sure major bound for $ \| \vb* \lambda_n \|^2 + B \dfrac{C_{n,N}}{\mu_n \Gamma_n} \| \vb* \lambda_n - \vb* \alpha_n \|^2$ as in Theorem \eqref{Theorem_AsymptoticBehavior}. Also, the index $N$ from Theorem \eqref{Theorem_AsymptoticBehavior} is a final iteration, while $n$ can be chosen randomly from the set $\lbrace 1, \ldots, N \rbrace$ with equal probabilities $p_n = \dfrac{1}{N}$. Then, the convergence rates of algorithm \eqref{AcceleratedGradientMethod} holds with probability at least $1 - \delta^\prime$ as below:
\begin{align*}
\min\limits_{n = 1, \ldots, N} \| \nabla F(\vb z_n) \|^2 &\leq \dfrac{\overline{M} + \dfrac{L}{2} \sigma^2 \sum\limits_{n = 1}^{N} R_n }{\delta^\prime \sum\limits_{n = 1}^{N} m_n}.
\end{align*}
\end{proposition}

We end the present section with the following simple remark regarding the Euclidean norm and entrywise Hadamard inequalities.
\begin{remark}\label{Remark_PositivityNormInequality}
Consider $\vb u, \vb v \in \mathbb{R}^d$, such that $0 \leq \vb u \leq \vb v$, namely $0 \leq u_{[j]} \leq v_{[j]}$ for every $j \in \lbrace 1, \ldots, d\rbrace$. Then, we obtain that $0 \leq u_{[1]}^2 + \ldots + u_{[d]}^2 \leq v_{[1]}^2 + \ldots + v_{[d]}^2$, i.e. $\| \vb u \|^2 \leq \| \vb v \|^2$ with respect to the Euclidean $\| \cdot \|_2$ norm.
\end{remark}

\section{Examples of new adaptive momentum methods}\label{Section_ExamplesAdaptiveMethods}

In this section, motivated by \cite{GhadimiLanAccelerated}, we will construct examples of some adaptive Adam-like methods which can be considered new in the Machine Learning community, by using \eqref{ChangesVariablesAlgorithms} from section \eqref{Section_AdaptiveHeavyBall}. We will also deduce the rate of convergence of these methods using the conclusions stated in Theorem \eqref{Theorem_AsymptoticBehavior}. Also, along with the assumptions \emph{$(A_1)$}, \emph{$(A_2)$}, \emph{$(A_3)$}, \emph{$(A_4)$} and \emph{$(A_5)$} that were used in Theorem \eqref{Theorem_AsymptoticBehavior}, we also impose \emph{$(A_6)$} for the almost sure bound of the gradient (a usual hypothesis used for the situations involving adaptive optimizers).
\begin{example}\label{Example_1}
Inspired by the Nesterov implementation method with constant momentum term from \texttt{PyTorch}, we consider $\mu_1 = 1$ and $\mu_n := \mu \in (0, 1)$ for each $n \geq 2$. For every $n \geq 1$, we impose the condition
\begin{align}\label{eq:conditions_adaptive_coeffs}
\vb* \lambda_n = M \vb* \alpha_n,
\end{align}
where \eqref{eq:conditions_adaptive_coeffs} hold in the sense of the Hadamard notations, and $M$ is chosen such that $M > 0$, hence $\vb* \lambda_n > 0$ for each $n \geq 1$ (since we will take $\vb* \alpha_n > 0$ for every $n \geq 1$). Furthermore, we mention that we can take $M \geq \dfrac{1}{\mu}$ similarly with \cite{GhadimiLanAccelerated}, hence $M \geq 1$ in order to have $\vb* \lambda_n - \vb* \alpha_n \geq 0$, for each $n \geq 1$. Also, we mention that the limiting case when $M = \dfrac{1}{\mu}$ is the one where $\vb g_{n-1}$ doesn't appear but $\vb g_n$ does in the formula of the \eqref{AcceleratedGradientMethod} iterative method. But, since our computations involve only the squared term $[\vb* \lambda_n - \vb* \alpha_n]^2$ and the norm $\| \vb* \lambda_n - \vb* \alpha_n \|^2$, then it does not matter if $\vb* \lambda_n$ is greater or not than $\vb* \alpha_n$. \\
By considering \eqref{eq:Gamma}, we have that $\Gamma_1 = 1$ and $\Gamma_2 = (1 - \mu) \Gamma_1$, i.e. $\Gamma_2 = 1 - \mu$. By induction, we can prove that $\Gamma_n = (1-\mu)^{n-1}$ for all $n \geq 2$ (evidently, this formula holds also for $n = 1$, hence it is valid from $n \geq 1$). Using \eqref{C_nN}, for every $n \geq 1$, it follows that
\begin{align*}
C_{n, N} = \sum\limits_{j = n}^{N} \Gamma_j = \sum\limits_{j = n}^{N} (1-\mu)^{j-1} = (1-\mu)^{n-1} \cdot \dfrac{1-(1-\mu)^{N-n}}{\mu} \leq \dfrac{(1-\mu)^{n-1}}{\mu}.
\end{align*}
Now, we turn our attention to
\begin{align*}
\vb D_n = \vb* \lambda_n - L [\vb* \lambda_n]^2 - \dfrac{LB}{2} \dfrac{C_{n,N}}{\mu_n \Gamma_n} [\vb* \lambda_n - \vb* \alpha_n]^2 &= \vb* \lambda_n \odot \left( 1 - L \vb* \lambda_n - \dfrac{LB}{2} \dfrac{C_{n,N}}{\mu_n \Gamma_n \vb* \lambda_n} \odot [\vb* \lambda_n - \vb* \alpha_n]^2 \right) \\
&= \vb* \lambda_n \odot \left( 1 - L \left( \vb* \lambda_n + B \dfrac{C_{n,N}}{2 \mu_n \Gamma_n \vb* \lambda_n} \odot [\vb* \lambda_n - \vb* \alpha_n]^2 \right) \right) := \vb* \lambda_n \odot \vb D_{n,N}.
\end{align*}
Now, $[\vb* \lambda_n - \vb* \alpha_n]^2 = (M - 1)^2 [\vb* \alpha_n]^2$, so we have, for every $n \geq 2$, the following almost sure computations:
\begin{align*}
\vb D_{n,N} &= 1 - L \left( \vb* \lambda_n + B \dfrac{C_{n,N}}{2 \mu \Gamma_n \vb* \lambda_n} \odot [\vb* \lambda_n - \vb* \alpha_n]^2 \right) \geq 1 - L \vb* \alpha_n \left( M + \dfrac{B}{2 \mu^2} \dfrac{(M-1)^2}{M} \right), \, (w.p.1),
\end{align*}
where we have used that $\dfrac{C_{n,N}}{\mu \Gamma_n} \leq \dfrac{1}{\mu^2}$ for every $n \geq 2$, and where the constant $B$ will be determined later by us. Further, let's impose the Hadamard-type almost sure inequality for $n \geq 1$
\begin{align}\label{eq:cond_3}
\vb* \alpha_n \leq \dfrac{1}{2 L \left( M + \dfrac{B}{2 \mu^2} \dfrac{(M-1)^2}{M} \right)}, \, (w.p.1),    
\end{align}
which leads, for every $n \geq 1$, to
\begin{align*}
\vb D_n = \vb* \lambda_n \odot \vb D_{n,N} \geq \dfrac{\vb* \lambda_n}{2} = \dfrac{M \vb* \alpha_n}{2}, \, (w.p.1),    
\end{align*}
along with $\vb D_{n,N} \geq 1 - L \left( M + \dfrac{B}{2 \mu^2} \dfrac{(M-1)^2}{M} \right) \vb* \alpha_n \geq \dfrac{1}{2} > 0$. It is trivial to observe if we make similar computations as before, taking into account the equality $\mu_1 = 1$, then it follows that $\vb D_{1,N} \geq 1 - \left( M + \dfrac{B}{2 \mu} \dfrac{(M - 1)^2}{M} \right) \vb* \alpha_1 \geq 1 - \left( M + \dfrac{B}{2 \mu^2} \dfrac{(M-1)^2}{M} \right) \vb* \alpha_1 \geq \dfrac{1}{2}$, since $\mu^2 \leq \mu \in (0, 1)$, hence \eqref{eq:cond_3} and the lower bound for $\vb D_{1,N}$ are satisfied also for $n = 1$. \\
Furthermore, in our next computations we will consider for each $n \geq 1$, by referring to Lemma \eqref{Lemma_AlmostSureConvergence}, the following scalar term: 
\begin{align}\label{eq:cond_4}
m_n \in \left(0, \dfrac{M \vb* \alpha_n}{2} \right].
\end{align}
By the fact that $\| \vb* \lambda_n - \vb* \alpha_n \|^2 = \| (M-1) \vb* \alpha_n \|^2 = (M - 1)^2 \| \vb* \alpha_n \|^2$ then, for our next computations related to \eqref{MinimumGradientRate}, we have that
\begin{align*}
\sum\limits_{n = 1}^{N} \left[ \| \vb* \lambda_n \|^2 + B \dfrac{C_{n,N}}{\mu_n \Gamma_n} \| \vb* \lambda_n - \vb* \alpha_n \|^2 \right] &\leq \sum\limits_{n = 1}^{N} \left( M^2 \| \vb* \alpha_n \|^2 + \dfrac{B (M-1)^2}{\mu^2} \| \vb* \alpha_n \|^2 \right) = \left( M^2 + \dfrac{B (M-1)^2}{\mu^2} \right) \cdot \| \vb* \alpha_n \|^2,
\end{align*}
where, as we have said before, $B$ will be specified later. In what follows, we will need to define $R_n$ to be an upper bound for $\left( M^2 + \dfrac{B (M-1)^2}{\mu^2} \right) \cdot \| \vb* \alpha_n \|^2$. Evidently, we have used that $\dfrac{C_{n,N}}{\mu_n \Gamma_n} \leq \dfrac{1}{\mu^2}$ for $n \geq 2$, while for the situation when $n = 1$, we obtain $\dfrac{C_{1,N}}{\mu_1 \Gamma_1} \leq \dfrac{1}{\mu} \leq \dfrac{1}{\mu^2}$, hence the major bound $\dfrac{1}{\mu^2}$ was used for all $n \geq 1$ in the above computations. \\
Now, let the \emph{adaptive stepsize} $\vb* \alpha_n = \nu \vb a_n$, where we have taken $\vb \nu_n = \nu$ and $\vb r_n = \vb a_n$. Furthermore, in the following, we introduce the so-called \emph{effective stepsize}
\begin{align}\label{a_n}\tag{EffectiveStepsize}
\vb a_n = \dfrac{\eta_n}{\varepsilon + \sqrt{\vb v_n}},    
\end{align}
where $\eta_n > 0$ for every $n \geq 1$, and where we define, as in \cite{ChenAdam}, $\vb v_n$ to be the exponential moving average of the squared gradients adjusted by a \emph{max} term related to an \emph{AMSGrad} formulation, namely
\begin{align}\label{v_n}\tag{AccumulatedSqGrad}
\begin{cases}
\tilde{\vb v}_n = \beta_2 \tilde{\vb v}_{n-1} + (1-\beta_2) [\vb g_n]^2 \\
\vb v_n = \max \lbrace \vb v_{n-1}, \tilde{\vb v}_n \rbrace,
\end{cases}
\end{align}
where we consider the initial data $\vb v_0 = \tilde{\vb v}_0 = 0$, and the coefficient $\beta_2 \in (0, 1)$. The first thing we need to do is to get an upper bound for \eqref{a_n}. Firstly, we will show that $\tilde{\vb v}_n$ is componentwise bounded by $\mathcal{K}^2$, under assumption \emph{$(A_6)$}. Working with the Euclidean norm, employing Proposition \eqref{Proposition_SquaredVector} and taking into consideration that every component of a random vector with non-negative values is less than its Euclidean norm (thus the vector is less than its Euclidean norm in a Hadamard sense), we have that $\tilde{\vb v}_1 = \beta_2 \tilde{\vb v}_0 + (1-\beta_2) [\vb g_1]^2 = (1-\beta_2) [\vb g_1]^2 \leq (1-\beta_2) \| [\vb g_1]^2 \| \leq (1-\beta_2) \| \vb g_1 \|^2 \leq \mathcal{K}^2$. Similarly, $\tilde{\vb v}_2 = \beta_2 \tilde{\vb v}_1 + (1-\beta_2) [\vb g_2]^2 \leq \beta_2 \mathcal{K}^2 + (1-\beta_2) \| [\vb g_2]^2 \| \leq \beta_2 \mathcal{K}^2 + (1 - \beta_2) \| \vb g_2 \|^2 \leq \mathcal{K}^2$. Now, we will proceed by induction. Suppose that $\tilde{\vb v}_n \leq \mathcal{K}^2$. Then, utilizing again Proposition \eqref{Proposition_SquaredVector}, we get that $\tilde{\vb v}_{n+1} = \beta_2 \tilde{\vb v}_n + (1-\beta_2) [\vb g_{n+1}]^2 \leq \beta_2 \mathcal{K}^2 + (1-\beta_2) \| [\vb g_{n+1}]^2 \| \leq \beta_2 \mathcal{K}^2 + (1-\beta_2) \| \vb g_{n+1} \|^2 \leq \beta_2 \mathcal{K}^2 + (1-\beta_2) \mathcal{K}^2 = \mathcal{K}^2$. Now, we shift our attention to determining a major bound for $\vb v_n$. First of all, we have that $\vb v_1 = \max \lbrace \vb v_0, \tilde{\vb v}_1 \rbrace = \max \lbrace 0, \tilde{\vb v}_1 \rbrace \leq \mathcal{K}^2$. Analogous, $\vb v_2 = \max \lbrace \vb v_1, \tilde{\vb v}_2 \rbrace \leq \mathcal{K}^2$, since $\vb v_1 \leq \mathcal{K}^2$ and also $\tilde{\vb v}_2 \leq \mathcal{K}^2$. It is easy to see that we can proceed as before and show by induction that $\vb v_n \leq \mathcal{K}^2$ for every $n \geq 1$. \\
By taking into consideration the above analysis, we compute
\begin{align*}
\| \vb* \alpha_n \|^2 = \nu^2 \| \vb a_n \|^2 = \nu^2 \Big\| \dfrac{\eta_n}{\sqrt{\vb v}_n + \varepsilon} \Big\|^2 = \nu^2 \eta_n^2 \Big\| \dfrac{1}{\sqrt{\vb v}_n + \varepsilon} \Big\|^2 \leq d \left( \dfrac{\nu}{\varepsilon} \right)^2 \eta_n^2,
\end{align*}
hence we can take $R_n := d \left( \dfrac{\nu}{\varepsilon} \right)^2 \left( M^2 + \dfrac{B(M - 1)^2}{\mu^2} \right) \eta_n^2$. We also have the following computations:
\begin{align}\label{eq:bound_mn}
\dfrac{M \vb* \alpha_n}{2} = \dfrac{M \nu \eta_n}{2 (\sqrt{\vb v}_n+\varepsilon)} \geq \dfrac{M \nu \eta_n}{2(\mathcal{K} + \varepsilon)} := m_n, 
\end{align}
We notice that $m_n > 0$ because $\nu > 0$ and $\eta_n > 0$. Then, the convergence rate of the minimum values of the true gradient given in \eqref{MinimumGradientRate} merely turns into
\begin{align}\label{OrderConv_Bound}\tag{OrderConv}
\mathbb{E} \left[ \min_{n = 1, \ldots, N} \| \nabla F(\vb z_n) \|^2 \right] &\leq \dfrac{\overline{M} + \dfrac{L \sigma^2}{2} \left( M^2 + \dfrac{B(M-1)^2}{\mu^2} \right) d \left( \dfrac{\nu}{\varepsilon} \right)^2 \sum\limits_{n = 1}^{N} \eta_n^2}{ \dfrac{M \nu}{2 (\mathcal{K} + \varepsilon)} \sum\limits_{n = 1}^{N} \eta_n}.
\end{align}
We will analyze two cases for our (possibly) adaptive methods as follows. For the finite-time horizon situation, we consider the constant stepsize that depends on the final iteration $N$, namely $\eta_n = \dfrac{C}{\sqrt{N}}$, where $C > 0$. Then, we obtain the rate of convergence
\begin{align*}
\mathbb{E} \left[ \min_{n = 1, \ldots, N} \| \nabla F(\vb z_n) \|^2 \right] &\leq \dfrac{\overline{M} + \dfrac{Ld \sigma^2}{2} \left( \dfrac{\nu}{\varepsilon} \right)^2 \left( M^2 + \dfrac{B(M-1)^2}{\mu^2} \right) \sum\limits_{n = 1}^{N} \dfrac{C^2}{N}}{ \dfrac{M \nu}{2 (\mathcal{K} + \varepsilon)} \sum\limits_{n = 1}^{N} \dfrac{C}{\sqrt{N}}} \\
&= 2 (\mathcal{K}+\varepsilon) \dfrac{\overline{M} + \dfrac{Ld \sigma^2}{2} C^2 \left( \dfrac{\nu}{\varepsilon} \right)^2 \left( M^2 + \dfrac{B(M-1)^2}{\mu^2} \right)}{M C \nu} \cdot \dfrac{1}{\sqrt{N}} \\
&= \mathcal{O} \left( \dfrac{1}{\sqrt{N}} \right).
\end{align*}
On the other hand, for the situation involving a decreasing stepsize, we consider $\eta_n = \dfrac{C}{\sqrt{n}}$ where $C > 0$. Then, taking into account that $\sum\limits_{n=1}^{N} \dfrac{1}{n} \leq 1 + log(N)$ and denoting $\overline{Q} := d \left( \dfrac{\nu}{\varepsilon} \right)^2 \left( M^2 + \dfrac{B(M-1)^2}{\mu^2} \right)$, it follows that
\begin{align*}
\mathbb{E} \left[ \min_{n = 1, \ldots, N} \| \nabla F(\vb z_n) \|^2 \right] &\leq \dfrac{\overline{M} + \dfrac{L \overline{Q} \sigma^2}{2} \sum\limits_{n = 1}^{N} \dfrac{C^2}{n}}{ \dfrac{M \nu}{2 (\mathcal{K}+\varepsilon)} \sum\limits_{n = 1}^{N} \dfrac{C}{\sqrt{n}}} \leq \dfrac{\overline{M} + \dfrac{L \overline{Q} C^2 \sigma^2}{2} (1 + log(N)) }{ \dfrac{M \nu}{2 ( \mathcal{K} + \varepsilon)} \sum\limits_{n = 1}^{N} \dfrac{C}{\sqrt{n}}} \\
&\leq \dfrac{2 (\mathcal{K}+\varepsilon)}{C M \nu} \dfrac{ \left( \overline{M} + \dfrac{L \overline{Q} C^2 \sigma^2}{2} \right) + \dfrac{L \overline{Q} C^2 \sigma^2}{2} log(N) }{\sum\limits_{n=1}^{N} \dfrac{1}{\sqrt{N}}} \\
&\leq \dfrac{1}{\sqrt{N}} \left( Q_1 + Q_2 log(N) \right) \leq \mathcal{O} \left( \dfrac{log(N)}{\sqrt{N}} \right),
\end{align*}
for some constants $Q_1$ and $Q_2$ independent of the final iteration $N$. \\
Finally, we end this example by emphasizing a bound for the sequence $(\eta_n)_{n \geq 1}$ related to the almost sure inequality \eqref{eq:cond_3}. For the case when $\eta_n = \dfrac{C}{\sqrt{N}}$ and using the fact that $\vb* \alpha_n = \nu \vb a_n = \nu \dfrac{\eta_n}{\sqrt{\vb v_n} + \varepsilon} \leq \dfrac{\nu \eta_n}{\varepsilon}$, we find that
\begin{align*}
\dfrac{C}{\sqrt{N}} \leq \dfrac{\varepsilon}{2 L \nu \left( M + \dfrac{B}{2 \mu^2} \dfrac{(M-1)^2}{M} \right)},
\end{align*}
hence $N \geq RHS$, where
\begin{align*}
RHS := \left[ \dfrac{2C L \nu}{\varepsilon} \left( M + \dfrac{B}{2 \mu^2} \dfrac{(M-1)^2}{M} \right) \right]^2.    
\end{align*}
Thus, the link between the last iteration $N$ and the constant $C > 0$ is obtained. On the other hand, for the case when $\eta_n = \dfrac{C}{\sqrt{n}}$, we can clearly see that we can impose $n \geq 1 \geq RHS$, which is a relationship between $\nu$ and $C$ that encompass the constant component of the learning rate of the algorithm.
\end{example}
Similarly to \cite{ChenAdam}, as a consequence of Theorem \eqref{Theorem_AsymptoticBehavior}, we get the following corollary which implies that the rate of convergence of \eqref{AdaptiveAcceleratedMomentumMethodShiftedUpdates} is the same as \emph{SGD} and \emph{AMSGrad} in the stochastic non-convex setting.
\begin{corollary}\label{Corollary_RateConvergence}
Let the adaptive algorithm \eqref{AdaptiveAcceleratedMomentumMethodShiftedUpdates} given through the effective stepsize $\vb a_n$ from \eqref{a_n} with the moving average of the accumulated squared gradients $\vb v_n$ given in \eqref{v_n}. Considering as in Lemma \eqref{Lemma_AlmostSureConvergence}, a natural number $N$ such that $N \geq n \geq 1$ to be a final time iteration, and a constant $C > 0$, then for $\eta_n = \dfrac{C}{\sqrt{N}}$, one has that $\mathbb{E} \left[ \min\limits_{n = 1, \ldots, N} \| \nabla F(\vb z_n) \|^2 \right] \leq \mathcal{O} \left( \dfrac{1}{\sqrt{N}} \right)$. On the other hand, for the choice of the non-constant stepsize $\eta_n = \dfrac{C}{\sqrt{n}}$, we have $\mathbb{E} \left[ \min\limits_{n = 1, \ldots, N} \| \nabla F(\vb z_n) \|^2 \right] \leq \mathcal{O} \left( \dfrac{log(N)}{\sqrt{N}} \right)$. 
\end{corollary}
At last, we consider the following two general remarks concerning adaptive methods with shifted updates.
\begin{remark}
In Example \eqref{Example_1} and in Lemma \eqref{Lemma_AlmostSureConvergence}, we have that $\dfrac{\| 1 - \tilde{\vb* \gamma}_n \|^2}{(1-\mu_n)} \leq B$. Taking into account Remark \eqref{Remark_NormBound}, and Example \eqref{Example_1}, we can take $\tilde{\vb* \gamma}_n \in \mathbb{R}$, such that $\dfrac{( 1 - \tilde{\vb* \gamma}_n )^2}{(1-\mu_n)} \leq B$ for every $n \geq 2$. Evidently, we take $\tilde{\vb* \gamma}_n \geq \mu_n$ for $n \geq 2$, and also the simplified version (as we shall do in the last section) $\tilde{\vb* \gamma}_n = \tilde{\gamma}$ for $n \geq 2$, such that $0 < \mu \leq \tilde{\gamma} < 1$ and $\tilde{\gamma}_1 = \mu_1 = 1$. Therefore, we get for each $n \geq 2$ that $\dfrac{( 1 - \tilde{\gamma} )^2}{(1-\mu)} = \left( \dfrac{1-\tilde{\gamma}}{1-\mu} \right)^2 \cdot (1 - \mu) \leq 1$, thus we find that the constant $B$ from Lemma \eqref{Lemma_AlmostSureConvergence} and Example \eqref{Example_1} satisfies $B = 1$ and does not depend on the dimension $d$. 
\end{remark}
\begin{remark}\label{Adaptive_Gamma_n}
For the case when the stepsize $\vb* \alpha_n$ is not adaptive, i.e. $\vb* \alpha_n \in \mathbb{R}$ is a deterministic given stepsize, it is trivial to observe that the rate of convergence from Corollary \eqref{Corollary_RateConvergence} remains the same. But, in the case when $\vb* \alpha_n$ is not adaptive, the only method in which we can make \eqref{AdaptiveAcceleratedMomentumMethodShiftedUpdates} to become an adaptive method is to take $\tilde{\vb* \gamma}_n$ to be adaptive, hence to depend on the moving average of accumulated squared gradients $\vb v_n$. Since $\tilde{\vb* \gamma}_n = \vb* \gamma_n \dfrac{\mu_{n-1}}{1 - \mu_{n-1}}$ for $n \geq 2$ from \eqref{ChangesVariablesAlgorithms} (using that $\mu_n = \mu$ for all $n \geq 2$), we get that $\tilde{\vb* \gamma}_n = \vb* \gamma_n \dfrac{\mu}{1 - \mu}$. Hence $\dfrac{\| 1 - \tilde{\vb* \gamma}_n \|^2}{(1-\mu_n)} \leq B$ becomes $\| 1 - \mu(1 + \vb* \gamma_n) \|^2 \leq B(1 - \mu)^3$ for every $n \geq 2$. Let's take the case when $1 - \mu(1 + \vb* \gamma_n) \geq 0$ which is equivalent to $\vb* \gamma_n \leq \dfrac{1-\mu}{\mu}$. On the other hand, $1 - \mu(1 + \vb* \gamma_n) \leq 1$ for each $n \geq 2$. Hence, applying Remark \eqref{Remark_PositivityNormInequality}, we impose $\| 1 - \mu(1 + \vb* \gamma_n) \|^2 \leq d \leq B(1 - \mu)^3$, hence we can take $B = \dfrac{d}{(1 - \mu)^3}$. By letting $\vb* \gamma_n = \dfrac{\pi_n}{\sqrt{\vb v_n} + \varepsilon} \leq \dfrac{\pi_n}{\varepsilon}$, where $\varepsilon > 0$ is small enough in order to not interfere with the adaptive part of $\vb* \gamma_n$, and $\pi_n > 0$, then utilizing the bound $\vb* \gamma_n \leq \dfrac{1-\mu}{\mu}$, we choose $\pi_n$ such that $\pi_n \leq \varepsilon \dfrac{1-\mu}{\mu}$, hence by taking $\pi_n := \varepsilon \dfrac{1-\mu}{\mu}$ for every $n \geq 2$, then $\vb* \gamma_n = \dfrac{1-\mu}{\mu} \cdot \dfrac{\varepsilon}{\sqrt{\vb v}_n + \varepsilon}$. Finally, it is worth mentioning that there are two disadvantages of this situation with respect to our theoretical analysis (this is the actual reason why we did not used this in our implementations): firstly, the constant $B$ is very large and affects the rate of convergence that we have obtained and second of all the constant coefficient of $\vb* \gamma_n$ is $\dfrac{\varepsilon(1-\mu)}{\mu}$ which has a very low value and practically neglects the adaptive structure of $\vb* \gamma_n$. More specifically, in Example \eqref{Example_1}, the coefficient from the numerator of the right hand side of the convergence bound \eqref{OrderConv_Bound} is proportional to $\dfrac{(d \nu)^2}{\varepsilon^2 \mu^2 (1-\mu)^3}$, hence we can take $\dfrac{(d \nu)^2}{\varepsilon^2 \mu^2 (1-\mu)^3} \leq 1$ in order to have a low value for the constant appearing in the right hand side of \eqref{OrderConv_Bound}, so $d \leq \dfrac{\varepsilon \mu (1-\mu)}{\nu} \cdot \sqrt{1 - \mu}$, which actually imposes a constraint on the dimension $d$ of the algorithms parameters (weights $\&$ biases in the cases of neural networks).
\end{remark}

We are finally arriving at our last result of the present section in which we focus on the complexity result for a particular case of our adaptive momentum algorithm associated with a given final time iteration as in \cite{GadatGavra}. In this result, we consider the choices of the coefficients that are used in Example \eqref{Example_1}.
\begin{proposition}\label{Proposition_ComplexityResult}
We take the algorithm \eqref{AcceleratedGradientMethod}, along with the assumptions from Theorem \eqref{Theorem_AsymptoticBehavior}. Consider a given integer $N \geq 1$ and $n$ an integer sampled uniformly from $\lbrace 1, \ldots, N \rbrace$. With respect to the section \eqref{Section_AssumptionSetting} we take $\sigma$ to possibly depend on $d$. Similarly with Example \eqref{Example_1}, we assume that $\dfrac{C_{n, N}}{\mu_n \Gamma_n} = \hat{C} $ for some $\hat{C} > 0$ and for every $n \geq 1$. Let $\eta_n$ from the definition of $\vb a_n$ from \eqref{a_n} with $\eta_n = \dfrac{C}{\sqrt{N}}$ for some $C > 0$ that may depend on the dimension $d$, and where $M > 0$. \\ 
In addition, for every $n \geq 1$, let $m_n = \dfrac{M \nu \eta_n}{2 (\mathcal{K}+\varepsilon)} \leq \dfrac{M}{2} \vb* \alpha_n$ as in \eqref{eq:bound_mn}, $\vb* \lambda_n = M \vb* \alpha_n$ and $\vb* \alpha_n = \nu \vb a_n$ where $\nu > 0$ and $\vb a_n$ is defined by \eqref{a_n}. Moreover, for a given $\delta > 0$, in order to guarantee that $\dfrac{1}{N} \sum\limits_{n = 1}^{N} \mathbb{E} \left[ \| \nabla F(\vb z_n) \|^2 \right] \leq \delta$, then the number of iterations needed is at least $\mathcal{O} \left( \dfrac{1}{\delta^2} \right)$.
\end{proposition}

Finally, focusing on the proof of the previous result, we mention the following remark regarding the dimension dependent coefficients used for the adaptive sequences in the algorithm \eqref{AcceleratedGradientMethod}.
\begin{remark}\label{RemarkIterationComplexity}
If we proceed as in \cite{GadatGavra} concerning our iteration complexity result, we can take $\sigma^2 = \dfrac{1}{d}$, hence the orders for the iteration $N$ do not depend on the dimension $d$. On the other hand, if we take $\sigma^2$ to be independent of $d$, then it is easy to observe that, from a practical point of view, we must consider $B$ or $\hat{C}$ to dependent on the dimension $d$ and to be proportional to $\dfrac{1}{d}$. This actually resembles the choice of the parameters from the iteration complexity of Theorem 2 of the aforementioned paper \cite{GadatGavra}, and it affects also the theoretical rate of convergence through the term $B$ and the choices of the parameters, since in Example \eqref{Example_1}, $\hat{C} \leq 1 / \mu$. More precisely, for the latter case we remind that $\dfrac{C_{n,N}}{\mu_n \Gamma_n}$ is less or equal than $\dfrac{1}{\mu} \leq \dfrac{1}{\mu^2}$ for $n = 1$, and less or equal than $\dfrac{1}{\mu^2}$ for $n \geq 2$. 
\end{remark}

\section{The \texttt{PyTorch} equivalent formulation for Algorithm \eqref{AdaptiveAcceleratedMomentumMethodShiftedUpdates}}\label{EquivalentForms}

In the present section we shall strictly follow the research paper \cite{DefazioCurved} of Defazio, which shows the equivalence between the original Nesterov method and some alternative formulations: the Sutskever algorithm, and also with the modern momentum form which currently is used in \texttt{PyTorch} \footnote{\url{https://github.com/pytorch/pytorch/blob/master/torch/optim/sgd.py}}. We mention that in our \texttt{PyTorch} implementation we shall use the Sutskever formulation (and not the modern momentum one) due to its simplicity. Our aim is to show that Algoritm \eqref{AdaptiveHeayBallMethod} can be restated into an equivalent form, and then we show that this form can also be used in the \emph{AMSGrad}-like algorithm \eqref{AdaptiveAcceleratedMomentumMethodShiftedUpdates}. Also, we specify that we shall use $\tilde{\vb* \gamma}_n \in \mathbb{R}^d$ even though in our numerical experiments we will take $\tilde{\vb* \gamma}_n$ to be a deterministic real number (we have already pointed out in Remark \eqref{Adaptive_Gamma_n}, and also see section \eqref{Section_NeuralNetworks}). Moreover, we will utilize the Sutskever approach to our \texttt{PyTorch} implementation involving the training of neural networks.  \\
In order to do this, we consider the adaptive accelerated algorithm \eqref{AdaptiveHeayBallMethod}. For simplicity, taking $n \geq 2$, we restate it below:
\begin{align*}
\begin{cases}
\vb y_n = \vb* \theta_n + \beta_n (\vb* \theta_n - \vb* \theta_{n-1}) - \mu_n \vb* \omega_n \odot \vb* g_{n-1} \\ \\
\vb z_n = \vb* \theta_n + \vb* \gamma_n \odot (\vb* \theta_n - \vb* \theta_{n-1}) - \tilde{\vb* \gamma}_n \odot \vb* \omega_n \odot \vb* g_{n-1} \\ \\
\vb* \theta_{n+1} = \vb y_n - \vb* \alpha_n \odot \vb* g_n
\end{cases}
\end{align*}
We take into account the following variable change:
$$ \vb m_{n+1}^S = \vb* \theta_{n+1} - \vb* \theta_n. $$
Then, for $n \geq 2$, we obtain that
\begin{align*}
\vb m_{n+1}^S &= [\vb y_n - \vb* \alpha_n \odot \vb g_n] - \vb* \theta_n = (\vb y_n - \vb* \theta_n) - \vb* \alpha_n \odot \vb g_n \\
&= \beta_n(\vb* \theta_n - \vb* \theta_{n-1}) - \mu_n \vb* \omega_n \odot \vb g_{n-1} - \vb* \alpha_n \odot \vb g_n \\
&= \beta_n \vb m_{n}^S - \mu_n \vb* \omega_n \odot \vb g_{n-1} - \vb* \alpha_n \odot \vb g_n.   
\end{align*}
By taking $\vb* \theta_n^S := \vb* \theta_n$, it implies that $\vb* \theta_{n+1}^S = \vb* \theta_n^S + \vb m_{n+1}^S.$ Thus, we are able to find the Sutskever form of Algorithm \eqref{AdaptiveHeayBallMethod}, where the key part is that the inertial term containing $\vb* \theta_n - \vb* \theta_{n-1}$ is replaced by a momentum term $\vb m_n^S$. By taking $\vb z_n^S := \vb z_n$, the algorithm takes the following form for every $n \geq 2$.
\begin{equation}\label{AdaptiveSutskever}\tag{SutskeverForm}
\begin{aligned}
\begin{cases}
& \vb z_n^S = \vb* \theta_n^S + \vb* \gamma_n \odot \vb m_n^S - \tilde{\vb* \gamma}_n \odot \vb* \omega_n \odot \vb g_{n-1} \\ \\
&\vb m_{n+1}^S = \beta_n \vb m_{n}^S - \mu_n \vb* \omega_n \odot \vb g_{n-1} - \vb* \alpha_n \odot \vb g_n \\ \\
&\vb* \theta_{n+1}^S = \vb* \theta_n^S + \vb m_{n+1}^S.
\end{cases}
\end{aligned}
\end{equation}
Now, the next step is to show that Algorithm \eqref{AdaptiveSutskever} can be written into a form that can be easily used for the implementations utilized in neural network training. Even though modern momentum change of variables is used in \texttt{PyTorch} and \texttt{Tensorflow} implementations, where the momentum term is $\vb m_n^M = - \dfrac{1}{\vb* \alpha_{n-1}} \odot \vb m_n^S$ for every $n \geq 2$, we stick with the Sutskever formulation since, for this method, the adaptive momentum term is relatively simple to implement. More precisely, in the modern momentum form, we will obtain ratios of the form $\dfrac{\vb* \alpha_{n-1}}{\vb* \alpha_n}$, i.e. $\vb m_{n+1}^M = \beta_n \dfrac{\vb* \alpha_{n-1}}{\vb* \alpha_n} \odot \vb m_n^M + \vb g_n + \dfrac{\mu_n \vb* \omega_n}{\vb* \alpha_n} \odot \vb g_{n-1}$, which shall complicate our neural network implementation. \\
We turn our focus to the fact that we know that, for every $n \geq 2$, we obtain that
$$ \vb* \theta_{n}^S = \vb z_n^S - \vb* \gamma_n \odot \vb m_n^S + \tilde{\vb* \gamma}_n \odot \vb* \omega_n \odot \vb g_{n-1}. $$
By the fact that, for every $n \geq 2$, we have
$$ \vb* \theta_{n+1}^S = \vb* \theta_n^S + \vb m_{n+1}^S, $$
then, for each $n \geq 2$, it follows that
\begin{align*}
\vb z_{n+1}^S - \vb* \gamma_{n+1} \odot \vb m_{n+1}^S + \tilde{\vb* \gamma}_{n+1} \odot \vb* \omega_{n+1} \odot \vb g_n = \vb z_n^S - \gamma_n \odot \vb m_n^S + \tilde{\vb* \gamma}_n \odot \vb* \omega_n \odot \vb g_{n-1} + \vb m_{n+1}^S,    
\end{align*}
which shows that
\begin{align*}
\vb z_{n+1}^S &= \vb z_n^S + \left( \vb m_{n+1}^S - \vb* \gamma_n \odot \vb m_n^S + \tilde{\vb* \gamma}_n \odot \vb* \omega_n \odot \vb g_{n-1} - \tilde{\vb* \gamma}_{n+1} \odot \vb* \omega_{n+1} \odot \vb g_n \right) + \vb* \gamma_{n+1} \odot \vb m_{n+1}^S \\
&= ( \vb z_n^S + ( \vb m_{n+1}^S - \beta_n \vb m_n^S ) ) + \vb* \gamma_{n+1} \odot \vb m_{n+1}^S + (\beta_n - \vb* \gamma_n) \odot \vb m_n^S + \left( \tilde{\vb* \gamma}_n \odot \vb* \omega_n \odot \vb g_{n-1} - \tilde{\vb* \gamma}_{n+1} \odot \vb* \omega_{n+1} \odot \vb g_n \right).
\end{align*}
Henceforth, for every $n \geq 2$, one has that
\begin{align*}
\vb z_{n+1}^S &= \vb z_n^S - \vb* \alpha_n \odot \vb g_n + \vb* \gamma_{n+1} \odot \vb m_{n+1}^S + ( \beta_n - \vb* \gamma_n) \odot \vb m_n^S + \left[ (\tilde{\vb* \gamma}_n - \mu_n) \odot \vb* \omega_n \odot \vb g_{n-1} - \tilde{\vb* \gamma}_{n+1} \odot \vb* \omega_{n+1} \odot \vb g_n \right].    
\end{align*}

\subsection{The shifted AMSGrad-like algorithm}\label{Subsection_AdamLikeAlgorithms}
We will start by describing the Sutskever formulation of our adaptive methods in an algorithmic framework. We then go on to present the connection to the \eqref{AcceleratedGradientMethod} and the selection of the non-adaptive coefficients. As noted in Remark \eqref{Remark_NormBound}  we consider $\gamma_n$ and $\tilde{\gamma}_n$ to be deterministic real numbers. By considering Example \eqref{Example_1}, we let $\vb* \lambda_n = M \vb* \alpha_n$, where $M > 0$ for every $n \geq 1$. Also, from the definition of the adaptive stepsize, we shall take $\vb* \alpha_n = \nu \vb a_n$ where $\nu \in (0, 1)$ with $\vb a_n$ being defined in \eqref{a_n} for each $n \geq 1$. Also, we will take $\eta_n = \eta$ for each $n \geq 1$ in order to be analogous with Remark \eqref{RemarkIterationComplexity} where $\eta_n = \dfrac{C}{\sqrt{N}}$. At the same time, because we will take $\vb* \gamma_n, \tilde{\vb* \gamma}_n \in \mathbb{R}$, then for these deterministic scalars we will not employ the Hadamard notations, i.e. entrywise multiplication.
From \eqref{AdaptiveHeayBallMethod}, we know that $\beta_n=\dfrac{\mu_n}{\mu_{n-1}}(1-\mu_{n-1})$ for each $n \geq 2$. On the other hand, we know that 
\begin{align*}
    \mu_n = 
    \begin{cases}
    1; & n = 1 \\
    \mu \in (0,1); & n \geq 2.
    \end{cases}
\end{align*}
The analysis that was done above led us to 
\begin{align*}
    \beta_n = 
    \begin{cases}
    0; & n = 2 \\
    1-\mu; & n \geq 3.
    \end{cases}
\end{align*}
Also, from \eqref{ChangesVariablesAlgorithms}, we have that $ \gamma_n = \dfrac{\tilde{ \gamma}_n}{\mu_{n-1}}(1-\mu_{n-1})$ for each $n \geq 2$. We know that $\tilde{\gamma}_1 = 1$ and that $\tilde{\gamma}_n = \tilde{\gamma}$ for every $n \geq 2$, namely
\begin{align*}
    \tilde{\gamma}_n = 
    \begin{cases}
    1; & n = 1 \\
    \tilde{\gamma}; & n \geq 2,
    \end{cases}
\end{align*}
where $\tilde{\gamma} \geq \mu$.
Since, for $ n = 2$, $\gamma_2 = \dfrac{\tilde{\gamma}_2}{\mu_1}(1-\mu_1) = 0$, then it follows that 
\begin{align*}
    \gamma_n = 
    \begin{cases}
    0; & n = 2 \\
    \tilde{\gamma} \dfrac{1-\mu}{ \mu }; & n \geq 3.
    \end{cases}
\end{align*}
From the definition of $\vb* \omega_n$ given in \eqref{omega_n}, we find that $ \omega_n = \left( M - \dfrac{1}{\mu_{n-1}} \right) \vb* \alpha_{n-1}$, for each $ n \geq 2$.\\
Now, for each $n \geq 2$, we have the following 
\begin{align*}
   \vb m_{n+1}^S = \beta_n \vb m_{n}^S - \mu_n \vb* \omega_n \odot \vb g_{n-1} - \vb* \alpha_n \odot \vb g_n.
\end{align*}
Using the calculation from above regarding $\vb* \omega_n$, for every $n \geq 2$ we get that 
\begin{align}\label{momentum_Sutskever}
    \vb m_{n+1}^S = \beta_n \vb m_{n}^S - \dfrac{\mu_n}{\mu_{n-1}}(M \mu_{n-1}-1) \vb* \alpha_{n-1} \odot \vb g_{n-1} - \vb* \alpha_n \odot \vb g_n.
\end{align}
We now concentrate on the shifted update of the Sutskever formulation, which, for every $n \geq 2$, it reads as follows:
\begin{align*}
    \vb z_{n+1}^S &= \vb z_n^S - \vb* \alpha_n \odot \vb g_n + \gamma_{n+1} \vb m_{n+1}^S +(\beta_n - \gamma_n) \vb m_n^S + \left[(\tilde{\gamma}_n - \mu_n) \vb* \omega_n \odot \vb g_{n-1} - \tilde{\gamma}_{n+1} \vb* \omega_{n+1} \odot \vb g_n \right] \\
    &= \vb z_n^S + \gamma_{n+1} \vb m_{n+1}^S + (\beta_n - \gamma_n ) \vb m_n^S + (\tilde{\gamma}_n - \mu_n ) \vb* \omega_n \odot \vb g_{n-1} - \left( \vb* \alpha_n + \tilde{\gamma}_{n+1} \vb* \omega_{n+1} \right) \odot \vb g_n.
\end{align*}
From the definitions of $\vb* \omega_n$ and $\vb* \omega_{n+1}$ (for every $n \geq 2$) it follows that 
\begin{align*}
    \vb z_{n+1}^S = \vb z_n^S + \gamma_{n+1} \vb m_{n+1}^S + (\beta_n - \gamma_n ) \vb m_n^S + (\tilde{\gamma}_n - \mu_n ) \dfrac{M\mu_{n-1}-1}{\mu_{n-1}} \vb* \alpha_{n-1} \odot \vb g_{n-1} - \left(1+\tilde{\gamma}_{n+1} \dfrac{M\mu_n - 1}{\mu_n}\right) \vb* \alpha_n \odot \vb g_n.
\end{align*}
Using \eqref{momentum_Sutskever} and performing some basic algebraic manipulations, we obtain for every $n \geq 2$ that
\begin{align}\label{shifted_Sutskever}
\vb z_{n+1}^S &= \vb z_n^S + \left[ \beta_n (1 + \gamma_{n+1}) - \gamma_n \right] \vb m_n^S - \dfrac{M\mu_{n-1}-1}{\mu_{n-1}} \left[ \mu_n(1+\gamma_{n+1}) - \tilde{\gamma}_n \right] \vb* \alpha_{n-1} \odot \vb g_{n-1} \nonumber \\
&- \left[ (1+\gamma_{n+1}) + \tilde{\gamma}_{n+1} \dfrac{M\mu_n-1}{\mu_n}\right] \vb* \alpha_n \odot \vb g_n.
\end{align}
Until we proceed further, we remind that $\vb z_n^S = \vb z_n$ is defined for each $n \geq 3$ in \eqref{zn_AAMMSU} with respect to the iterative process \eqref{AdaptiveAcceleratedMomentumMethodShiftedUpdates}. This is equivalent to saying that $\vb z_{n+1}^S$ is defined for every $n \geq 2$ for \eqref{AdaptiveAcceleratedMomentumMethodShiftedUpdates}. On the other, hand, for every $n \geq 3$, from \eqref{eq:9} it follows that $\vb a_{n-1} \odot \vb p_n = - (\vb* \theta_n - \vb* \theta_{n-1})$ (where we have used that $\vb a_{n-1} = \vb r_{n-1}$). But, utilizing that $\vb m_n^S = \vb* \theta_n - \vb* \theta_{n-1}$, we get the relationship between the original momentum term and the one in the Sutskever formulation, i.e. $\vb m_n^S = - \vb a_{n-1} \odot \vb p_n$, for every $n \geq 3$, hence $\vb m_{n+1}^S$ is defined for every $n \geq 2$.
Therefore, from the above computations, we have that the Sutskever formulation of \eqref{AdaptiveAcceleratedMomentumMethodShiftedUpdates} is defined for every $n \geq 2$ through \eqref{momentum_Sutskever} and \eqref{shifted_Sutskever}, hence the pair $(\vb m_{n+1}^S, \vb z_{n+1}^S)$ is defined for every $n \geq 2$. From \eqref{AdaptiveHeayBallMethod} we have that 
\begin{align*}
    \vb z_2^S = \vb z_2 &= \vb* \theta_2 + \gamma_2( \vb* \theta_2 - \vb* \theta_1 ) - \tilde{\gamma}_2 \vb* \omega_2 \odot \vb g_1 \\
    &= \vb* \theta_2 - \tilde{\gamma} \left( M - \dfrac{1}{\mu_1}\right) \vb* \alpha_1 \odot \vb g_1\\
    &= \vb* \theta_2 - \tilde{\gamma}(M-1) \vb* \alpha_1 \odot \vb g_1.
\end{align*}
Knowing that \eqref{AcceleratedGradientMethod} is equivalent to \eqref{AdaptiveHeayBallMethod}, for each $n \geq 1$, we have that $\vb* \theta_2 = \vb y_1 - \vb* \alpha_1 \odot \vb g_1$. According to the calculations presented above, it follows that 
\begin{align*}
    \vb z_2^S &= \vb y_1 - \left[ 1 + \tilde{\gamma}(M-1)\right] \vb* \alpha_1 \odot \vb g_1 \\
    &= (1- \mu_1) \vb* \theta_1 + \mu_1 \vb w_1 - \left[ 1+ \tilde{\gamma} (M-1) \right] \vb* \alpha_1 \odot \vb g_1 \\
    &= \vb w_1 - \left[ 1+ \tilde{\gamma} (M-1) \right] \vb* \alpha_1 \odot \vb g_1.
\end{align*}
On the other hand, for the first iteration we take the \eqref{AcceleratedGradientMethod} formulation, which implies that $\vb z_1^S = (1 - \tilde{\gamma}_1) \vb* \theta_1 + \tilde{\gamma}_1 \vb w_1 = \vb w_1$, where $\vb w_1$ is arbitrarily given. \\
Combining the calculations mentioned above, we obtain the algorithmic description of our optimizer shown below in Algorithm \eqref{AlgDescription_Sutskever_AAMMSU} entitled \emph{Sutskever formulation of AAMMSU}.

\begin{algorithm}[H]\label{AlgDescription_Sutskever_AAMMSU}
\SetAlgoLined
\SetKwFunction{FMain}{SutskeverAAMMSU}
\SetKwProg{Fn}{Function}{:}{}
  
\KwIn{ $\varepsilon \in (0, 1)$, $N \geq 1$, $\beta_2 \in (0,1)$, $\nu \in (0, 1)$, $\mu \in (0, 1)$, $M > 0$, $\tilde{\gamma} \in [\mu, 1)$, $\eta > 0$ }

\Fn{\FMain{$\varepsilon$, $N$, $\beta_2$, $\nu$, $\mu$, $M$, $\tilde{\gamma}$, $\eta$}}{
    \SetKwInput{kwInit}{Initialization}
    \kwInit{$\vb z_1^S$, $\vb m_2^S$, $\vb v_0 = \tilde{\vb v}_0 = 0$, $\tilde{\gamma}_1 = \mu_1 = 1$}
\For{$n = 1, 2, \ldots, N$}{
\If{$n = 2$}{
$\beta_n = \gamma_n = 0 \,$, $\, \tilde{\gamma}_n = \tilde{\gamma}\, $ and $\, \mu_n = \mu$\;
$\vb z_n^S = \vb z_{n-1}^S - \left[ 1 + \tilde{\gamma} (M - 1) \right] \vb* \alpha_{n-1} \odot \vb g_{n-1}$ \;
}
\If{$n \geq 3$}{
$\tilde{\gamma}_n = \tilde{\gamma}\, $\ and $\, \gamma_n = \tilde{\gamma} \dfrac{1-\mu}{ \mu }$ \;
$\vb m_{n}^S = \beta_{n-1} \vb m_{n-1}^S - \dfrac{\mu_{n-1}}{\mu_{n-2}}(M \mu_{n-2}-1) \vb* \alpha_{n-2} \odot \vb g_{n-2} - \vb* \alpha_{n-1} \odot \vb g_{n-1}$\;
$\vb z_{n}^S = \vb z_{n-1}^S + \left[ \beta_{n-1} (1 + \gamma_{n}) - \gamma_{n-1} \right] \vb m_{n-1}^S - \dfrac{M \mu_{n-2}-1}{\mu_{n-2}} \left[ \mu_{n-1} (1 + \gamma_{n}) - \tilde{\gamma}_{n-1} \right] \vb* \alpha_{n-2} \odot \vb g_{n-2} - \left[ (1 + \gamma_{n}) + \tilde{\gamma}_{n} \dfrac{M \mu_{n-1} - 1}{\mu_{n-1}}\right] \vb* \alpha_{n-1} \odot \vb g_{n-1}$\;
$\beta_n = 1 - \mu \,$ and $\, \mu_n = \mu$\;
}
$\vb g_n := \nabla F(\vb z_n^S) + \vb* \delta_n$\;
$\tilde{\vb v}_n = \beta_2 \tilde{\vb v}_{n-1} + (1-\beta_2) [\vb g_n]^2$\;
$\vb v_n = \max \lbrace \vb v_{n-1}, \tilde{\vb v}_n \rbrace$\;
$\vb a_n = \dfrac{\eta}{\varepsilon + \sqrt{\vb v_n}}$\;
$\vb* \alpha_n = \nu \vb a_n$\;
}
  \KwOut{ \footnotesize{$\vb z_{N}^S$}}
  } 
 \caption{Sutskever formulation of AAMMSU}
\end{algorithm}

We end the present section with a remark concerning the convergence of the adaptive algorithm \eqref{AlgDescription_Sutskever_AAMMSU}.
\begin{remark}
In the \texttt{PyTorch} formulation, namely the Sutskever form \eqref{AlgDescription_Sutskever_AAMMSU} of our AMSGrad-type algorithms, we have considered the main iteration to be $\vb z_{n}^S$, which actually represents $\vb z_n$ of \eqref{AdaptiveAcceleratedMomentumMethodShiftedUpdates}. This is consistent with the findings from the Theorem \eqref{Theorem_AsymptoticBehavior} where the rate of convergence for \eqref{AdaptiveAcceleratedMomentumMethodShiftedUpdates} is given for the squared norm of the gradient of the objective function in the shifted updates $\vb z_n = \vb z_n^S$.
\end{remark}

\section{Neural network training}\label{Section_NeuralNetworks}
\begin{figure}[htb]
    \centering
\begin{subfigure}{0.3\textwidth}
  \includegraphics[scale=0.3]{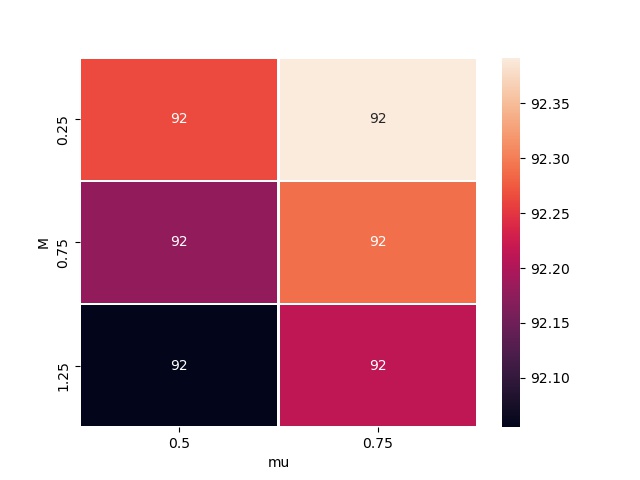}
\end{subfigure}\hfil\hfil\hfil
\begin{subfigure}{0.3\textwidth}
  \includegraphics[scale=0.3]{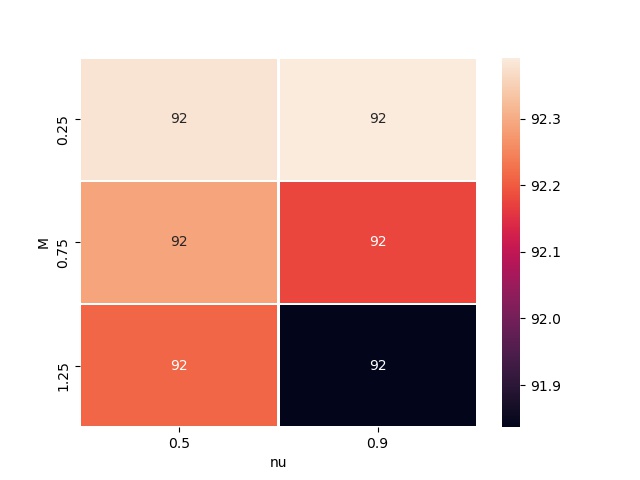}
\end{subfigure}\hfil\hfil\hfil
\begin{subfigure}{0.3\textwidth}
  \includegraphics[scale=0.3]{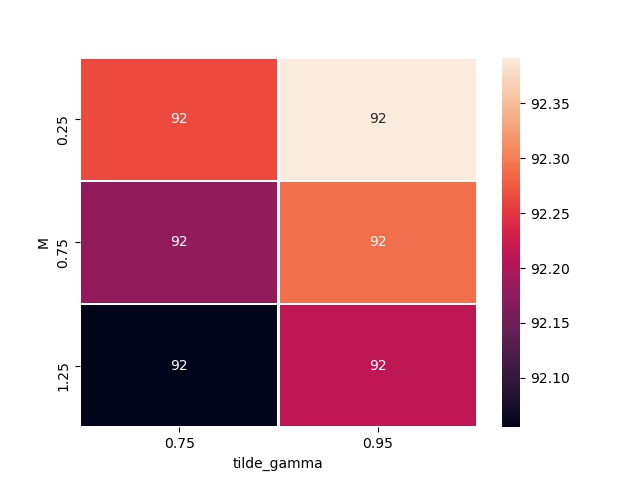}
\end{subfigure}

\begin{subfigure}{0.3\textwidth}
  \includegraphics[scale=0.3]{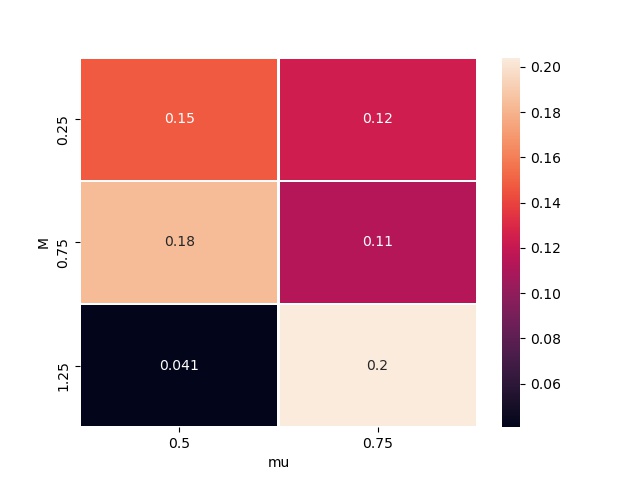}
\end{subfigure}\hfil \hfil\hfil
\begin{subfigure}{0.3\textwidth}
  \includegraphics[scale=0.3]{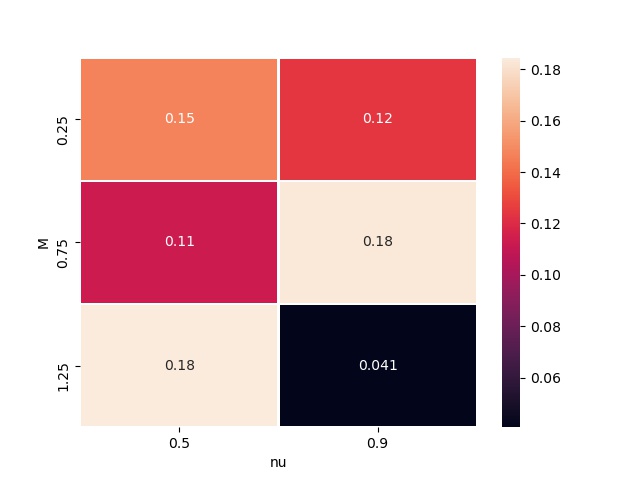}
\end{subfigure}\hfil\hfil\hfil
\begin{subfigure}{0.3\textwidth}
  \includegraphics[scale=0.3]{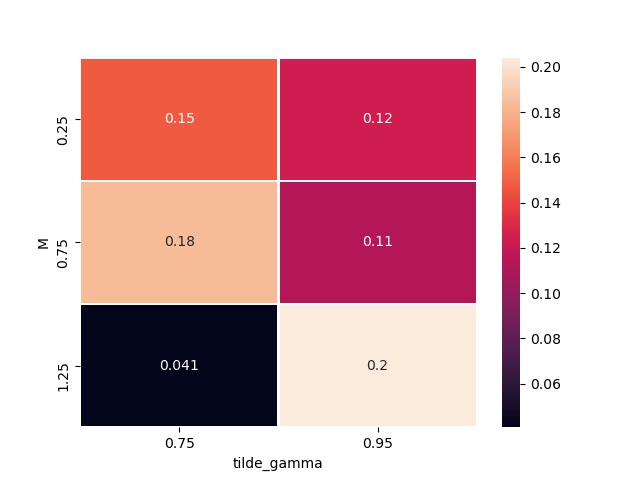}
\end{subfigure}
\caption{Heatmaps for \texttt{LR-MNIST} and \eqref{AdaptiveAcceleratedMomentumMethodShiftedUpdates}}
\label{fig:LR-MNIST_heatmaps_AAMMSU}
\end{figure}

\begin{figure}[htpb]
    \centering
\begin{subfigure}{0.3\textwidth}
  \includegraphics[scale=0.3]{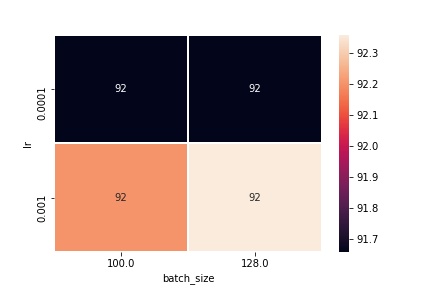}
\end{subfigure}
\begin{subfigure}{0.3\textwidth}
  \includegraphics[scale=0.3]{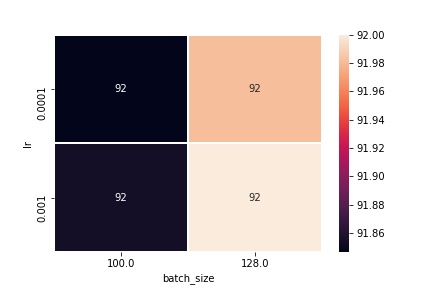}
\end{subfigure}

\begin{subfigure}{0.3\textwidth}
  \includegraphics[scale=0.3]{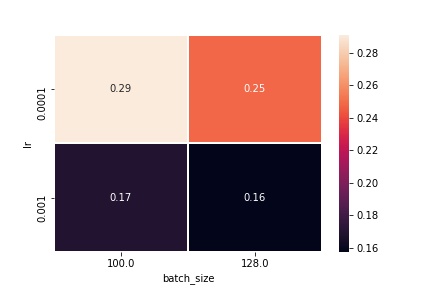}
\end{subfigure}
\begin{subfigure}{0.3\textwidth}
  \includegraphics[scale=0.3]{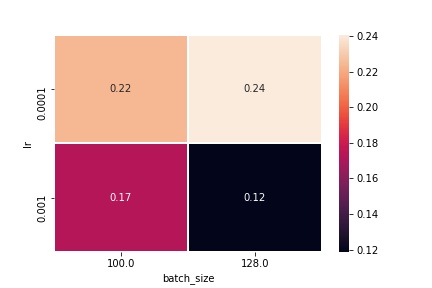}
\end{subfigure}
\caption{Batch size - learning rate evolution for \texttt{LR-MNIST} }
\label{fig:LR-MNIST_batch_size_evolution}
\end{figure}

Various numerical simulations relating to the training of neural networks are carried out in this section. Due to limited available resources, we have considered only two types of classical benchmark datasets, i.e. \emph{MNIST} and \emph{CIFAR10}, along with four types of neural networks, namely \emph{Logistic Regression} denoted as \emph{LR} (although in the well-known statistical framework it is a regression algorithm, we consider using the classifier which is put on top of the regression method through the thresholding of the probabilities, since this approach is usually employed in Machine Learning), a convolutional model (\emph{CNN}), \emph{VGG-11} and \emph{ResNet18}, respectively. In order to simplify our presentation, we shall use the following short names for the combination of models and datasets: \texttt{LR-MNIST}, \texttt{CNN-CIFAR10}, \texttt{VGG-CIFAR10} and \texttt{ResNet-CIFAR10}. For the implementation of the aforementioned models \footnote{\url{https://github.com/CDAlecsa/AAMMSU}}, we have used the \emph{CrossEntropy} loss function for multi-classification. For our experiments we have compared the \eqref{AdaptiveAcceleratedMomentumMethodShiftedUpdates} given through Algorithm \eqref{AlgDescription_Sutskever_AAMMSU} entitled \emph{Sutskever formulation of AAMMSU} with the \emph{AMSGrad} optimizer, since they are closely related through the choice of the accumulated squared gradients. For both of them we have set $\beta_2 = 0.999$ and $\varepsilon = 1\text{e-}8$, and in most cases we have considered as the default learning rates the values $1\text{e-}4$ and $1\text{e-}3$ since these are the ones that are used for Adam-type optimizers (this is natural, because from a theoretical perspective, for most choices of the number of epochs $N$, $1\text{e-}4$ and $1\text{e-}3$ are close to $1/\sqrt{N}$). The full results of our experiments (where we have printed the mean and the standard deviation) are deferred to the subsection \eqref{subsection_additional_numerical_results} of the Appendix section \eqref{Section_Appendix}, where the blue color from the tables represents the best validation accuracy and the red color is used for the best test accuracy. Along with these tables, we will present, for the each chosen model, some heatmaps regarding the evolution of the validation accuracy. We made some numerical simulations for various values of the optimizer's coefficients, and for the pairs of parameters to be presented in a heatmap, we have removed the duplicates. More precisely, we have chosen the values of the parameters which were not present in the heatmap and which led to the highest validation mean along with the lowest standard deviation, for a specific type of choices for the parameters which belong to the heatmap. It is noteworthy that we did not set in our experiments a random seed as in many neural network experiments, and instead of this, we have taken a more realistic approach, by repeating the same experiments for a predefined number of runs (briefly \texttt{n$\_$runs}), in order to retain the mean and the standard deviation of the results. In each experiment and for every simulation run, we have generated the \texttt{PyTorch} dataset splitted into training $\&$ validation ($80\%$, $20\%$), along with the corresponding \texttt{PyTorch} dataset. It should be noted that, since we have not used a random seed approach, every time we have generated a new experiment, we have obtained different training/validation datasets but with the same established percentages $80\%$ and $20\%$, respectively (we have considered this process to be more reasonable from a practical point of view). Until we begin explaining in detail our numerical results, we note that we have examined the comparison between \eqref{AdaptiveAcceleratedMomentumMethodShiftedUpdates} and \emph{AMSGrad} through the lens of the baseline metrics, namely the \emph{loss} and the \emph{accuracy} of the models (as in the majority of research papers in Machine Learning). Even though this choice differs from our theoretical results from Theorem \eqref{Theorem_AsymptoticBehavior} and from \eqref{OrderConv_Bound} where we have established bounds for the gradient of the objective function, we have considered that the usage of \emph{loss} and \emph{accuracy} metrics are more relevant for practical applications. We also remind that the results from the grid search of the parameters are different than the ones from the evolution of the learning versus the batch sizes (see the tables from subsection \eqref{subsection_additional_numerical_results}), since for the latter case we have re-run the models with the best evaluation coefficients. \\
Our initial findings include the \texttt{LR-MNIST} model where we have employed $5$ runs per experiment, and which has been chosen by us since it is composed of a neural network which can be trained fast along with the basic \emph{MNIST} dataset. Furthermore, \texttt{LR-MNIST} was mainly considered in order to make a grid search for the parameters $M$, $\mu$, $\nu$ and $\tilde{\gamma}$. The complete results for \eqref{AdaptiveAcceleratedMomentumMethodShiftedUpdates} can be found in the table \eqref{table_LR_MNIST_AAMMSU}, while the results for the \emph{AMSGrad} algorithm are presented in the table \eqref{table_models_lr_cnn_AMSGrad}. Since from our preliminary results we have observed that the \texttt{LR-MNIST} model can be trained quite fast and can eventually lead to overfitting, we have considered the number of epochs in $\lbrace 15, 35, 50 \rbrace$. The range of all the parameters and the best validation $\&$ test \emph{accuracy} results are reported in table \eqref{table_LR_MNIST_AAMMSU}. Our initial conclusion is that the accuracy values are stable with respect to the parameter $M$, while the inertial coefficient $\tilde{\gamma}$ must take a value closer to $1$ in order to obtain a good mean result (but not necessarily the lowest standard deviation value). This is also due to the fact that $\tilde{\gamma}$ must be greater than $\mu$. A simplification of the results on the validation datasets from table \eqref{table_LR_MNIST_AAMMSU} is presented in the figure \eqref{fig:LR-MNIST_heatmaps_AAMMSU}, where the mean values are represented in the first row, while the second row contains the standard deviation results.
It is easy to observe for a large value of $M$, i.e. $1.25$ in our case, the standard deviation of the results is, in general, lower than the one from the other coefficient choices (when $\mu$ and $\tilde{\gamma}$ have values closer to $0.5$). In this case when $M = 1.25$, if $\tilde{\gamma}$ is set to a lower value, namely $0.75$, the standard deviation gets lower along with also a lower mean value of the validation accuracy. On the other hand, we can assert that the results of the accuracy when $\nu = 0.5$ are stable enough, while the values of $\mu$ are much more diverse, in the sense that they appear to be related to the values of $M$. \\
For the evolution of the learning rate with respect to the batch size, as in table \eqref{table_models_AAMMSU_AMSGrad_batch_evolution}, we present the heatmap \eqref{fig:LR-MNIST_batch_size_evolution} for the validation accuracy, where \eqref{AdaptiveAcceleratedMomentumMethodShiftedUpdates} is presented in the first column, while the \emph{AMSGrad} results are showed in the second one. The \emph{AMSGrad} algorithm seems more stable when the batch size is $128$, than \eqref{AdaptiveAcceleratedMomentumMethodShiftedUpdates}, but the latter one achieves a better validation accuracy when we choose a higher learning rate $1e-3$. The standard deviation of the results, on the other hand, indicates that both algorithms are fairly stable in terms of simulation runs.

\begin{figure}[ht!]
\hspace{-4em}
    \centering
\begin{subfigure}{0.28\textwidth}
  \includegraphics[scale=0.28]{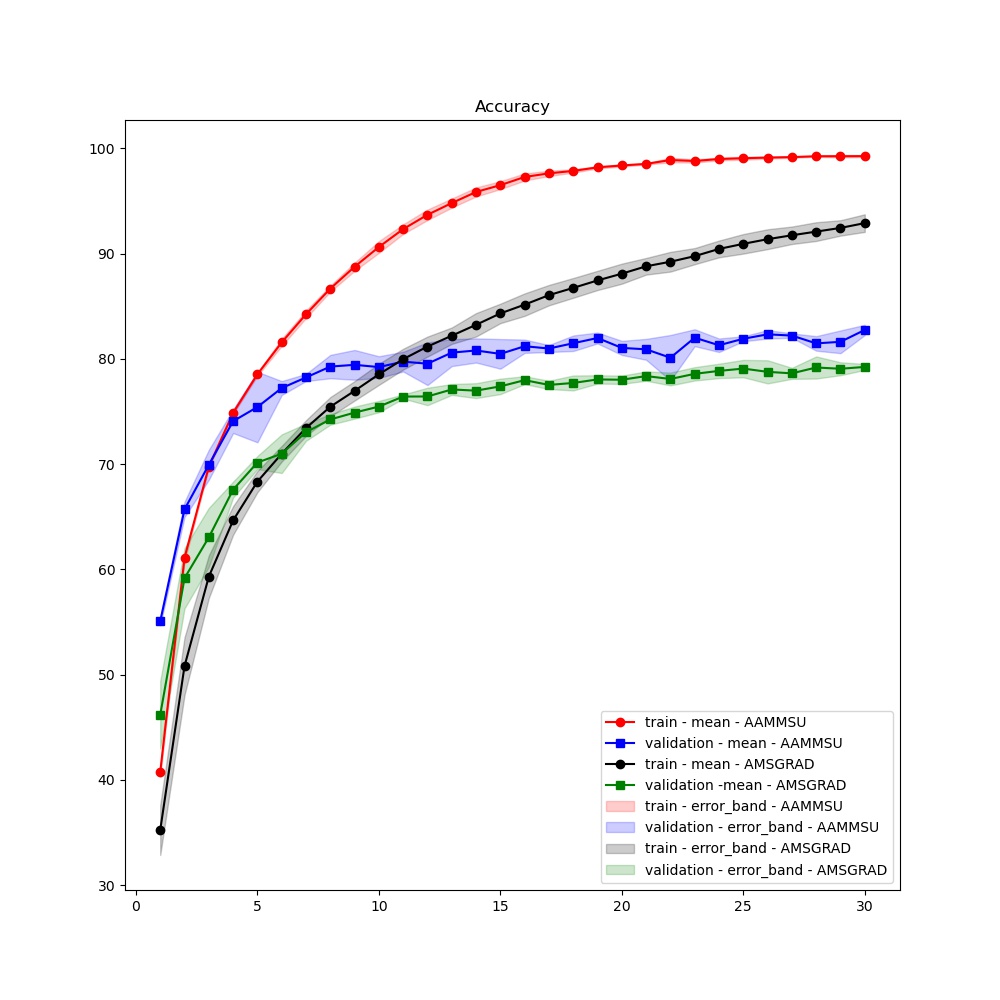}
\end{subfigure}\hfil
\begin{subfigure}{0.28\textwidth}
  \includegraphics[scale=0.28]{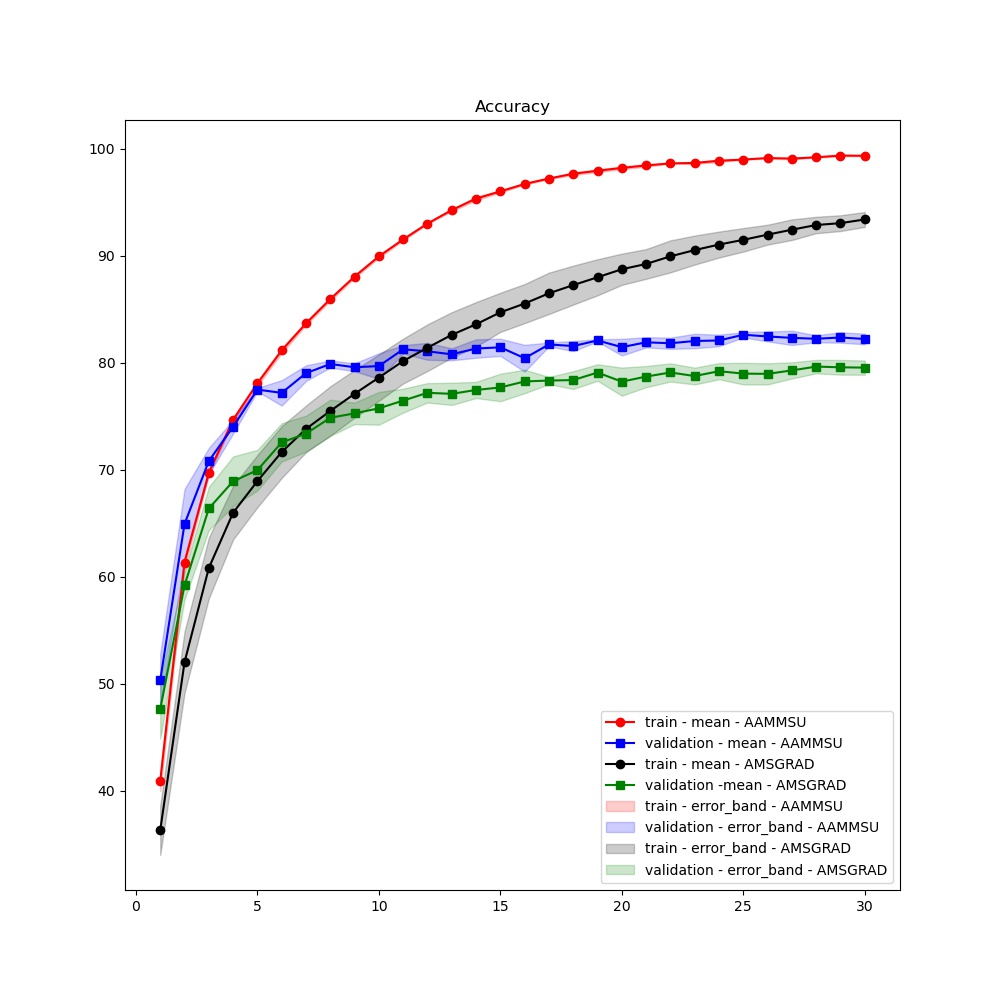}
\end{subfigure}

\hspace{-4em}
\begin{subfigure}{0.28\textwidth}
  \includegraphics[scale=0.28]{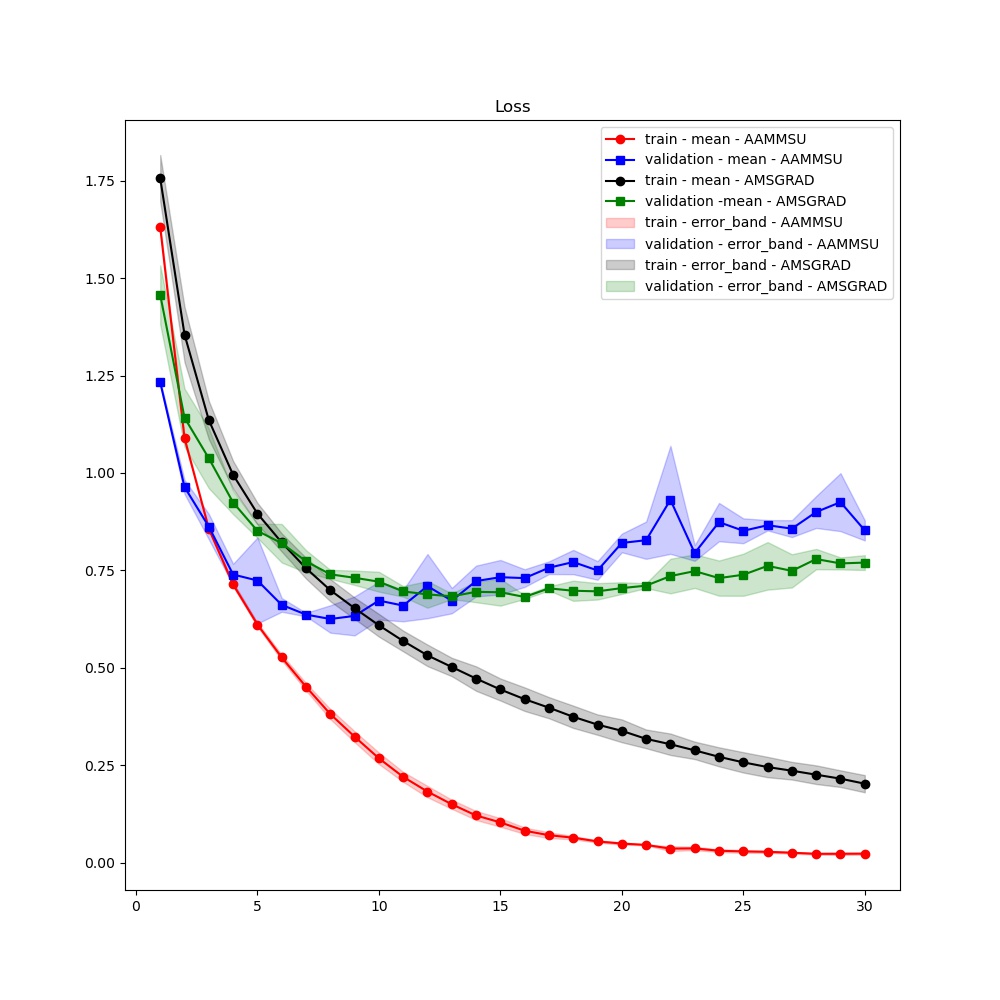}
\end{subfigure}\hfil
\begin{subfigure}{0.28\textwidth}
  \includegraphics[scale=0.28]{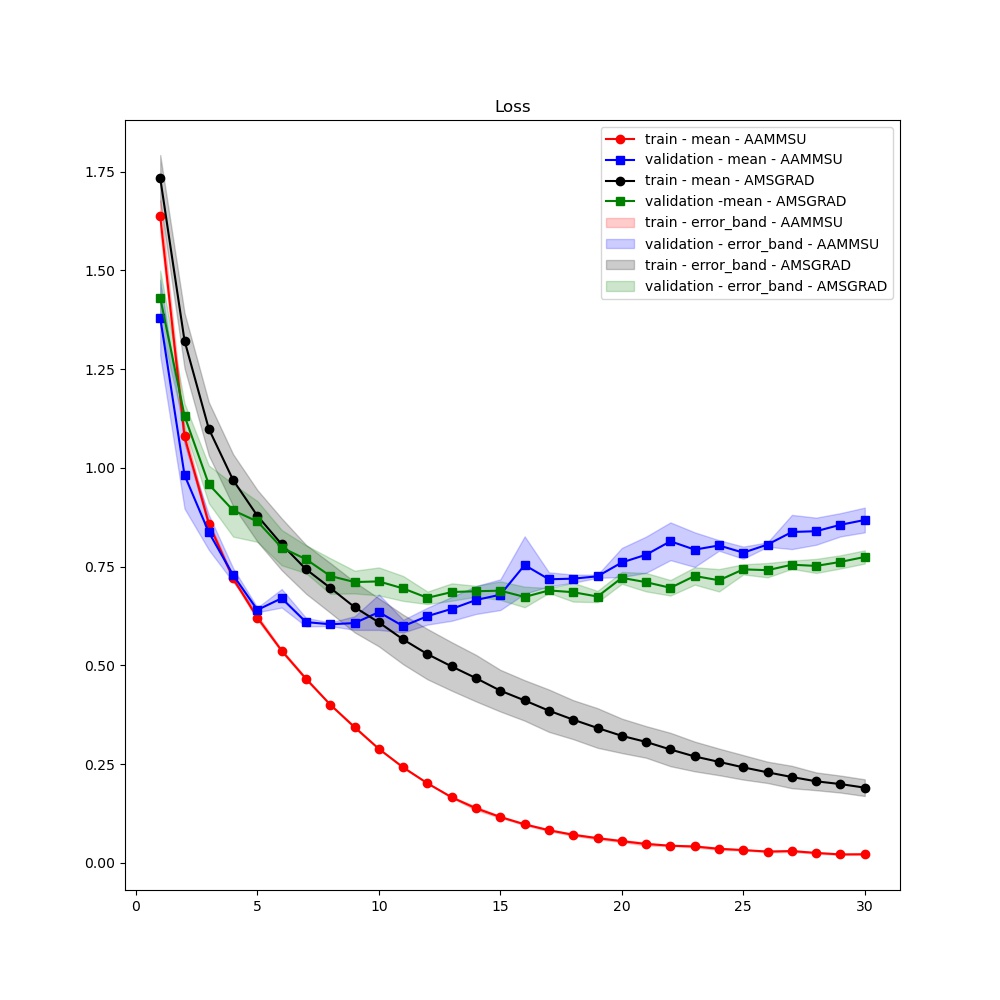}
\end{subfigure}
\caption{Convergence profiles for \texttt{CNN-CIFAR10}}
\label{fig:CNN_CIFAR10_loss_accuracy}
\end{figure}

Now, we turn our attention to the \texttt{CNN-CIFAR10} model \footnote{\url{https://shonit2096.medium.com/cnn-on-cifar10-data-set-using-pytorch-34be87e09844}} where we have employed $3$ runs per experiment, where we have made a grid search for epochs in $ \lbrace 10, 17, 30 \rbrace$ (the limitation of the model with respect to the number of epochs was studied by us in some preliminary numerical results) as reported in table \eqref{table_CNN_CIFAR10_AAMMSU}. We observe that different values of the pairs $(M, \mu)$, $(M, \nu)$ and $(M, \tilde{\gamma})$ have a higher impact on the convolutional model than for the previous model, which reveals that the \texttt{CNN-CIFAR10} model is much less robust than the \texttt{LR-MNIST} model for our adaptive optimizers. It is critical to emphasize that by following table \eqref{table_CNN_CIFAR10_AAMMSU} one observes two cases. When $M$ is close to $0$ and the elements from the pair $(\mu, \nu, \tilde{\gamma})$ are close to $0.5$ or $0.75$, then one obtains good values for the validation $\&$ test \emph{accuracy}. On the other hand, when $M$ is relatively large, i.e. $M = 1.75$, then our optimizer converges really slow and requires much more epochs. Now, let's take a look at figure \eqref{fig:CNN_CIFAR10_loss_accuracy} representing the \emph{loss} and \emph{accuracy} values for both adaptive optimizers, where the validation results are given in the first column, while the test results are given in the second one. The colors red and black are used for the training metrics, while the blue and green are used for the evaluation metrics (validation $\&$ test). Here, the red and blue colors are utilized for \eqref{AdaptiveAcceleratedMomentumMethodShiftedUpdates} while blue and green are used for \emph{AMSGrad}. We note that the \emph{accuracy} values are higher than those of the \emph{AMSGrad}, despite the fact that our adaptive algorithm overfits more in terms of the \emph{loss} metric.

\begin{figure}[ht!]
    \centering
\begin{subfigure}{0.3\textwidth}
  \includegraphics[scale=0.3]{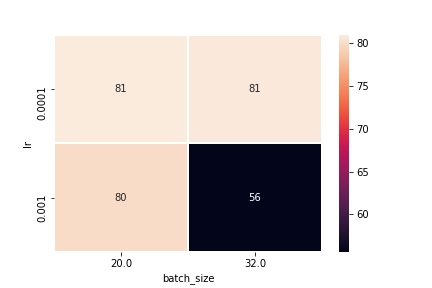}
\end{subfigure}
\begin{subfigure}{0.3\textwidth}
  \includegraphics[scale=0.3]{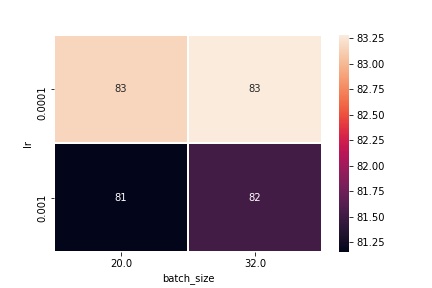}
\end{subfigure}

\begin{subfigure}{0.3\textwidth}
  \includegraphics[scale=0.3]{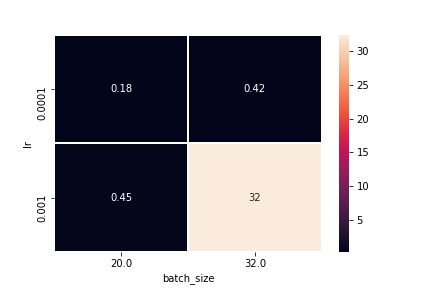}
\end{subfigure}
\begin{subfigure}{0.3\textwidth}
  \includegraphics[scale=0.3]{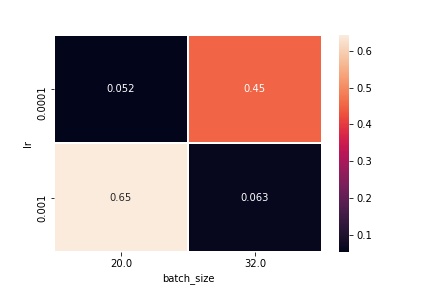}
\end{subfigure}
\caption{Batch size - learning rate evolution for \texttt{CNN-CIFAR10}}
\label{fig:CNN-CIFAR10_batch_size_evolution}
\end{figure}
\begin{figure}[hbt!]
    \centering
\hspace{-4em}
\begin{subfigure}{0.28\textwidth}
  \includegraphics[scale=0.28]{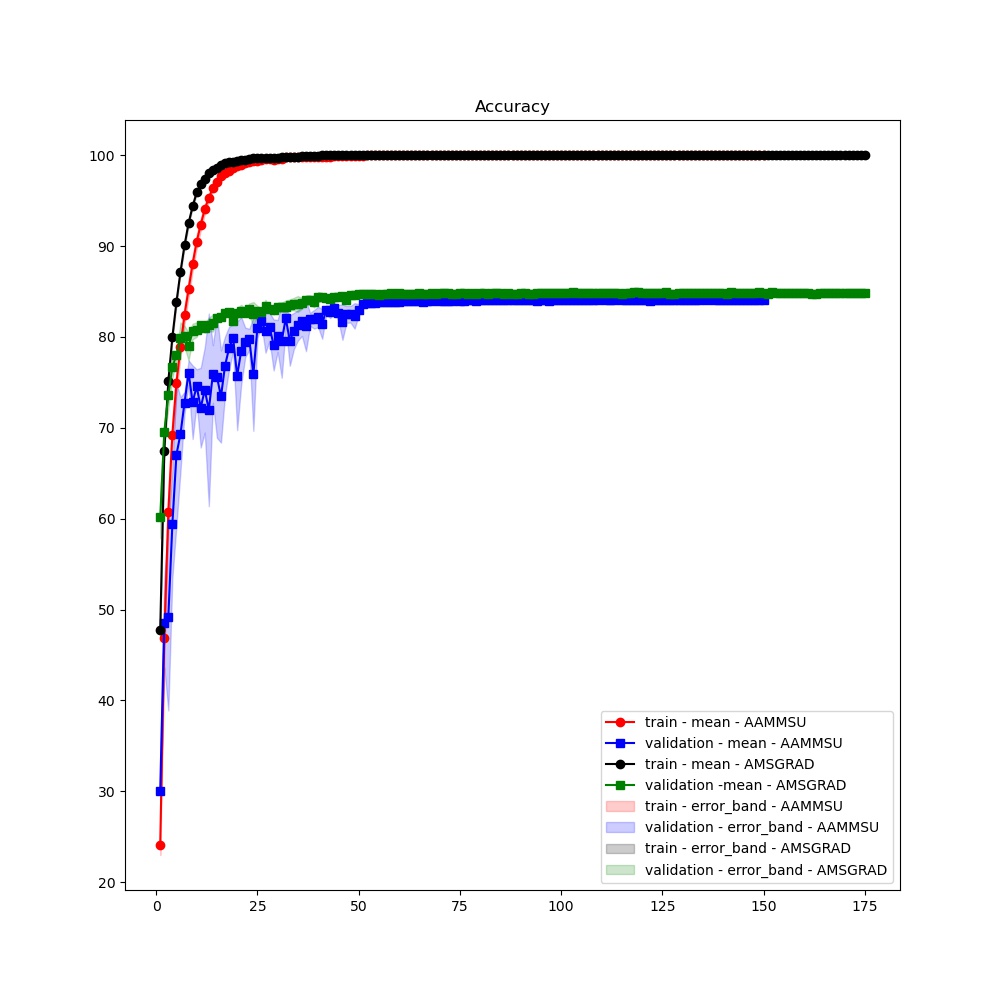}
\end{subfigure}\hfil
\begin{subfigure}{0.28\textwidth}
  \includegraphics[scale=0.28]{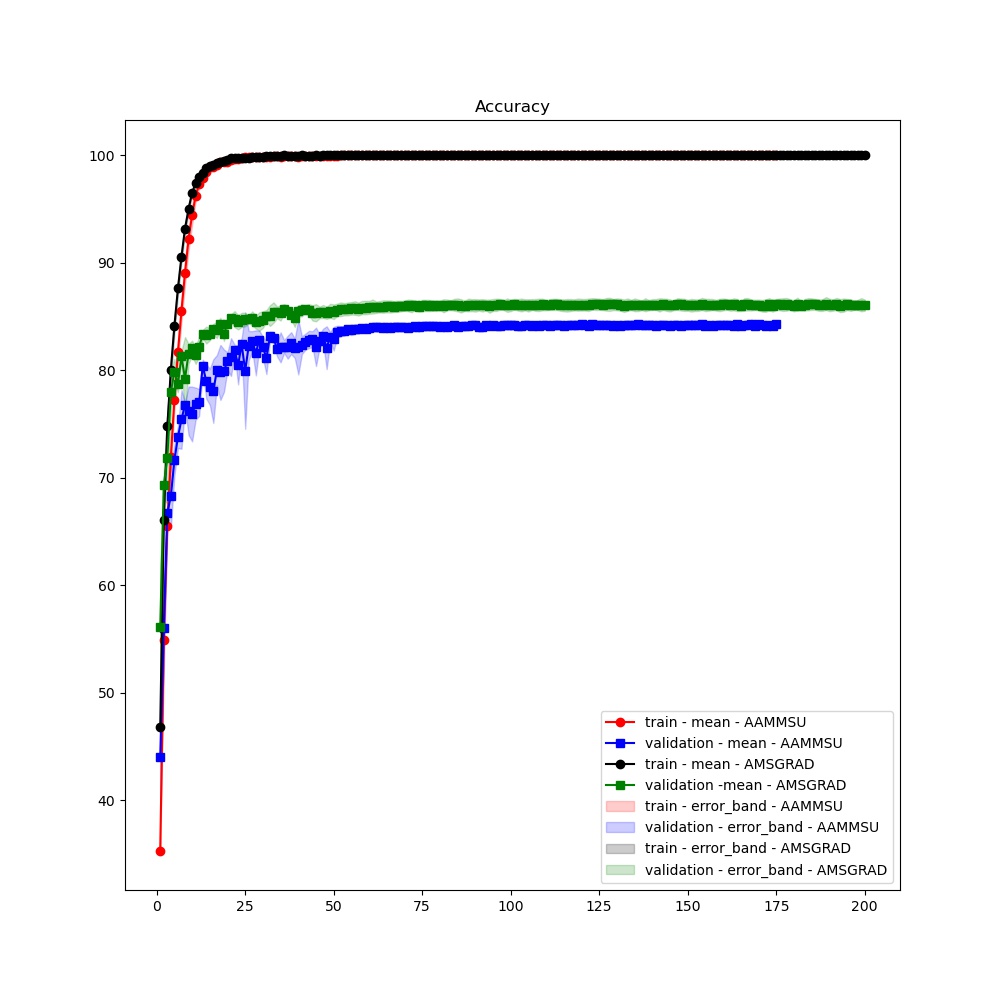}
\end{subfigure}

\hspace{-4em}
\begin{subfigure}{0.28\textwidth}
  \includegraphics[scale=0.28]{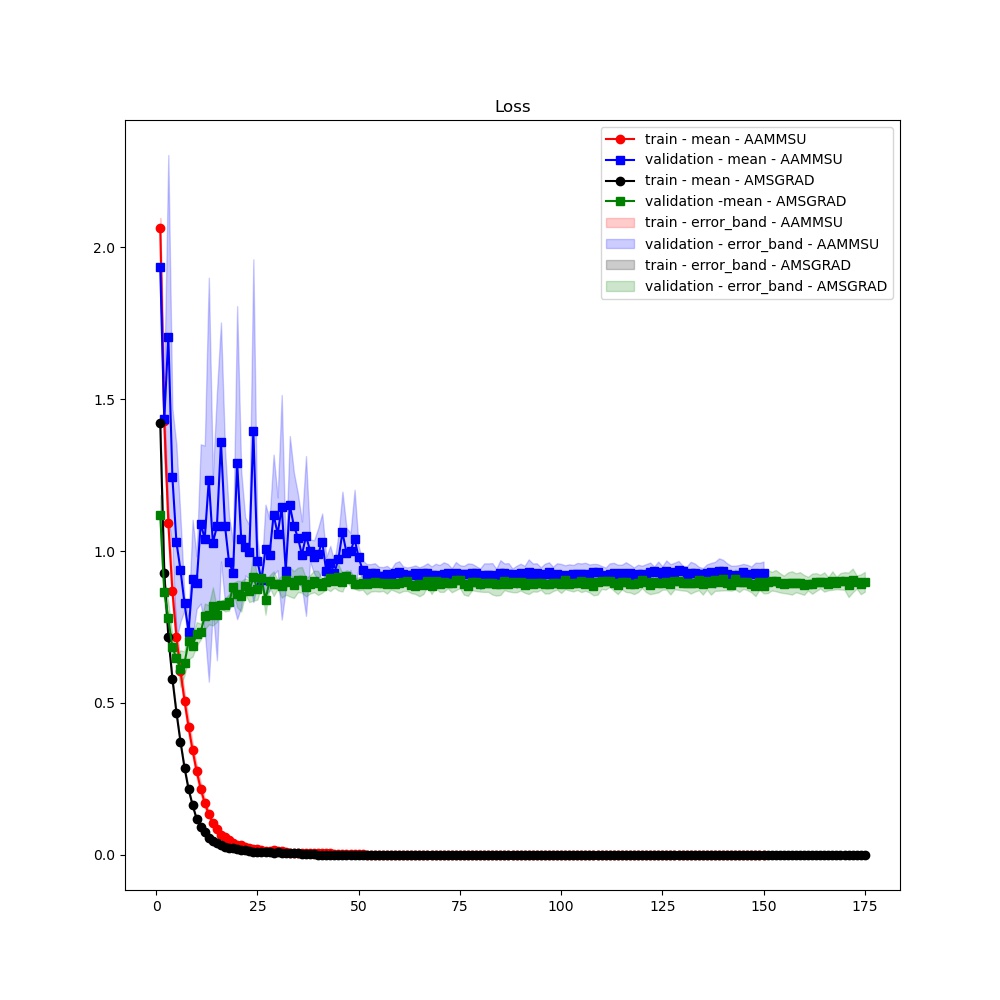}
\end{subfigure}\hfil
\begin{subfigure}{0.28\textwidth}
  \includegraphics[scale=0.28]{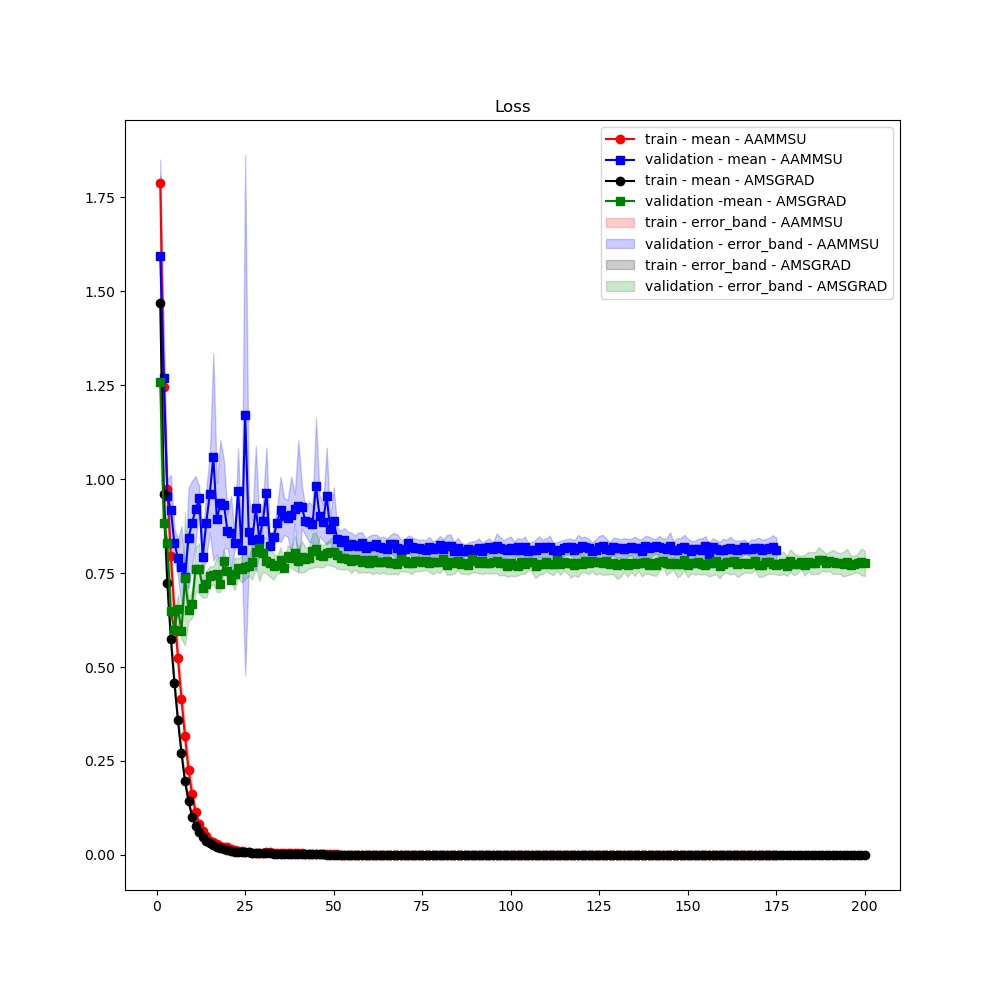}
\end{subfigure}
\caption{Convergence profiles for \texttt{VGG-CIFAR10} and \texttt{ResNet-CIFAR10}}
\label{fig:VGG_RESNET_loss_accuracy}
\end{figure}

Last, for the evolution of the learning rates over batches for the model \texttt{CNN-CIFAR10} with respect to the validation accuracy, as in table \eqref{table_models_AAMMSU_AMSGrad_batch_evolution}, we have considered the figure \eqref{fig:CNN-CIFAR10_batch_size_evolution}. The number of epochs was set in this case to $30$, while $(M, \mu, \nu, \tilde{\gamma})$ was chosen $(0.75, 0.5, 0.5, 0.75)$. The \eqref{AdaptiveAcceleratedMomentumMethodShiftedUpdates} is presented in the first column, while \emph{AMSGrad} in the second one. In the top row we have the means over $3$ runs where we see that our adaptive algorithm achieves smaller validation accuracy for different batch sizes and learning rate values, while in the bottom row we have the reported standard deviations over those $3$ runs, where we see that \eqref{AdaptiveAcceleratedMomentumMethodShiftedUpdates} is much more stable over a set of given simulation runs (this must be related also the small number of simulation runs due to our limited available resources). The only exception is when the batch size \texttt{bs}=$32$ and the learning rate is $1\text{e-}3$, in which situation we obtain a very large standard deviation (this case is also reflected in the top row where the reported mean is very small compared to the other cases). \\
At last, we present the results regarding the models \texttt{VGG-CIFAR10} and \texttt{ResNet-CIFAR10}, respectively \footnote{\url{https://github.com/uclaml/Padam/tree/master/models}}, where the number of simulation runs per experiment was set to $5$. Tables \eqref{table_VGG_ResNet_AAMMSU} and \eqref{table_grid_search_vgg_resnet_AMSGrad} show the complete results of the parameter search, where we have chosen epochs in $\lbrace 50, 75, 100, 150, 175, 200 \rbrace$. For these models we did not used a constant learning rate as in the previous experiments, but we have considered a \texttt{PyTorch} \emph{multi-step learning rate scheduler} similar to \cite{Padam}, namely we have decreased the learning rate of both adaptive algorithms with $1\text{e-1}$ (it was initially set with the value $1\text{e-}3$) at the epochs $50$, $100$ and $150$. At the same time, we have set $M$ in $\lbrace 1, 2 \rbrace$ while $(\mu, \nu, \tilde{\gamma}) = (0.5, 0.5, 0.5)$. In the plots from figure \eqref{fig:VGG_RESNET_loss_accuracy}, the black and the red colors represent the training metrics, while the green and blue represent the validation metrics. At the same time, for \eqref{AdaptiveAcceleratedMomentumMethodShiftedUpdates} we have employed the colors red and blue. It is interest to note that for the best validation parameters of our adaptive algorithms for both models, the \eqref{AdaptiveAcceleratedMomentumMethodShiftedUpdates} optimizer overfits in a similar manner as \emph{AMSGrad} but requires less epochs than the latter one. This benefit can also be used in conjunction with an additional early stopping technique. Moreover, from the same figure and for both models, we infer that although the standard deviation of the results is very high in the first epochs, it stabilizes quickly to the one of the \emph{AMSGrad} algorithm. Finally, similar to the case of the \texttt{CNN-CIFAR10} model, when $M$ has a large value (in our case $M = 2$), then we have a large standard deviation along with a very slow convergence that requires a large predefined number of epochs.

\section{Conclusions}\label{Conclusions}
In this final section we present a brief overview on our adaptive methods along with the limitations and the possible extensions for future research.
\subsection{Findings}
In the present study, our primary contribution is the introduction of an \emph{AMSGrad}-type approach with shifted updates, which is completely novel from both a theoretical and practical standpoint. Moreover, our algorithms which we have introduced can be traced back to the Nesterov-type algorithm with two inertial steps from \cite{ALP} and the flexibility of the inertial methods from \cite{ALV}. Second of all, we have shown that the changes of variables relating \emph{AMSGrad}-type algorithms to heavy-ball methods gives us a consistent theoretical treatment, hence our \emph{AMSGrad}-like methods can be faithfully seen as adaptive Nesterov type extensions with two momentum terms, i.e. $\tilde{\vb* \gamma}_n$ (that depends on $\vb* \gamma_n$) and $\mu_n$. More specifically, we have justified that these algorithms maintain the rate of convergence of popular adaptive approaches in non-convex stochastic settings. But, if we make a comparison with the theoretical result of Chen et al. from \cite{ChenAdam}, the bound of $\mathbb{E} \left[ \min\limits_{n = 1, \ldots, N} \| \nabla F(\vb z_n) \|^2 \right]$ given in \eqref{OrderConv_Bound} is tighter than the one from \cite{ChenAdam} which contain extra terms depending on the $l^2$ and $l^1$ norms of the difference of effective stepsize values, i.e. $\| \vb a_n - \vb a_{n-1} \|^2$.

\subsection{Research limitations}
\begin{itemize}
    \item The aforementioned work of Chen et al. \cite{ChenAdam} contains a different adaptive method variant called \emph{AdaForm} that we did not pursue studying in relation to our algorithms with shifted updates because it goes beyond the theoretical and numerical treatments from our paper. Nevertheless, it is worthwhile to take it into consideration in a subsequent research paper.
    \item In contrast with the \emph{AMSGrad} optimizer, the \emph{AAMMSU} method given in Algorithm \eqref{AlgDescription_Sutskever_AAMMSU} entitled \emph{Sutskever formulation of AAMMSU} contains, beside $\beta_2$, $\varepsilon$ and $\nu$ (where $\nu$ is in fact often denoted by $\beta_1$ in the \emph{AMSGrad} formulation), $3$ hyper-parameters, namely $M$, $\mu$ and $\tilde{\vb* \gamma}$. This actually means that the grid search for tuning these terms is more time consuming when comparing with the hyper-parameters of other adaptive optimizers. Despite the fact that the presence of shifted updates is new from a theoretical perspective for adaptive optimizers, the numerical simulations from section \eqref{Section_NeuralNetworks} reveal that the empirical performance in the neural network training is competitive but does not differ significantly from the Adam-type optimizers.
\end{itemize}

\subsection{Recommendations for future research}
\begin{itemize}
    \item It is worth pointing out that in the present paper the sequence $(\mu_n)_{n \geq 1}$ is limited to be from $\mathbb{R}$ and not from $\mathbb{R}^d$. The reason behind this is that the inequality that we have applied and that which has led to \eqref{ineq_norm_theta_w} is the classical Jensen inequality, which was employed in \cite{GhadimiLanAccelerated}. We highlight that a future research investigation can cover the more general case when $(\mu_n)_{n \geq 1} \in \mathbb{R}^d$, where one can create a multi-dimensional type Jensen inequality via Hadamard multiplications.
    \item Finally, we also like to add that there are other possible extensions which can be made to our adaptive methods. As an example, one can extend our analysis to non-convex and possibly non-smooth objective functions as in the recent work of Xiao \cite{XiaoAdam} or to investigate convergence guarantees of adaptive stochastic algorithms using the KL property (which stands for Kurdyka-\L{}ojasiewicz) as in the paper of Chouzenoux et al. \cite{ChouzenouxKL} (see also \cite{ALV} for the deterministic case of some accelerated optimization algorithms endowed with backward inertial steps with respect to the KL assumption).
\end{itemize}

\section{Acknowledgements}\label{Section_Acknowledgements}
We would like to thank J\' er\^ ome Bolte for pointing out several recent articles to us that are closely related to the connection between Machine Learning and Optimization problems, as well as anonymous reviewers who contributed to the development the current paper.

\section{Declarations}
\underline{Ethical approval}:\\
Not applicable.
\vskip+0.15cm
\underline{Competing interests}:\\
The author declares that he has no financial interests.
\vskip+0.15cm
\underline{Funding}:\\
This work was supported by a grant of the Romanian Ministry of Research and Innovation, CCCDI - UEFISCDI, project number 178PCE/2021, PN-III-P4-ID-PCE-2020-0788, Object PErception and Reconstruction with deep neural Architectures (OPERA), within PNCDI III. 
\vskip+0.15cm
\underline{Availability of data and materials}:\\
The datasets used in this paper can be accessed through the open source \texttt{PyTorch} module.
\vskip+0.15cm
\underline{Authors' contributions}:\\
Not applicable.

\section{Data availability statement}
Data sharing not applicable to this article as all the datasets used were generated using the \texttt{Pytorch} module.

\bibliographystyle{plain}  
\bibliography{references}

\section{Appendix}\label{Section_Appendix}

\subsection{Proof of Proposition \eqref{Proposition_SquaredVector}}

It is obvious that the inequality $\| [\vb u]^2 \| \leq \| \vb u \|^2$ is equivalent to $\| [\vb u]^2 \|^2 \leq \| \vb u \|^4 = \left( \| \vb u \|^2 \right)^2$, due to the non-negativity of the terms. Using the definition of the Euclidean inner product with respect to $\| \cdot \|_2$, the inequality becomes $\langle [\vb u]^2, [\vb u]^2 \rangle \leq \left( \langle \vb u, \vb u \rangle \right)^2$, which is identical to $\sum\limits_{j = 1}^{d} u_{[j]}^2 u_{[j]}^2 \leq \left( \sum\limits_{j = 1}^{d} u_{[j]} u_{[j]} \right)^2$, which is in fact $\sum\limits_{j = 1}^{d} (u_{[j]}^2)^2 \leq \left( \sum\limits_{j = 1}^{d} u_{[j]}^2 \right)^2$. In order to show the previous inner product inequality, we denote $v_j := u_{[j]}^2$, hence $v_j \geq 0$ for each $j \in \lbrace 1, \ldots, d \rbrace$. This implies that we must prove $\sum\limits_{j = 1}^{d} v_j^2 \leq \left( \sum\limits_{j = 1}^{d} v_j \right)^2$. This can be easily shown by induction as follows. If $d = 1$, then the conclusion is obvious. On the other hand, if $d > 1$, we consider to be true the inequality $\sum\limits_{j = 1}^{d} v_j^2 \leq \left( \sum\limits_{j = 1}^{d} v_j \right)^2$. Then $\sum\limits_{j = 1}^{d + 1} v_j^2 = \sum\limits_{j = 1}^{d} v_j^2 + v_{d + 1}^2 \leq \left( \sum\limits_{j = 1}^{d} v_j \right)^2 + v_{d + 1}^2 = (v_1 + \ldots + v_d)^2 + v_{d + 1}^2 = (v_1 + \ldots + v_{d+1})^2 - 2 v_{d+1} \cdot (v_1 + \ldots + v_d) \leq (v_1 + \ldots + v_{d+1})^2 = \left( \sum\limits_{j = 1}^{d + 1} v_j \right)^2$, where we have used the fact that each $v_j$ is non-negative.

\subsection{Proof of Proposition \eqref{Proposition_HadamardProduct}}

Since the norms of the vectors take non-negative values, then it is enough to show that $\| \vb u \odot \vb v \|^2 \leq \| \vb u \|^2 \cdot \| \vb v \|^2$, which takes the form $\langle \vb u \odot \vb v, \vb u \odot \vb v \rangle \leq \langle \vb u, \vb u \rangle \cdot \langle \vb v, \vb v \rangle$. Since we are working with the Euclidean norm $\| \cdot \|_2$, we must prove the equivalent inequality $\sum\limits_{j = 1}^{d} u_{[j]}^2 v_{[j]}^2 \leq \left( \sum\limits_{j = 1}^{d} u_{[j]}^2 \right) \cdot \left( \sum\limits_{j = 1}^{d} v_{[j]}^2 \right)$. By using the CBS inequality, we find that
\begin{align*}
\sum\limits_{j = 1}^{d} u_{[j]}^2 v_{[j]}^2 &\leq \sqrt{\left( \sum\limits_{j = 1}^{d} u_{[j]}^4 \right) \cdot \left( \sum\limits_{j = 1}^{d} v_{[j]}^4 \right)} = \sqrt{\sum\limits_{j = 1}^{d} u_{[j]}^2 u_{[j]}^2} \cdot \sqrt{\sum\limits_{j = 1}^{d} v_{[j]}^2 v_{[j]}^2} \\
&= \sqrt{\langle [\vb u]^2, [\vb u]^2 \rangle} \cdot \sqrt{\langle [\vb v]^2, [\vb v]^2 \rangle} = \sqrt{\| [\vb u]^2 \|^2} \cdot \sqrt{\| [\vb v]^2 \|^2},
\end{align*}
therefore
\begin{align*}
\| \vb u \odot \vb v \|^2 \leq  \| [\vb u]^2 \| \cdot \| [\vb v]^2 \|.
\end{align*}
Finally, using Proposition \eqref{Proposition_SquaredVector}, we thus obtain $\| [\vb u]^2 \| \leq \| \vb u \|^2$ and $\| [\vb v]^2 \| \leq \| \vb v \|^2$, hence the proof is finished.

\subsection{Proof of Proposition \eqref{Proposition_SwitchSummation}}

First and foremost, we have $\sum\limits_{n = 1}^{N} e_n \sum\limits_{j = 1}^{n} f_j = e_1 f_1 + e_2 (f_1 + f_2) + \ldots + e_N (f_1 + \ldots + f_N)$. At the same time, we compute $\sum\limits_{n = 1}^{N} \left( \sum\limits_{j = n}^{N} e_j \right) f_n = f_1(e_1 + \ldots + e_N) + \ldots + f_{N-1}(e_{N-1} + e_N) + f_N e_N$. Then, $\sum\limits_{n = 1}^{N} \left( \sum\limits_{j = n}^{N} e_j \right) f_n = e_1 f_1 + \ldots + e_{N-1}(f_1 + \ldots + f_{N-1}) + e_N(f_1 + \ldots + f_{N})$ and the proof is complete.

\subsection{Proof of Lemma \eqref{Lemma_AlmostSureConvergence}}

We will begin by considering, as in the previous section, the error term between the stochastic gradient and the gradient of the objective function in the shifted updates, namely $\vb* \delta_n := \vb g_n - \nabla F(\vb z_n)$. Also, let $\vb* \Delta_n := \nabla F(\vb w_n) - \nabla F(\vb z_n)$ to be the error between the iteration $\vb w_n$ and the shifted update $\vb z_n$ with respect to $F$. In order to simplify the present proof, we will break it up into several steps as follows. \\
\underline{\emph{\textbf{Step I.}}} By applying \eqref{DL_Inequality} for each $n \geq 1$, we get
\begin{align*}
F(\vb w_{n+1}) &\leq F(\vb w_n) + \langle \nabla F(\vb w_n), \vb w_{n+1} - \vb w_n \rangle + \dfrac{L}{2} \| \vb w_{n+1} - \vb w_n \|^2 \\
&= F(\vb w_n) + \langle \vb \Delta_n + \nabla F(\vb z_n), \vb w_{n+1} - \vb w_n \rangle + \dfrac{L}{2} \| \vb w_{n+1} - \vb w_n \|^2 \\
&= F(\vb w_n) + \langle \vb \Delta_n + \nabla F(\vb z_n), - \vb* \lambda_n \odot \vb g_n \rangle + \dfrac{L}{2} \| \vb* \lambda_n \odot \vb g_n \|^2 \\
&= F(\vb w_n) + \langle \vb \Delta_n + \nabla F(\vb z_n), - \vb* \lambda_n \odot \left( \vb* \delta_n + \nabla F(\vb z_n) \right) \rangle + \dfrac{L}{2} \| \vb* \lambda_n \odot \left( \vb* \delta_n + \nabla F(\vb z_n) \right) \|^2 \\
&= F(\vb w_n) - \langle \vb \Delta_n + \nabla F(\vb z_n), \vb* \lambda_n \odot \nabla F(\vb z_n) \rangle + \langle \vb \Delta_n + \nabla F(\vb z_n), - \vb* \lambda_n \odot \vb* \delta_n \rangle \\
&+ \dfrac{L}{2} \| \vb* \lambda_n \odot \left( \vb* \delta_n + \nabla F(\vb z_n) \right) \|^2, \, (w.p.1).
\end{align*}
Taking into account that $\nabla F(\vb w_n) = \vb* \Delta_n + \nabla F(\vb z_n)$, it follows that
\begin{align*}
F(\vb w_{n+1}) &\leq F(\vb w_n) - \langle \vb* \Delta_n, \vb* \lambda_n \odot \nabla F(\vb z_n) \rangle - \langle \vb* \lambda_n, [\nabla F(\vb z_n)]^2 \rangle - \langle \vb* \lambda_n \odot \vb* \delta_n, \nabla F(\vb w_n) \rangle \\
&+ \dfrac{L}{2} \| \vb* \lambda_n \odot \vb* \delta_n \|^2 + \dfrac{L}{2} \| \vb* \lambda_n \odot \nabla F(\vb z_n) \|^2 + L \langle \vb* \lambda_n \odot \vb* \delta_n, \vb* \lambda_n \odot \nabla F(\vb z_n) \rangle, \, (w.p.1).
\end{align*}
Hence
\begin{align*}
F(\vb w_{n+1}) &\leq F(\vb w_n) - \langle \vb* \Delta_n, \vb* \lambda_n \odot \nabla F(\vb z_n) \rangle - \langle \vb* \lambda_n, [\nabla F(\vb z_n)]^2 \rangle - \langle \vb* \lambda_n \odot \vb* \delta_n, \nabla F(\vb w_n) \rangle \\
&+ \dfrac{L}{2} \| \vb* \lambda_n \odot \vb* \delta_n \|^2 + \dfrac{L}{2} \| \vb* \lambda_n \odot \nabla F(\vb z_n) \|^2 + L \langle [\vb* \lambda_n]^2, \vb* \delta_n \odot \nabla F(\vb z_n) \rangle, \, (w.p.1).
\end{align*}
As a result, the computations made above are equivalent to
\begin{align*}
F(\vb w_{n+1}) &\leq F(\vb w_n) - \langle \vb* \Delta_n, \vb* \lambda_n \odot \nabla F(\vb z_n) \rangle - \langle \vb* \lambda_n, [\nabla F(\vb z_n)]^2 \rangle - \langle \vb* \lambda_n \odot \vb* \delta_n, \nabla F(\vb w_n) \rangle \\
&+ \dfrac{L}{2} \| \vb* \lambda_n \odot \vb* \delta_n \|^2 + \dfrac{L}{2} \langle [\vb* \lambda_n]^2, [\nabla F(\vb z_n)]^2 \rangle + L \langle [\vb* \lambda_n]^2, \vb* \delta_n \odot \nabla F(\vb z_n) \rangle, \, (w.p.1).
\end{align*}
Thus, we obtain
\begin{align}\label{eq:12}
F(\vb w_{n+1}) &\leq F(\vb w_n) - \Big\langle [\nabla F(\vb z_n)]^2, \vb* \lambda_n - \dfrac{L}{2} [\vb* \lambda_n]^2 \Big\rangle + \dfrac{L}{2} \| \vb* \lambda_n \odot \vb* \delta_n \|^2 - \langle \vb* \lambda_n \odot \vb* \Delta_n, \nabla F(\vb z_n) \rangle \nonumber \\
&- \langle \vb* \delta_n, \vb* \lambda_n \odot \left( \nabla F(\vb w_n) - L \vb* \lambda_n \odot \nabla F(\vb z_n) \right) \rangle, \, (w.p.1).
\end{align}
On the other hand, by using the \emph{CBS inequality} and assumption \emph{$(A_2)$}, we have the following evaluation:
\begin{align}\label{eq:13}
- \langle \vb* \lambda_n \odot \vb* \Delta_n, \nabla F(\vb z_n) \rangle &= - \langle \vb* \Delta_n, \vb* \lambda_n \odot \nabla F(\vb z_n) \rangle \leq | \langle \vb* \Delta_n, \vb* \lambda_n \odot \nabla F(\vb z_n) \rangle | \leq \| \vb* \Delta_n \| \cdot \| \vb* \lambda_n \odot \nabla F(\vb z_n) \| \nonumber \\
&= \| \nabla F(\vb w_n) - \nabla F(\vb z_n) \| \cdot \| \vb* \lambda_n \odot \nabla F(\vb z_n) \| \leq L \cdot \| \vb w_n - \vb z_n \| \cdot \| \vb* \lambda_n \odot \nabla F(\vb z_n) \|, \, (w.p.1).
\end{align}
At the same time, we have that
\begin{align}\label{eq:14}
\| \vb z_n - \vb w_n \| = \| (1 - \tilde{\vb* \gamma}_n) \odot \vb* \theta_n + \tilde{\vb* \gamma}_n \odot \vb w_n - \vb w_n \| = \| (1 - \tilde{\vb* \gamma}_n) \odot (\vb* \theta_n - \vb w_n) \|.    
\end{align}
Then, \eqref{eq:13} and \eqref{eq:14} imply
\begin{align}\label{eq:15}
- \langle \vb* \lambda_n \odot \vb* \Delta_n, \nabla F(\vb z_n) \rangle \leq L \cdot \| (1 - \tilde{\vb* \gamma}_n) \odot (\vb* \theta_n - \vb w_n) \| \cdot \| \vb* \lambda_n \odot \nabla F(\vb z_n) \|, \, (w.p.1). 
\end{align}
Furthermore, by using the basic inequality $ab \leq \dfrac{a^2}{2} + \dfrac{b^2}{2}$ and Proposition \eqref{Proposition_HadamardProduct} when $\tilde{\vb* \gamma}_n \in \mathbb{R}^d$ (when $\tilde{\vb* \gamma}_n \in \mathbb{R}$, then we can use $|1 - \tilde{\vb* \gamma}_n|$ instead of the Euclidean norm, as specified in Remark \eqref{Remark_NormBound}), we find
\begin{align}\label{eq:16}
\| (1 - \tilde{\vb* \gamma}_n) \odot (\vb* \theta_n - \vb w_n) \| \cdot \| \vb* \lambda_n \odot \nabla F(\vb z_n) \| &\leq \dfrac{1}{2} \| (1 - \tilde{\vb* \gamma}_n) \odot (\vb* \theta_n - \vb w_n) \|^2 + \dfrac{1}{2} \langle [\vb* \lambda_n]^2, [\nabla F(\vb z_n)]^2 \rangle \nonumber \\
&\leq \dfrac{\| 1 - \tilde{\vb* \gamma}_n \|^2}{2} \|\vb* \theta_n - \vb w_n \|^2 + \dfrac{1}{2} \langle [\vb* \lambda_n]^2, [\nabla F(\vb z_n)]^2 \rangle
, \, (w.p.1).
\end{align}
By combining \eqref{eq:15} with \eqref{eq:16}, one has that
\begin{align}\label{eq:17}
- \langle \vb* \lambda_n \odot \vb* \Delta_n, \nabla F(\vb z_n) \rangle \leq  \dfrac{L \| 1 - \tilde{\vb* \gamma}_n \|^2}{2} \| \vb* \theta_n - \vb w_n \|^2 + \dfrac{L}{2} \langle [\vb* \lambda_n]^2, [\nabla F(\vb z_n)]^2 \rangle, \, (w.p.1).
\end{align}
Putting together \eqref{eq:12} and \eqref{eq:17}, for each $n \geq 1$ we infer that
\begin{align}\label{eq:18}
F(\vb w_{n+1}) &\leq F(\vb w_n) - \langle [\nabla F(\vb z_n)]^2, \vb* \lambda_n - L [\vb* \lambda_n]^2 \rangle + \dfrac{L}{2} \| \vb* \lambda_n \odot \vb* \delta_n \|^2 + \dfrac{L \| 1 - \tilde{\vb* \gamma}_n\|^2}{2} \| \vb* \theta_n - \vb w_n \|^2 \nonumber \\
&- \langle \vb* \delta_n, \vb* \lambda_n \odot \left( \nabla F(\vb w_n) - L \vb* \lambda_n \odot \nabla F(\vb z_n) \right) \rangle, \, (w.p.1).
\end{align}
\underline{\emph{\textbf{Step II.}}} In the following, we will consider the evaluation of the term $\vb* \theta_n - \vb w_n$. Taking $n \geq 1$, we proceed below in the following way.
\begin{align*}
\vb* \theta_{n+1} - \vb w_{n+1} &= (\vb y_n - \vb* \alpha_n \odot \vb g_n) - (\vb w_n - \vb* \lambda_n \odot \vb g_n) \\
&= (\vb y_n - \vb w_n) + (\vb* \lambda_n - \vb* \alpha_n) \odot \vb g_n \\
&= \left( (1 - \mu_n) \vb* \theta_n + \mu_n \vb w_n - \vb w_n \right) + (\vb* \lambda_n - \vb* \alpha_n) \odot \vb g_n \\
&= (1 - \mu_n) (\vb* \theta_n - \vb w_n) + (\vb* \lambda_n - \vb* \alpha_n) \odot \vb g_n.
\end{align*}
Similarly to \cite{GhadimiLanAccelerated}, we consider
\begin{align}\label{eq:Gamma}
\Gamma_n :=
\begin{cases}
1; & n = 1 \\
(1 - \mu_n) \Gamma_{n-1}; & n \geq 2.
\end{cases}
\end{align}
Using the fact that $\mu_1 = 1$ and that $\mu_n \in (0,1)$ for each $n \geq 2$, we find that $\Gamma_n > 0$ for every $n \geq 1$. Furthermore, for the first index, i.e. $j = 1$, we get
\begin{align*}
\dfrac{\vb* \theta_2 - \vb w_2}{\Gamma_1} &= \dfrac{(1 - \mu_1)(\vb* \theta_1 - \vb w_1)}{\Gamma_1} + \dfrac{\vb* \lambda_1 - \vb* \alpha_1}{\Gamma_1} \odot \vb g_1 = \dfrac{\vb* \lambda_1 - \vb* \alpha_1}{\Gamma_1} \odot \vb g_1,
\end{align*}
hence
\begin{align}\label{eq:j_1}
\dfrac{\vb* \theta_2 - \vb w_2}{\Gamma_1} = \dfrac{\vb* \lambda_1 - \vb* \alpha_1}{\Gamma_1} \odot \vb g_1.    
\end{align}
For each $j \geq 2$ we have that $\dfrac{1 - \mu_{j}}{\Gamma_j} = \dfrac{1}{\Gamma_{j-1}}$, then taking $j \geq 2$ we get
\begin{align*}
\dfrac{\vb* \theta_{j+1} - \vb w_{j+1}}{\Gamma_j} = \dfrac{(1 - \mu_{j}) (\vb* \theta_{j} - \vb w_{j})}{\Gamma_j} + \dfrac{(\vb* \lambda_{j} - \vb* \alpha_{j})}{\Gamma_{j}} \odot \vb g_{j} = \dfrac{\vb* \theta_{j} - \vb w_{j}}{\Gamma_{j-1}} + \dfrac{(\vb* \lambda_{j} - \vb* \alpha_{j})}{\Gamma_{j}} \odot \vb g_{j}.
\end{align*}
Making the summation from $2$ to $n$ (for $n \geq 2$), it follows that
\begin{align*}
\sum\limits_{j = 2}^{n} \dfrac{\vb* \theta_{j+1} - \vb w_{j+1}}{\Gamma_j} - \sum\limits_{j = 2}^{n} \dfrac{\vb* \theta_{j} - \vb w_{j}}{\Gamma_{j-1}} = \sum\limits_{j = 2}^{n} \dfrac{(\vb* \lambda_{j} - \vb* \alpha_{j})}{\Gamma_{j}} \odot \vb g_{j},
\end{align*}
therefore
\begin{align*}
\dfrac{\vb* \theta_{n+1} - \vb w_{n+1}}{\Gamma_n} = \dfrac{\vb* \theta_2 - \vb w_2}{\Gamma_1} + \sum\limits_{j = 2}^{n} \dfrac{(\vb* \lambda_{j} - \vb* \alpha_{j})}{\Gamma_{j}} \odot \vb g_{j},
\end{align*}
which, by using \eqref{eq:j_1}, it means that
\begin{align*}
\dfrac{\vb* \theta_{n+1} - \vb w_{n+1}}{\Gamma_n} = \dfrac{\vb* \lambda_1 - \vb* \alpha_1}{\Gamma_1} \odot \vb g_1 + \sum\limits_{j = 2}^{n} \dfrac{(\vb* \lambda_{j} - \vb* \alpha_{j})}{\Gamma_{j}} \odot \vb g_{j},
\end{align*}
concluding for each $n \geq 2$ that
\begin{align}\label{eq:19}
\vb* \theta_{n+1} - \vb w_{n+1} = \Gamma_n \cdot \sum\limits_{j = 1}^{n} \dfrac{(\vb* \lambda_{j} - \vb* \alpha_{j})}{\Gamma_{j}} \odot \vb g_{j}.
\end{align}
From \eqref{eq:j_1}, we obtain that \eqref{eq:19} holds also for every $n \geq 1$. Now, for each $n \geq 2$, using the definition of $\Gamma_n$, we have that $\Gamma_n = (1 - \mu_n) \Gamma_{n-1}$, hence $\mu_n = 1 - \dfrac{\Gamma_n}{\Gamma_{n-1}}$. Since $\mu_1 = 1$, we obtain
\begin{align}\label{eq:20}
\sum\limits_{j = 1}^{n} \dfrac{\mu_{j}}{\Gamma_{j}} = \dfrac{\mu_1}{\Gamma_1} + \sum\limits_{j = 2}^{n} \dfrac{\mu_{j}}{\Gamma_{j}} = \dfrac{\mu_1}{\Gamma_1} + \sum\limits_{j = 2}^{n} \dfrac{1}{\Gamma_{j}} \left( 1 - \dfrac{\Gamma_{j}}{\Gamma_{j-1}} \right) = \dfrac{1}{\Gamma_1} + \sum\limits_{j = 2}^{n} \left( \dfrac{1}{\Gamma_{j}} - \dfrac{1}{\Gamma_{j-1}} \right) = \dfrac{1}{\Gamma_{n}},  
\end{align}
thus, for every $n \geq 2$, one has that
\begin{align}\label{sum_Gamma_n}
\Gamma_n \sum\limits_{j = 1}^{n} \dfrac{\mu_j}{\Gamma_j} = 1.    
\end{align}
At the same time, it is trivial to show that \eqref{sum_Gamma_n} holds also for $n = 1$ since $\mu_1 = 1$, hence \eqref{sum_Gamma_n} is valid for every $n \geq 1$. Using the computations and the above explanations, and utilizing \eqref{eq:19} and \eqref{eq:20} and proceeding as in \cite{GhadimiLanAccelerated} by applying Jensen's inequality for the squared Euclidean norm $\| \cdot \|_2$, for every $n \geq 1$ and by taking into account that $\mu_j \neq 0$ for $j \in \lbrace 1, \ldots, n \rbrace$, we get that
\begin{align*}
\| \vb* \theta_{n+1} - \vb w_{n+1} \|^2 &= \left\| \Gamma_n \cdot \sum\limits_{j = 1}^{n} \dfrac{\vb* \lambda_j - \vb* \alpha_j}{\Gamma_j} \odot \vb g_j \right\|^2 = \left\| \Gamma_n \cdot \sum\limits_{j = 1}^{n} \dfrac{\mu_j}{\Gamma_j} \cdot \left( \dfrac{\vb* \lambda_j - \vb* \alpha_j}{\mu_j} \odot \vb g_j \right) \right\|^2 \\
&\leq \Gamma_n \cdot \sum\limits_{j = 1}^{n} \dfrac{\mu_j}{\Gamma_j} \cdot \left\| \dfrac{\vb* \lambda_j - \vb* \alpha_j}{\mu_j} \odot \vb g_j \right\|^2 = \Gamma_n \cdot \sum\limits_{j = 1}^{n} \dfrac{1}{\mu_j \Gamma_j} \cdot \left\| (\vb* \lambda_j - \vb* \alpha_j) \odot \vb g_j \right\|^2, \, (w.p.1).
\end{align*}
Therefore, for each $n \geq 2$
\begin{align}\label{ineq_norm_theta_w}
\| \vb* \theta_n - \vb w_n \|^2 \leq \Gamma_{n-1} \cdot \sum\limits_{j = 1}^{n-1} \dfrac{1}{\mu_j \Gamma_j} \cdot \left\| (\vb* \lambda_j - \vb* \alpha_j) \odot \vb g_j \right\|^2, \, (w.p.1).
\end{align}
Since for every $n \geq 2$ one has $\Gamma_{n-1} = \dfrac{\Gamma_n}{1 - \mu_n}$ and $\mu_n \in (0, 1)$ for $n \geq 2$, and by the fact that all the terms are non-negative, we infer that
\begin{align}\label{eq:21}
\| \vb* \theta_n - \vb w_n \|^2 \leq \dfrac{\Gamma_n}{1 - \mu_n} \cdot \sum\limits_{j = 1}^{n-1} \dfrac{1}{\mu_j \Gamma_j} \cdot \left\| (\vb* \lambda_j - \vb* \alpha_j) \odot \vb g_j \right\|^2 \leq \dfrac{\Gamma_n}{1 - \mu_n} \cdot \sum\limits_{j = 1}^{n} \dfrac{1}{\mu_j \Gamma_j} \cdot \left\| (\vb* \lambda_j - \vb* \alpha_j) \odot \vb g_j \right\|^2, \, (w.p.1).
\end{align}
\underline{\emph{\textbf{Step III.}}} Combining \eqref{eq:18} and \eqref{eq:21} for every $n \geq 2$ and taking into consideration that $\mu_n \neq 1$ for every $n \geq 2$, we get that
\begin{align*}
F(\vb w_{n+1}) &\leq F(\vb w_n) - \langle [\nabla F(\vb z_n)]^2, \vb* \lambda_n - L [\vb* \lambda_n]^2 \rangle + \dfrac{L}{2} \| \vb* \lambda_n \odot \vb* \delta_n \|^2 + \dfrac{L \Gamma_n \|1 - \tilde{\vb* \gamma}_n\|^2}{2(1-\mu_n)} \sum\limits_{j = 1}^{n} \dfrac{1}{\mu_j \Gamma_j} \cdot \left\| (\vb* \lambda_j - \vb* \alpha_j) \odot \vb g_j \right\|^2 \\
&- \langle \vb* \delta_n, \vb* \lambda_n \odot \left( \nabla F(\vb w_n) - L \vb* \lambda_n \odot \nabla F(\vb z_n) \right) \rangle, \, (w.p.1).
\end{align*}
From the theorem's hypotheses we find that $ \dfrac{\|1-\tilde{\vb* \gamma}_n\|^2}{1 - \mu_n} \leq B$ for every $n \geq 2$. Then, for $n \geq 2$, the above almost sure inequality is equivalent to
\begin{align}\label{eq:DescentPropertyObjectiveFunction}
F(\vb w_{n+1}) &\leq F(\vb w_n) - \langle [\nabla F(\vb z_n)]^2, \vb* \lambda_n - L [\vb* \lambda_n]^2 \rangle + \dfrac{L}{2} \| \vb* \lambda_n \odot \vb* \delta_n \|^2 + \dfrac{L B \Gamma_n}{2} \sum\limits_{j = 1}^{n} \dfrac{1}{\mu_j \Gamma_j} \cdot \langle [\vb* \lambda_j - \vb* \alpha_j]^2, [\vb g_j]^2 \rangle \nonumber \\
&- \langle \vb* \delta_n, \vb* \lambda_n \odot \left( \nabla F(\vb w_n) - L \vb* \lambda_n \odot \nabla F(\vb z_n) \right) \rangle, \, (w.p.1),
\end{align}
Taking into account that \eqref{eq:18} holds also for $n = 1$, and since $\tilde{\vb* \gamma}_1 = 1$ then $\dfrac{L \|1-\tilde{\vb* \gamma}_1\|^2}{2} \| \vb* \theta_1 - \vb w_1 \|^2 = 0$, which implies that \eqref{eq:DescentPropertyObjectiveFunction} is valid for $n \geq 1$. By considering $N \geq n$, where $n \geq 1$, and taking the sum from $1$ to $N$, we obtain
\begin{align*}
\sum\limits_{n = 1}^{N} \left( F(\vb w_{n+1}) - F(\vb w_n) \right) &\leq - \sum\limits_{n = 1}^{N} \langle [\nabla F(\vb z_n)]^2, \vb* \lambda_n - L [\vb* \lambda_n]^2 \rangle + \dfrac{L}{2} \sum\limits_{n = 1}^{N} \langle [\vb* \delta_n]^2, [\vb* \lambda_n]^2 \rangle \\
&- \sum\limits_{n = 1}^{N} \langle \vb* \delta_n, \vb* \lambda_n \odot \left( \nabla F(\vb w_n) - L \vb* \lambda_n \odot \nabla F(\vb z_n) \right) \rangle \\
&+ \dfrac{LB}{2} \sum\limits_{n = 1}^{N} \Gamma_n \sum\limits_{j = 1}^{n} \dfrac{1}{\mu_j \Gamma_j} \cdot \langle [\vb* \lambda_j - \vb* \alpha_j]^2, [\vb g_j]^2 \rangle, \, (w.p.1),
\end{align*}
Using Proposition \eqref{Proposition_SwitchSummation} for the above almost sure inequality, it implies that
\begin{align*}
\sum\limits_{n = 1}^{N} \left( F(\vb w_{n+1}) - F(\vb w_n) \right) &\leq - \sum\limits_{n = 1}^{N} \langle [\nabla F(\vb z_n)]^2, \vb* \lambda_n - L [\vb* \lambda_n]^2 \rangle + \dfrac{L}{2} \sum\limits_{n = 1}^{N} \langle [\vb* \delta_n]^2, [\vb* \lambda_n]^2 \rangle \\
&- \sum\limits_{n = 1}^{N} \langle \vb* \delta_n, \vb* \lambda_n \odot \left( \nabla F(\vb w_n) - L \vb* \lambda_n \odot \nabla F(\vb z_n) \right) \rangle \\
&+ \dfrac{LB}{2} \sum\limits_{n = 1}^{N} \left( \sum\limits_{j = n}^{N} \Gamma_j \right)  \dfrac{1}{\mu_n \Gamma_n} \cdot \langle [\vb* \lambda_n - \vb* \alpha_n]^2, [\vb g_n]^2 \rangle, \, (w.p.1),
\end{align*}
For the forthcoming computations, we will use that
\begin{align*}
\langle [\vb* \lambda_n - \vb* \alpha_n]^2, [\vb g_n]^2 \rangle &= \langle [\vb* \lambda_n - \vb* \alpha_n]^2, [\vb* \delta_n + \nabla F(\vb z_n)]^2 \rangle = \langle [\vb* \lambda_n - \vb* \alpha_n]^2, [\vb* \delta_n]^2 + [\nabla F(\vb z_n)]^2 + 2 \vb* \delta_n \odot \nabla F(\vb z_n) \rangle \\
&= \langle [\vb* \lambda_n - \vb* \alpha_n]^2, [\vb* \delta_n]^2 \rangle + \langle [\vb* \lambda_n - \vb* \alpha_n]^2, [\nabla F(\vb z_n)]^2 \rangle + 2 \langle [\vb* \lambda_n - \vb* \alpha_n]^2, \vb* \delta_n \odot \nabla F(\vb z_n) \rangle.
\end{align*}
Furthermore, by employing for every $n \geq 1$ the notation
\begin{align}\label{C_nN}
C_{n, N} := \sum\limits_{j = n}^{N} \Gamma_j,    
\end{align}
and utilizing Proposition \eqref{Proposition_HadamardProduct} almost surely, such that $\langle [\vb* \delta_n]^2, [\vb* \lambda_n]^2 \rangle = \| \vb* \lambda_n \odot \vb* \delta_n \|^2 \leq \left( \| \vb* \lambda_n \| \cdot \| \vb* \delta_n \| \right)^2 = \| \vb* \lambda_n \|^2 \cdot \| \vb* \delta_n \|^2$, we get
\begin{align*}
F(\vb w_{N+1}) - F(\vb w_1) &\leq - \sum\limits_{n = 1}^{N} \langle [\nabla F(\vb z_n)]^2, \vb* \lambda_n - L [\vb* \lambda_n]^2 \rangle + \dfrac{L}{2} \sum\limits_{n = 1}^{N} \| \vb* \delta_n\|^2 \cdot \| \vb* \lambda_n \|^2 \\
&- \sum\limits_{n = 1}^{N} \langle \vb* \delta_n, \vb* \lambda_n \odot \left( \nabla F(\vb w_n) - L \vb* \lambda_n \odot \nabla F(\vb z_n) \right) \rangle +  LB \sum\limits_{n = 1}^{N}  \dfrac{C_{n,N}}{\mu_n \Gamma_n} \cdot \langle [\vb* \lambda_n - \vb* \alpha_n]^2, \vb* \delta_n \odot \nabla F(\vb z_n) \rangle \\
&+ \dfrac{LB}{2} \sum\limits_{n = 1}^{N}  \dfrac{C_{n,N}}{\mu_n \Gamma_n} \cdot \langle [\vb* \lambda_n - \vb* \alpha_n]^2, [\vb* \delta_n]^2 \rangle + \dfrac{LB}{2} \sum\limits_{n = 1}^{N}  \dfrac{C_{n,N}}{\mu_n \Gamma_n} \cdot \langle [\vb* \lambda_n - \vb* \alpha_n]^2, [\nabla F(\vb z_n)]^2 \rangle, \, (w.p.1),
\end{align*}
For simplicity, we denote
\begin{align*}
\vb Q_n := \langle \vb* \delta_n, \vb* \lambda_n \odot \left( \nabla F(\vb w_n) - L \vb* \lambda_n \odot \nabla F(\vb z_n) \right) - LB \dfrac{C_{n,N}}{\mu_n \Gamma_n} [\vb* \lambda_n - \vb* \alpha_n]^2 \odot \nabla F(\vb z_n) \rangle.    
\end{align*}
Utilizing Proposition \eqref{Proposition_SquaredVector} along with the CBS inequality, i.e. $\langle [\vb* \lambda_n - \vb* \alpha_n]^2, [\vb* \delta_n]^2 \rangle \leq \| [\vb* \lambda_n - \vb* \alpha_n]^2 \| \cdot \| [\vb* \delta_n]^2 \| \leq \| \vb* \lambda_n - \vb* \alpha_n \|^2 \cdot \| \vb* \delta_n \|^2$ it follows
\begin{align*}
F(\vb w_{N+1}) - F(\vb w_1) &\leq - \sum\limits_{n = 1}^{N} \langle [\nabla F(\vb z_n)]^2, \vb* \lambda_n - L [\vb* \lambda_n]^2 - \dfrac{LB}{2} \dfrac{C_{n,N}}{\mu_n \Gamma_n} [\vb* \lambda_n - \vb* \alpha_n]^2 \rangle - \sum\limits_{n = 1}^{N} \vb Q_n \\
&+ \dfrac{L}{2} \sum\limits_{n = 1}^{N} \| \vb* \delta_n \|^2 \cdot \left( \| \vb* \lambda_n \|^2 + B \dfrac{C_{n,N}}{\mu_n \Gamma_n} \| \vb* \lambda_n - \vb* \alpha_n \|^2 \right) , \, (w.p.1),
\end{align*}
Denoting
\begin{align*}
\vb D_n := \vb* \lambda_n - L [\vb* \lambda_n]^2 - \dfrac{LB}{2} \dfrac{C_{n,N}}{\mu_n \Gamma_n} [\vb* \lambda_n - \vb* \alpha_n]^2,    
\end{align*}
and using the theorem's assumptions from which we know that for each $n \geq 1$, there exists $m_n > 0$, such that $\vb D_n \geq m_n$ almost surely, along with the fact that $[\nabla F(\vb z_n)]^2 \geq 0$ with respect to Hadamard notations, we thus have
\begin{align*}
F(\vb w_{N+1}) - F(\vb w_1) + \sum\limits_{n = 1}^{N} m_n \cdot \| \nabla F(\vb z_n) \|^2 &\leq - \sum\limits_{n = 1}^{N} \vb Q_n + \dfrac{L}{2} \sum\limits_{n = 1}^{N} \| \vb* \delta_n \|^2 \cdot \left( \| \vb* \lambda_n \|^2 + B \dfrac{C_{n,N}}{\mu_n \Gamma_n} \| \vb* \lambda_n - \vb* \alpha_n \|^2 \right) , \, (w.p.1),
\end{align*}
therefore we have finished the proof.

\subsection{Proof of Theorem \eqref{Theorem_AsymptoticBehavior}}

By taking into account the hypothesis that $\| \vb* \lambda_n \|^2 + B \dfrac{C_{n,N}}{\mu_n \Gamma_n} \| \vb* \lambda_n - \vb* \alpha_n \|^2 \leq R_n$ almost surely and in the Hadamard sense for $n \geq 1$, along with \eqref{eq:Results_AlmostSureInequality} we have that
\begin{align*}
F(\vb w_{N+1}) - F(\vb w_1) + \sum\limits_{n = 1}^{N} m_n \cdot \| \nabla F(\vb z_n) \|^2 &\leq - \sum\limits_{n = 1}^{N} \vb Q_n + \dfrac{L}{2} \sum\limits_{n = 1}^{N} R_n \| \vb* \delta_n \|^2, \, (w.p.1),
\end{align*}
where $m_n$, $C_{n,N}$ and $\vb Q_n$ were properly defined in Lemma \eqref{Lemma_AlmostSureConvergence}. Applying the conditional expectation (along with its linearity property) with respect to $\mathcal{F}_N$ (which was defined in section \eqref{Section_AssumptionSetting}), it follows that
\begin{align*}
\mathbb{E} \left[ F(\vb w_{N+1}) \, | \, \mathcal{F}_N \right] + \sum\limits_{n = 1}^{N} m_n \cdot \mathbb{E} \left[ \| \nabla F(\vb z_n) \|^2 \, | \, \mathcal{F}_N \right] &\leq F(\vb w_1) - \sum\limits_{n = 1}^{N} \mathbb{E} \left[ \vb Q_n \, | \, \mathcal{F}_N \right] + \dfrac{L}{2} \sum\limits_{n = 1}^{N} R_n \cdot \mathbb{E} \left[ \| \vb* \delta_n \|^2 \, | \, \mathcal{F}_N \right], \, (w.p.1).
\end{align*}
Taking the total expectation with respect to the underlying probability space (along with its linearity property), leads to
\begin{align}\label{eq:mean_totalExp}
\mathbb{E} \left[ \mathbb{E} \left[ F(\vb w_{N+1}) \, | \, \mathcal{F}_N \right] \right] + \sum\limits_{n = 1}^{N} m_n \cdot \mathbb{E} \left[ \mathbb{E} \left[ \| \nabla F(\vb z_n) \|^2 \, | \, \mathcal{F}_N \right] \right] &\leq \mathbb{E} \left[ F(\vb w_1) \right] - \sum\limits_{n = 1}^{N} \mathbb{E} \left[ \mathbb{E} \left[ \vb Q_n \, | \, \mathcal{F}_N \right] \right] \\
&+ \dfrac{L}{2} \sum\limits_{n = 1}^{N} R_n \cdot \mathbb{E} \left[ \mathbb{E} \left[ \| \vb* \delta_n \|^2 \, | \, \mathcal{F}_N \right] \right].
\end{align}
For the sake of clarity, we consider 
$\vb Y_{n, N} := \mathbb{E}[ \| \vb* \delta_n \|^2 \, | \, \mathcal{F}_N].$  
Since $(\mathcal{F}_n)_{n \geq 1}$ is an increasing sequence of $\sigma$ fields, so $\mathcal{F}_n \subset \mathcal{F}_N$, by employing the tower property of the conditional expectation we get that
\begin{align}\label{eq:23}
\mathbb{E}[\vb Y_{n,N} \, | \, \mathcal{F}_n] = \mathbb{E} [\| \vb* \delta_n \|^2 \, | \, \mathcal{F}_n] \leq \sigma^2, \, (w.p.1).
\end{align}
On the other hand, according to the law of total expectation,
\begin{align}\label{eq:24}
\mathbb{E} [\mathbb{E}[\vb Y_{n,N} \, | \, \mathcal{F}_n]] = \mathbb{E} [\vb Y_{n,N}].    
\end{align}
From \eqref{eq:23} and \eqref{eq:24}, we obtain
\begin{align}\label{eq:25}
\mathbb{E}[\vb Y_{n,N}] \leq \sigma^2.
\end{align}
In a similar manner as in \eqref{eq:25}, we denote $\vb W_{n,N} := \mathbb{E} \left[ \| \nabla F(\vb z_n) \|^2 \, | \, \mathcal{F}_N \right]$. Hence, using a similar analysis as before, we get that $\mathbb{E} \left[ \vb W_{n,N} \right] = \mathbb{E} \left[ \mathbb{E} \left[ \vb W_{n,N} \, | \, \mathcal{F}_n \right] \right] = \mathbb{E} \left[ \| \nabla F(\vb z_n) \|^2 \right]$, where we have used that $\mathbb{E} \left[ \vb W_{n,N} \, | \, \mathcal{F}_n \right] = \mathbb{E} \left[ \| \nabla F(\vb z_n) \|^2 \, | \, \mathcal{F}_n \right]$, therefore $\mathbb{E} \left[ \mathbb{E} \left[ \vb W_{n,N} \, | \, \mathcal{F}_n \right] \right] = \mathbb{E} \left[ \mathbb{E} \left[ \| \nabla F(\vb z_n) \|^2 \, | \, \mathcal{F}_n \right] \right]$. This reasoning is true from the measurability assumption \emph{$(A_5)$}. On the other hand, denote 
$$\vb B_n := \vb* \lambda_n \odot \left( \nabla F(\vb w_n) - L \vb* \lambda_n \odot \nabla F(\vb z_n) \right) - LB \dfrac{C_{n,N}}{\mu_n \Gamma_n} [\vb* \lambda_n - \vb* \alpha_n]^2 \odot \nabla F(\vb z_n).$$
Regardless of the fact that the above term $\vb B_n$ depends also on $N$ through $C_{n,N}$ we have used the above notation in order to simplify the technicalities. Then $\vb Q_n = \langle \vb* \delta_n, \vb B_n \rangle$. Denoting $\vb Z_{n,N} := \mathbb{E} \left[ \vb Q_n \, | \, \mathcal{F}_N \right]$, using a similar technique as before, we find that
$$ \mathbb{E} \left[ \vb Z_{n,N} \right] = \mathbb{E} \left[ \mathbb{E} \left[ \langle \vb* \delta_n, \vb B_n \rangle \, | \, \mathcal{F}_n \right] \right].$$
We will make use of the measurability of the random vectors with respect to the assumption \emph{$(A_5)$}. Therefore, we know that $\vb* \lambda_n$, $\vb* \alpha_n$, $\vb w_n$ and $\vb z_n$ are all of them $\mathcal{F}_n$-measurable, hence $\vb B_n$ is also $\mathcal{F}_n$-measurable. This means that each of its components is also a measurable random variable with respect to $\mathcal{F}_n$. Due to this measurability aspect of $\vb B_n$ with respect to $\mathcal{F}_n$, and employing the notations from subsection \eqref{Subsection_Notations}, we find that
\begin{align*}
\mathbb{E} \left[ \vb Q_n \, | \, \mathcal{F}_n \right] = \mathbb{E} \left[ \langle \vb* \delta_n, \vb B_n \rangle \, | \, \mathcal{F}_n \right] &= \mathbb{E} \left[ \sum\limits_{k=1}^{d} [\vb* \delta_n]_{[k]} \cdot [\vb B_n]_{[k]}  \, \Big| \, \mathcal{F}_n \right] = \sum\limits_{k=1}^{d} \mathbb{E} \left[ \vb* \delta_{n, [k]} \cdot B_{n, [k]}  \, \Big| \, \mathcal{F}_n \right] \\
&= \sum\limits_{k=1}^{d}  B_{n, [k]} \cdot \mathbb{E}  \left[ \vb* \delta_{n, [k]}  \, \Big| \, \mathcal{F}_n \right] = \langle \vb B_n, \mathbb{E} [\vb* \delta_n \, | \, \mathcal{F}_n] \rangle = 0,
\end{align*}
which is well defined since the random vectors $\vb B_n$ and $\vb* \delta_n$ take values in $\mathbb{R}^d$, hence the inner products with respect to these random vectors take value in $\mathbb{R}$. Then, we get that $\mathbb{E} \left[ \vb Z_{n, N} \right] = 0$. This, along with the above remarks and the property of total expectation imply that \eqref{eq:mean_totalExp} becomes
\begin{align*}
\mathbb{E} \left[ F(\vb w_{N+1}) \right] + \sum\limits_{n = 1}^{N} m_n \cdot \mathbb{E} \left[ \| \nabla F(\vb z_n) \|^2 \right] &\leq \mathbb{E} \left[ F(\vb w_1) \right] + \dfrac{L}{2} \sigma^2 \sum\limits_{n = 1}^{N} R_n.
\end{align*}
From assumption \emph{$(A_1)$}, there exists $\underline{M} > 0$, such that $F(\vb w_{N+1}) \geq \underline{M}$ almost surely, hence $\mathbb{E}[F(\vb w_{N+1})] \geq \underline{M}$, so denoting $\overline{M} := \mathbb{E}[F(\vb w_1)] - \underline{M}$ and employing the notations from Lemma \eqref{Lemma_AlmostSureConvergence}, we finally find that
\begin{align}\label{eq:26}
\sum\limits_{n = 1}^{N} m_n \cdot \mathbb{E} [\| \nabla F(\vb z_n) \|^2] &\leq \overline{M} + \dfrac{L}{2} \sigma^2 \sum\limits_{n = 1}^{N} R_n,
\end{align}
This means that \eqref{eq:26} leads to
\begin{align*}
\mathbb{E} \left[ \min_{n = 1, \ldots, N} \| \nabla F(\vb z_n) \|^2 \right] &\leq \dfrac{\overline{M} + \dfrac{L}{2} \sigma^2 \sum\limits_{n = 1}^{N} R_n}{\sum\limits_{n = 1}^{N} m_n},
\end{align*}
Taking the theorem's assumptions into consideration, and making $N \to +\infty$ in \eqref{eq:26}, we finally arrive at
\begin{align*}
\sum\limits_{n = 1}^{+\infty} m_n \cdot \mathbb{E} [\| \nabla F(\vb z_n) \|^2] &\leq \overline{M} + \dfrac{L}{2} \sigma^2 \sum\limits_{n = 1}^{+\infty} R_n < +\infty,
\end{align*}
and the proof is complete.

\subsection{Proof of Proposition \eqref{ChebyshevProbabilityGradient}}

Denote 
\begin{align*}
A_{n, N} := \dfrac{\overline{M} + \dfrac{L}{2} \sigma^2 \sum\limits_{n = 1}^{N} R_n}{\sum\limits_{n = 1}^{N} m_n}.
\end{align*}
Let $\vb X_n$ stand for $\min\limits_{n = 1, \ldots, N} \| \nabla F(\vb z_n) \|^2$, hence $\mathbb{E}[\vb X_n] \leq A_{n, N}$ due to \eqref{MinimumGradientRate}. 
Using $\mathbb{P}$ to represent the probability distribution of the underlying probability space, by \emph{Markov inequality}, it follows that
\begin{align*}
\mathbb{P} \left( |\vb X_n| > \dfrac{A_{n,N}}{\delta^\prime} \right) \leq \dfrac{\mathbb{E}[\vb X_n]}{\dfrac{A_{n,N}}{\delta^\prime}} \leq \delta^\prime.  
\end{align*}
Hence
\begin{align*}
P \left( |\vb X_n| \leq \dfrac{A_{n,N}}{\delta^\prime} \right) \geq 1 - \delta^\prime,    
\end{align*}
so the proof is finished.

\subsection{Proof of Proposition \eqref{Proposition_ComplexityResult}}

From the hypotheses, we have that $\vb* \lambda_n - \vb* \alpha_n = (M-1) \vb* \alpha_n$ (not necessarily greater than $0$), in addition to the fact that $\vb* \lambda_n \leq \dfrac{M C \nu}{\varepsilon \sqrt{N}}$, because $\vb* \alpha_n \leq \dfrac{\nu}{\varepsilon} \eta_n$. We denote $Q:= M^2 + B \hat{C} (M-1)^2$ where $B$ is the term that appears in Lemma \eqref{Lemma_AlmostSureConvergence}. Applying Remark \eqref{Remark_PositivityNormInequality}, along with $\| \vb* \lambda_n - \vb* \alpha_n \|^2 = (M - 1)^2 \| \vb* \alpha_n \|^2 \leq \dfrac{d(M - 1)^2 C^2 \nu^2}{\varepsilon^2 N}$ then, from simple computations, we infer that
\begin{align*}
\| \vb* \lambda_n \|^2 + B \dfrac{C_{n,N}}{\mu_n \Gamma_n} \| \vb* \lambda_n - \vb* \alpha_n \|^2 \leq \dfrac{dM^2C^2 \nu^2}{\varepsilon^2 N} + \dfrac{dB \hat{C}(M-1)^2 C^2 \nu^2}{\varepsilon^2 N} = \dfrac{d Q C^2 \nu^2}{\varepsilon^2 N},    
\end{align*}
hence we define $R_n := \dfrac{d Q C^2 \nu^2}{\varepsilon^2 N}$. Taking into account \eqref{eq:26}, it implies that 
\begin{align*}
\dfrac{1}{N} \sum\limits_{n = 1}^{N} m_n \cdot \mathbb{E} [\| \nabla F(\vb z_n) \|^2] &\leq \dfrac{ \overline{M} + \dfrac{L}{2} \sigma^2 \sum\limits_{n = 1}^{N} R_n }{N},
\end{align*}
and since $m_n = \dfrac{M \nu \eta_n}{2 (\mathcal{K} + \varepsilon)} = \dfrac{M C \nu}{2(\mathcal{K}+\varepsilon) \sqrt{N}}$ and using the definition of $R_n$, we get that
\begin{align*}
\dfrac{1}{N} \sum\limits_{n = 1}^{N} \mathbb{E} [\| \nabla F(\vb z_n) \|^2] &\leq 2 (\mathcal{K}+\varepsilon) \dfrac{ \overline{M} + \dfrac{L \sigma^2}{2} \dfrac{d Q C^2 \nu^2}{\varepsilon^2} }{M C \nu \sqrt{N}} \leq \delta,
\end{align*}
therefore $N \geq 4 \left( (\mathcal{K}+\varepsilon) \dfrac{ \overline{M} + \dfrac{L \sigma^2}{2} \dfrac{d Q C^2 \nu^2}{\varepsilon^2} }{M C \nu \delta} \right)^2 = \mathcal{O} \left( \dfrac{1}{\delta^2} \right)$ and the proof is complete.

\clearpage
\subsection{Additional numerical results}\label{subsection_additional_numerical_results}


\setlength{\tabcolsep}{10pt}
\renewcommand{\arraystretch}{0.75}
\begin{table}[h!]
\begin{tabular}{>{\small}l || >{\scriptsize}c | >{\scriptsize}c | >{\scriptsize}c || >{\scriptsize}c | >{\scriptsize}c | >{\scriptsize}c}
\toprule
 & \multicolumn{3}{c}{\texttt{LR-MNIST}} & \multicolumn{3}{c}{\texttt{CNN-CIFAR10}}\\
\cmidrule(lr){2-2} \cmidrule(lr){3-3} \cmidrule(lr){4-4} \cmidrule(lr){5-5} \cmidrule(lr){6-6} \cmidrule(lr){7-7}
\textbf{Metrics} &AAMMSU & AMSGrad & \texttt{($\eta$, bs)} & AAMMSU & AMSGrad & \texttt{($\eta$, bs)} \\
\midrule
\midrule

train acc.  & 92.101  $\pm$  0.083 & 92.763	$\pm$ 0.074 &  & 99.574 $\pm$ 0.021 & 99.598 $\pm$ 0.065   \\
val. acc. & 91.658 $\pm$ 0.291 & 91.847	$\pm$ 0.224 & & \textbf{ \color{blue} 80.970	$\pm$ 0.185 \color{black} } & 83.167 $\pm$ 0.052
& \\
test acc.   & 92.018 $\pm$ 0.103 & \textbf{ \color{red} 92.378 $\pm$ 0.073 \color{black} } & $(1\textit{e-}4,100)$ & \textbf{ \color{red} 80.820 $\pm$ 0.477 \color{black} } & \textbf{ \color{red} 82.690	$\pm$ 0.136 \color{black} } & $(1\textit{e-}4,20)$ \\
train loss  & 0.281	$\pm$ 0.002 & 0.260	$\pm$ 0.001 & & 0.014 $\pm$ 0.001 & 0.012 $\pm$ 0.001
&  \\
val. loss   & 0.296	$\pm$ 0.007 & 0.288	$\pm$ 0.005 &  & 0.969 $\pm$ 0.014 & 0.971 $\pm$ 0.032 \\
\bottomrule

train acc.  & 91.953 $\pm$ 0.077 & 92.603 $\pm$ 0.048 &  & 99.588 $\pm$ 0.048 & 99.656 $\pm$ 0.056  \\
val. acc. & 91.660 $\pm$ 0.249 & 91.982	$\pm$ 0.241& & 80.800 $\pm$ 0.421 & \textbf{ \color{blue} 83.283 $\pm$ 0.452 \color{black} } \\
test acc.   & 92.020 $\pm$ 0.054 & 92.340 $\pm$ 0.070 & $(1\textit{e-}4,128)$ & 80.303 $\pm$ 0.196& 82.607 $\pm$ 0.225 & $(1\textit{e-}4,32)$ \\
train loss  & 0.286 $\pm$ 0.002 & 0.264	$\pm$ 0.001&  & 0.014 $\pm$ 0.001& 0.011 $\pm$ 0.001 \\
val. loss   & 0.300	$\pm$ 0.006 & 0.287	$\pm$ 0.005&  & 0.954 $\pm$ 0.013 & 0.915 $\pm$ 0.013  \\
\bottomrule

train acc.  & 93.443 $\pm$ 0.048 & 93.465 $\pm$ 0.077 &  & 96.708 $\pm$ 1.159 & 95.832 $\pm$ 0.407  \\
val. acc. & 92.203 $\pm$ 0.170 & 91.855	$\pm$ 0.173 &  & 79.910	$\pm$ 0.452 & 81.157 $\pm$ 0.645 \\
test acc.   & 92.534 $\pm$ 0.104 & 92.360 $\pm$ 0.058 & $(1\textit{e-}3,100)$ & 79.690 $\pm$ 0.329& 80.847 $\pm$ 0.348
& $(1\textit{e-}3,20)$ \\
train loss  & 0.235	$\pm$ 0.001 & 0.234	$\pm$ 0.001 &  & 0.095 $\pm$ 0.032& 0.121 $\pm$ 0.011 \\
val. loss   & 0.286	$\pm$ 0.005 & 0.299 $\pm$ 0.008 &  & 0.925 $\pm$ 0.142& 0.808 $\pm$ 0.018 \\
\bottomrule

train acc.  & 93.392 $\pm$ 0.045 & 93.463 $\pm$ 0.057 &  & 64.700 $\pm$ 38.854& 96.707 $\pm$ 0.465 \\
val. acc. & \textbf{ \color{blue} 92.358 $\pm$ 0.157 \color{black} } & \textbf{ \color{blue} 92.000 $\pm$ 0.119 \color{black} } &  & 55.620 $\pm$	32.474& 81.537 $\pm$ 0.063 \\
test acc.   & \textbf{ \color{red} 92.540 $\pm$ 0.089 \color{black} } &92.310 $\pm$ 0.053 & $(1\textit{e-}3,128)$ & 55.313 $\pm$ 32.044& 81.370	$\pm$ 0.156& $(1\textit{e-}3,32)$ \\
train loss  & 0.239	$\pm$ 0.001 & 0.234	$\pm$ 0.001&  & 0.915 $\pm$ 0.981& 0.096 $\pm$ 0.014 \\
val. loss   & 0.274	$\pm$ 0.004 & 0.295	$\pm$ 0.004&  & 1.271 $\pm$ 0.730& 0.829 $\pm$ 0.016 \\
\bottomrule

\end{tabular}
\centering
\vspace{+0.5em}
\caption{Batch-size $\&$ learning rate evolution for \texttt{LR-MNIST}, \texttt{CNN-CIFAR10}, \eqref{AdaptiveAcceleratedMomentumMethodShiftedUpdates} and \emph{AMSGrad}. We have set (\texttt{epochs}, \texttt{n$\_$runs}) to $(40, 5)$ for \texttt{LR-MNIST} and $(30, 3)$ for \texttt{CNN-CIFAR10}. For \eqref{AdaptiveAcceleratedMomentumMethodShiftedUpdates} we considered $(M, \mu, \nu, \tilde{\gamma}) = (0.75, 0.5, 0.5, 0.75)$. }\label{table_models_AAMMSU_AMSGrad_batch_evolution}
\end{table}


\setlength{\tabcolsep}{10pt}
\renewcommand{\arraystretch}{0.75}
\begin{table}[htbp]
\begin{tabular}{>{\small}l || >{\scriptsize}c | >{\scriptsize}c  | >{\scriptsize}c ||
>{\small}c  }
\toprule
 & \multicolumn{1}{c}{\texttt{epochs=15}} & \multicolumn{1}{c}{\texttt{epochs=35}}& \multicolumn{1}{c}{\texttt{epochs=50}}\\
\cmidrule(lr){2-2} \cmidrule(lr){3-3} \cmidrule(lr){4-4} \cmidrule(lr){5-5}
\textbf{Metrics} &Results  & Results   & Results  & \texttt{$(M, \mu, \nu, \tilde{\gamma})$}      \\
\midrule
\midrule
train acc.  & 92.373 $\pm$ 0.033& 92.912 $\pm$ 0.02& 93.114 $\pm$ 0.024				
&   \\
val. acc. & 91.822 $\pm$ 0.084 &92.13 $\pm$ 0.087  & 92.230 $\pm$ 0.147	&  \\
test acc.   &92.202 $\pm$ 0.059& 92.392 $\pm$ 0.079 &92.486	$\pm$ 0.107	 & $(0.25,0.5,0.5,0.75)$ \\
train loss  & 0.273 $\pm$ 0.001 & 0.254 $\pm$ 0.001 & 0.247 $\pm$ 0.002 &   \\
val. loss   & 0.291 $\pm$ 0.005 & 0.283 $\pm$ 0.006 & 0.281 $\pm$ 0.006  &   \\
\bottomrule

train acc.  &92.330 $\pm$ 0.056&  92.903 $\pm$ 0.056& 93.125 $\pm$ 0.060 &   \\
val. acc    &91.957 $\pm$ 0.266& 92.267 $\pm$ 0.265 & 92.380 $\pm$ 0.295& \\
test acc.   &92.198 $\pm$ 0.070&  92.496 $\pm$ 0.056 & 92.560 $\pm$ 0.056&$(0.25, 0.75, 0.5, 0.95)$\\
train loss  &0.274 $\pm$ 0.002& 0.255 $\pm$	0.002 &0.248 $\pm$ 0.002 & \\
val. loss   &0.284 $\pm$ 0.007& 0.276 $\pm$	0.007 & 0.275 $\pm$	0.007&\\
\bottomrule

train acc.  & 92.988 $\pm$ 0.060 & 93.369 $\pm$ 0.042& 93.493 $\pm$	0.037&  \\
val. acc    & 92.098 $\pm$ 0.203  & 92.202 $\pm$ 0.167& 92.117 $\pm$ 0.209& \\
test acc.   & 92.488 $\pm$ 0.092 &92.65 $\pm$ 0.059 & \textbf{ \color{red} 92.658 $\pm$ 0.084 \color{black} } & $(0.75, 0.5, 0.5, 0.75)$\\
train loss  & 0.252 $\pm$ 0.001 & 0.238 $\pm$ 0.001& 	0.233 $\pm$	0.001& \\
val. loss   & 0.285	$\pm$ 0.004  &	0.285 $\pm$	0.005 & 	0.288 $\pm$	0.005&\\
\bottomrule

train acc.  & 92.993 $\pm$ 0.053 & 93.375 $\pm$	0.033				
& 93.544 $\pm$	0.055 & \\
val. acc    & 92.315 $\pm$ 0.105 &92.312 $\pm$	0.096 & 92.290 $\pm$ 0.113& \\
test acc.   & 92.486 $\pm$ 0.056  & 92.57 $\pm$	0.101&92.592 $\pm$ 0.090 &$(0.75, 0.75,	0.5, 0.95)
$ \\
train loss  &0.252 $\pm$ 0.001 & 0.238 $\pm$ 0.001& 0.232 $\pm$	0.001& \\
val. loss   &0.285 $\pm$ 0.006  & 0.286	$\pm$ 0.006& 0.288 $\pm$ 0.007&\\
\bottomrule
train acc.  & 93.073 $\pm$ 0.064 & 93.501 $\pm$	0.055& 93.595 $\pm$	0.044&  \\
val. acc    & 92.072 $\pm$ 0.150  & 92.092 $\pm$ 0.188& 92.055 $\pm$ 0.185& \\
test acc.   & 92.492 $\pm$ 0.047  & 92.432 $\pm$ 0.148&92.382 $\pm$	0.126 & $(1.25, 0.5, 0.5, 0.75)
$ \\
train loss  & 	0.248 $\pm$	0.002& 0.235 $\pm$ 0.002&0.230 $\pm$ 0.002 & \\
val. loss   & 0.283	$\pm$ 0.006 & 0.288	$\pm$ 0.009& 0.291 $\pm$ 0.009 &\\
\bottomrule

train acc.  & 93.104 $\pm$ 0.115 &93.468 $\pm$ 0.061& 93.568 $\pm$ 0.072& \\
val. acc    & 92.232 $\pm$ 0.232 & 92.142 $\pm$	0.316& 92.213 $\pm$	0.307& \\
test acc.   &92.496 $\pm$ 0.099 &92.508 $\pm$ 0.093 & 92.468 $\pm$ 0.135& $(1.25, 0.75, 0.5, 0.95)
$\\
train loss  &0.249 $\pm$ 0.003 & 0.236 $\pm$ 0.003& 0.230 $\pm$	0.003& \\
val. loss   & 	0.282 $\pm$	0.011 & 0.287 $\pm$	0.012& 0.290 $\pm$	0.013&\\
\bottomrule

train acc.  & 92.764 $\pm$ 0.052& 93.248 $\pm$ 0.056&93.395	$\pm$ 0.071 & \\
val. acc    & 91.948 $\pm$ 0.176 &92.182 $\pm$ 0.196 & 92.265 $\pm$	0.215& \\
test acc.   & 92.410 $\pm$ 0.124  & 92.514 $\pm$ 0.066& 92.608 $\pm$ 0.077& $(0.25, 0.5,	0.9, 0.75)$\\
train loss  & 	0.260 $\pm$	0.002& 0.244 $\pm$ 0.002&0.238 $\pm$ 0.002 & \\
val. loss   & 0.287 $\pm$ 0.008 &0.283 $\pm$ 0.009 & 0.284 $\pm$ 0.009&\\
\bottomrule

train acc.  &92.723 $\pm$ 0.069 & 93.22	$\pm$ 0.076				
 &93.401 $\pm$ 0.048 & \\
val. acc    & 92.173 $\pm$ 0.203 & 92.383 $\pm$ 0.125& \textbf{ \color{blue} 92.392 $\pm$ 0.124 \color{black} } & \\
test acc.   & 92.432 $\pm$ 0.057 &92.584 $\pm$ 0.095 & 92.632 $\pm$	0.075& $(0.25, 0.75, 0.9, 0.95)
$\\
train loss  &	0.259 $\pm$	0.002   & 0.243 $\pm$ 0.002&0.238 $\pm$	0.002 & \\
val. loss   & 0.287	$\pm$ 0.007  & 0.285 $\pm$ 0.007&0.285 $\pm$ 0.008 &\\
\bottomrule
train acc.  &93.048 $\pm$ 0.029& 93.435 $\pm$ 0.098				
&93.523 $\pm$ 0.049
 &  \\
val. acc    & 92.283 $\pm$ 0.255 & 92.152 $\pm$	0.176& 	92.178 $\pm$ 0.183& \\
test acc.   & 92.434 $\pm$ 0.047 & 92.418 $\pm$	0.206& 92.442 $\pm$	0.079	& $(0.75, 0.5,	0.9, 0.75)
$\\
train loss  &0.249 $\pm$ 0.002 &0.236 $\pm$	0.002 &0.230 $\pm$ 0.002 & \\
val. loss   & 0.283 $\pm$ 0.008& 0.29 $\pm$	0.007 & 0.295 $\pm$	0.008&\\
\bottomrule

train acc.  &93.103 $\pm$ 0.063& 93.48 $\pm$ 0.039& 93.575 $\pm$ 0.055&  \\
val. acc    & 92.020 $\pm$ 0.153  & 91.965 $\pm$ 0.19&91.913 $\pm$	0.212 & \\
test acc.   & 92.428 $\pm$ 0.112 & 	92.388 $\pm$ 0.162& 92.306 $\pm$ 0.123& $(0.75, 0.75, 0.9, 0.95)$\\
train loss  &0.248 $\pm$ 0.002  & 0.235 $\pm$ 0.002& 0.230 $\pm$ 0.002& \\
val. loss   &0.289	$\pm$ 0.007  & 0.293 $\pm$ 0.009&0.297 $\pm$ 0.011 &\\
\bottomrule

train acc.  & 93.025 $\pm$ 0.051&93.402	$\pm$ 0.077& 93.445	$\pm$ 0.101& \\
val. acc    & 92.115 $\pm$ 0.243  & 92.008 $\pm$ 0.175&91.837 $\pm$	0.041 & \\
test acc.   & 92.298 $\pm$ 0.123  & 92.322 $\pm$ 0.041&92.144 $\pm$	0.201 & $(1.25, 0.5,	0.9, 0.75)$\\
train loss  & 0.248	$\pm$ 0.001 &0.236 $\pm$ 0.001 & 0.232 $\pm$ 0.001& \\
val. loss   & 0.293	$\pm$ 0.006 &	0.298 $\pm$	0.003 & 0.306 $\pm$	0.005&\\
\bottomrule
train acc.  & 93.027 $\pm$ 0.055 & 93.403 $\pm$	0.064			
&93.492	$\pm$ 0.066 &  \\
val. acc    & 91.960 $\pm$ 0.189& 91.94	$\pm$ 0.279& 91.818	$\pm$ 0.204& \\
test acc.   & 92.394 $\pm$ 0.131 &92.406 $\pm$ 0.131 & 92.310 $\pm$	0.169& $(1.25, 0.75, 0.9, 0.95)$\\
train loss  & 0.248	$\pm$ 0.002 &0.236 $\pm$ 0.002 & 0.231	$\pm$ 0.002& \\
val. loss   & 0.297 $\pm$ 0.007 &	0.301 $\pm$	0.009 &0.310 $\pm$	0.009   &\\

\bottomrule
\end{tabular}
\centering
\vspace{+0.5em}
\caption{Grid search for \texttt{LR-MNIST} and \eqref{AdaptiveAcceleratedMomentumMethodShiftedUpdates}. The batch size \texttt{bs} was set to $128$, the learning rate $\eta$ to $1\text{e-}3$ and \texttt{n$\_$runs}=$5$.}\label{table_LR_MNIST_AAMMSU}
\end{table}


\begin{table}[htbp]
\begin{tabular}{>{\small}l || >{\scriptsize}c | >{\scriptsize}c  | >{\scriptsize}c ||
>{\small}c  }
\toprule
 & \multicolumn{1}{c}{\texttt{epochs=10}} & \multicolumn{1}{c}{\texttt{epochs=17}}& \multicolumn{1}{c}{\texttt{epochs=30}}\\
\cmidrule(lr){2-2} \cmidrule(lr){3-3} \cmidrule(lr){4-4} \cmidrule(lr){5-5}
\textbf{Metrics} &Results  & Results   & Results  & \texttt{$(M, \mu, \nu, \tilde{\gamma})$}      \\
\midrule
\midrule

train acc.  
&91.040 $\pm$	0.477
&97.793 $\pm$	0.169
&99.426 $\pm$	0.145
&  
\\
val. acc. 
& 80.540 $\pm$	0.925
&81.827 $\pm$	0.262
& \textbf{ \color{blue} 82.243 $\pm$ 0.331 \color{black} }
&  \\
test acc.  
&80.220 $\pm$	1.017
&81.407 $\pm$ 	0.194
& \textbf{ \color{red} 81.960 $\pm$	0.434 \color{black} }
& $(0.15, 0.5, 0.5, 0.65)$ \\
train loss  
&0.259 $\pm$	0.015
&0.067 $\pm$	0.006
&	0.019 $\pm$	0.004
&\\
val. loss   
& 0.609 $\pm$	0.044
&0.737 $\pm$	0.023
&	0.874 $\pm$	0.025
&\\

\bottomrule

train acc.  
&91.069 $\pm$	0.114
& 97.441 $\pm$	0.035
&99.277 $\pm$	0.053
&   \\
val. acc. 
& 80.400 $\pm$	0.658
&81.103 $\pm$	0.271
& 82.090 $\pm$	0.706
&  \\
test acc.  
& 80.583 $\pm$	0.456
& 81.093 $\pm$	0.141 
&81.763 $\pm$	0.152
& $(0.15, 0.75, 0.5, 0.85)$ \\
train loss  
&	0.254 $\pm$	0.002
&	0.076 $\pm$	0.001
& 0.022 $\pm$	0.001
&\\
val. loss   
& 0.620 $\pm$	0.023
&0.746 $\pm$	0.018
&0.880 $\pm$	0.029
&\\

\bottomrule
train acc.
&82.785 $\pm$	1.482 

&90.398 $\pm$	1.583
&95.676 $\pm$	1.104
&   \\
val. acc. 
&78.967 $\pm$	0.127 
&80.063 $\pm$0.334	
& 80.467 $\pm$	0.344
&  \\
test acc.  
&78.697 $\pm$	0.625
&	79.917 $\pm$	0.413	
&80.597 $\pm$	0.217
& $(0.65, 0.5, 0.5, 0.65)$ \\
train loss  
&0.486 $\pm$	0.043
&	0.269 $\pm$	0.046
&	0.123 $\pm$	0.030
&\\
val. loss   
& 0.621 $\pm$	0.002
&0.656 $\pm$	0.033
&0.787 $\pm$	0.055
&\\
\bottomrule

train acc.  
&77.468 $\pm$	3.858
&86.149 $\pm$	2.966
&93.683 $\pm$	1.568
&   \\
val. acc. 
& 75.650 $\pm$	2.151
&78.213 $\pm$	0.769
& 78.943 $\pm$	1.034
&  \\
test acc.  
&75.613 $\pm$	1.986
&78.333 $\pm$	0.684
&	79.133 $\pm$	1.088
&  $(0.65, 0.75, 0.5,	0.85)$\\
train loss  
&0.638 $\pm$	0.108
&0.392 $\pm$	0.083
&0.182 $\pm$	0.046
&\\
val. loss   
& 0.703 $\pm$	0.056
&0.660 $\pm$	0.010
&0.752 $\pm$	0.025
&\\

\bottomrule
train acc.  
&49.383 $\pm$	27.911
&54.792 $\pm$	31.909
&60.232 $\pm$	35.769
&   \\
val. acc. 
& 49.470 $\pm$	27.856
&52.807 $\pm$	30.496
& 	54.310 $\pm$	31.460
&  \\
test acc.  
&49.533 $\pm$	27.956
&52.873 $\pm$	30.317
&54.207 $\pm$	31.260
&  $(1.75, 0.5, 0.5,	0.65)$\\
train loss  
&1.349 $\pm$	0.675
&1.192 $\pm$	0.785
&1.038 $\pm$	0.895
&\\
val. loss   
& 1.354 $\pm$	0.671
& 1.260 $\pm$	0.738
&1.252 $\pm$	0.743
&\\
\bottomrule

train acc.  
&47.854 $\pm$	27.033
&54.528 $\pm$	31.776
&60.871 $\pm$	36.142
&   \\
val. acc. 
& 47.997 $\pm$	27.159
&52.657 $\pm$	30.236
& 54.223 $\pm$	31.476
&  \\
test acc.  
&	48.153 $\pm$	27.114
&52.420 $\pm$	30.006
&	54.167 $\pm$	31.239
&$(1.75, 0.75, 0.5,	0.85)$  \\
train loss  
&	1.390 $\pm$	0.649
&	1.201 $\pm$	0.782
&1.028 $\pm$	0.905 
&\\
val. loss   
& 1.393 $\pm$	0.645
&	1.285 $\pm$	0.720
&1.300 $\pm$	0.710
&\\

\bottomrule
train acc.  
&87.625 $\pm$	0.423
& 95.444 $\pm$	0.451
&98.803 $\pm$	0.123
&   \\
val. acc. 
& 77.387 $\pm$	1.994
&	81.223 $\pm$	0.092
& 81.417 $\pm$	0.343
&  \\
test acc.  
& 77.563 $\pm$      1.694
&81.567 $\pm$	0.379
&81.350 $\pm$	0.100
& $(0.15, 0.5, 0.75,	0.65)$ \\
train loss  
&	0.354 $\pm$	0.014
&0.132 $\pm$	0.012
&0.037 $\pm$	0.004
&\\
val. loss   
& 	0.694 $\pm$	0.068
& 	0.691 $\pm$	0.017
&	0.879 $\pm$	0.024
&\\
\bottomrule

train acc.  
&88.790 $\pm$	0.211
&95.964 $\pm$	0.201
&98.596 $\pm$	0.061
&   \\
val. acc. 
& 80.193 $\pm$	0.435
&80.527 $\pm$	0.176
& 81.707 $\pm$	0.283
&  \\
test acc.  
&80.397 $\pm$	0.601
&80.640 $\pm$	0.396
&81.613 $\pm$	0.274
&  $(0.15, 0.75, 0.75, 0.85)$\\
train loss  
&0.317 $\pm$	0.005
&0.118 $\pm$	0.007
&	0.041 $\pm$	0.002
&\\
val. loss   
& 0.610 $\pm$	0.016
&0.749 $\pm$	0.013
&0.870 $\pm$	0.013
&\\

\bottomrule
train acc.  

&71.059 $\pm$	17.238
&84.168 $\pm$	11.909
&91.692 $\pm$	7.770
&\\
val. acc. 
& 68.017 $\pm$	11.928
&74.087 $\pm$	4.956
& 75.970 $\pm$	3.026
&  \\
test acc.  
&	67.677 $\pm$	12.133
& 73.993 $\pm$	4.365
&	75.920 $\pm$	3.176
& $(0.65, 0.5, 0.75,	0.65)$ \\
train loss  
&0.793 $\pm$	0.458
&	0.444 $\pm$	0.334
&0.236 $\pm$	0.218
&\\
val. loss   
& 0.927 $\pm$	0.277
&0.874 $\pm$	0.052
&	1.016 $\pm$	0.155
&\\

\bottomrule

train acc.  
&72.145 $\pm$	3.372
&80.387 $\pm$	2.780
&88.836 $\pm$	2.528
&   \\
val. acc. 
& 72.347 $\pm$	2.705
&	76.380 $\pm$	1.565
& 78.457 $\pm$	1.336
&  \\
test acc.  
&	71.990 $\pm$	2.408
&75.790 $\pm$	1.764
&78.183 $\pm$	1.385
&$(0.65, 0.75, 0.75, 0.85)$  \\
train loss  
&0.781 $\pm$	0.091
&0.550 $\pm$	0.077
&0.313 $\pm$	0.068
&\\
val. loss   
& 0.802 $\pm$	0.074
&	0.698 $\pm$	0.050
&0.709 $\pm$	0.012
&\\
\bottomrule
 
train acc.  
&25.558 $\pm$	22.030
&29.871 $\pm$	28.057
&34.493 $\pm$	34.763
&   \\
val. acc. 
& 26.060 $\pm$	23.427
&28.090 $\pm$	25.958
& 28.953 $\pm$	27.695
&  \\
test acc.  
&26.323 $\pm$	23.085
&28.153 $\pm$	25.673
&29.020 $\pm$	26.898
&  $(1.75, 0.5, 0.75, 0.65)$\\
train loss  
&1.941 $\pm$	0.512
&1.823 $\pm$	0.679
&1.694 $\pm$	0.860
&\\
val. loss   
& 1.921 $\pm$	0.540
&	1.880 $\pm$	0.598
&	1.922 $\pm$	0.539
&\\	
\bottomrule
train acc.  
&38.581	$\pm$ 20.711
&44.968 $\pm$	25.348
&51.437 $\pm$	29.805
&   \\
val. acc. 
& 40.413 $\pm$	21.829
&45.950 $\pm$	25.717
& 	49.767 $\pm$	28.492
&  \\
test acc.  
&39.917 $\pm$	21.647
&	45.553 $\pm$	25.402
&49.357 $\pm$	28.114
& $(1.75, 0.75, 0.75, 0.85)$ \\
train loss  
&1.642 $\pm$	0.480
&1.465 $\pm$	0.605
&	1.289 $\pm$	0.728
&\\
val. loss   
& 1.596 $\pm$	0.511
&	1.447 $\pm$	0.611
&	1.350 $\pm$	0.679 
&\\
\bottomrule
\end{tabular}
\centering
\vspace{+0.5em}
\caption{Grid search for \texttt{CNN-CIFAR10} and \eqref{AdaptiveAcceleratedMomentumMethodShiftedUpdates}. The batch size \texttt{bs} was set to $20$, the learning rate $\eta$ to $1\text{e-}3$ and \texttt{n$\_$runs}=$3$.}\label{table_CNN_CIFAR10_AAMMSU}
\end{table}


\setlength{\tabcolsep}{10pt}
\renewcommand{\arraystretch}{0.75}
\begin{table}[htbp]
\begin{tabular}{>{\small}l || >{\scriptsize}c | >{\scriptsize}c || >{\scriptsize}c | >{\scriptsize}c  }
\toprule
& \multicolumn{1}{c}{\texttt{VGG-CIFAR10}} 
& \multicolumn{1}{c}{\texttt{ResNet-CIFAR10}}\\
\cmidrule(lr){1-1} \cmidrule(lr){2-2} \cmidrule(lr){3-3} \cmidrule(lr){4-4} \cmidrule(lr){5-5}
\textbf{Metrics}  & \textbf{Results} & \textbf{Results} & $(M, \mu, \nu, \tilde{\gamma})$ & \texttt{epochs} \\
\midrule
\midrule

train acc.  & 99.838 $\pm$  0.134 & 99.925 $\pm$ 0.042 & &    \\
val. acc. & 82.116 $\pm$ 1.154 & 81.878 $\pm$ 1.293 & & \\
test acc.  & 81.774	$\pm$ 0.995 & 81.696 $\pm$ 1.241 & (2, 0.5, 0.5, 0.5) & 50 \\
train loss & 0.005 $\pm$ 0.004 & 0.003 $\pm$ 0.001 & & \\
val. loss & 1.047 $\pm$ 0.069 & 1.017 $\pm$ 0.079 & & \\
\bottomrule

train acc.  & 99.943 $\pm$ 0.036 & 99.959 $\pm$ 0.023 & &    \\
val. acc. & 82.248 $\pm$ 1.379 & 83.774 $\pm$ 0.964 & & \\
test acc.  & 81.686	$\pm$ 1.439 & 83.084 $\pm$ 0.772 & (1, 0.5, 0.5, 0.5) & 50 \\
train loss & 0.002 $\pm$ 0.001 & 0.002 $\pm$ 0.001 & & \\
val. loss & 1.019 $\pm$ 0.081 & 0.838 $\pm$ 0.076 & & \\
\bottomrule
\bottomrule

train acc.  & 99.998 $\pm$ 0.003 & 99.997 $\pm$ 0.003 & &    \\
val. acc. & 83.802 $\pm$ 0.281 & 83.326 $\pm$ 0.660 & & \\
test acc.  & 83.434	$\pm$ 0.169 & 82.998 $\pm$  0.583 & (2, 0.5, 0.5, 0.5) & 75 \\
train loss & 0.000 $\pm$ 0.000 & 0.000 $\pm$ 0.000 & & \\
val. loss & 0.959 $\pm$ 0.035 & 0.922 $\pm$ 0.029 & & \\
\bottomrule

train acc.  & 99.999 $\pm$  0.001 & 99.994 $\pm$ 0.005 & &    \\
val. acc. & 84.140 $\pm$ 0.203 & 84.712	 $\pm$ 0.268 & & \\
test acc.  & 83.600 $\pm$ 0.182 & 84.206 $\pm$ 0.214 & (1, 0.5, 0.5, 0.5) & 75 \\
train loss & 0.000 $\pm$ 0.000 & 0.000 $\pm$ 0.000 & & \\
val. loss & 0.906 $\pm$ 0.012 & 0.783 $\pm$ 0.019 & & \\
\bottomrule
\bottomrule

train acc.  & 99.998 $\pm$ 0.002 & 99.998 $\pm$ 0.002 & &    \\
val. acc. & 83.946 $\pm$ 0.290 & 83.474	 $\pm$ 0.754 & & \\
test acc.  & 83.472 $\pm$ 0.194 & 83.160 $\pm$ 0.564 & (2, 0.5, 0.5, 0.5) & 100 \\
train loss & 0.000 $\pm$ 0.000 & 0.000 $\pm$ 0.000 & & \\
val. loss & 0.975 $\pm$ 0.036 & 0.929 $\pm$ 0.027 & & \\
\bottomrule

train acc.  & 99.999 $\pm$ 0.001 & 99.995 $\pm$ 0.003 & &    \\
val. acc. & 84.208 $\pm$ 0.303 & 84.770 $\pm$ 0.243 & & \\
test acc.  & 83.746 $\pm$ 0.288 & 84.220 $\pm$ 0.305 & (1, 0.5, 0.5, 0.5) & 100 \\
train loss & 0.000 $\pm$ 0.000 & 0.000 $\pm$ 0.000 & & \\
val. loss & 0.916  $\pm$ 0.016 & 0.787 $\pm$ 0.025 & & \\
\bottomrule
\bottomrule

train acc.  & 99.999 $\pm$ 0.002 & 99.997 $\pm$ 0.003 & &    \\
val. acc. & 83.940 $\pm$ 0.303 & 83.458 $\pm$ 0.670 & & \\
test acc.  & 83.546 $\pm$ 0.184 & 83.132 $\pm$ 0.594 & (2, 0.5, 0.5, 0.5) & 150 \\
train loss & 0.000  $\pm$  0.000 & 0.000 $\pm$ 0.000 & & \\
val. loss & 0.975 $\pm$ 0.036 & 0.916 $\pm$ 0.024 & & \\
\bottomrule

train acc.  & 100.000 $\pm$ 0.000 & 99.997 $\pm$ 0.003 & &    \\
val. acc. & \textbf{ \color{blue} 84.268 $\pm$ 0.144 \color{black} } & 84.822 $\pm$ 0.284 & & \\
test acc.  & \textbf{ \color{red} 83.808  $\pm$  0.232 \color{black} } &  84.318	$\pm$ 0.223 & (1, 0.5, 0.5, 0.5) & 150 \\
train loss & 0.000 $\pm$ 0.000 & 0.000 $\pm$ 0.000 & & \\
val. loss & 0.920 $\pm$ 0.010 & 0.784 $\pm$ 0.021 & & \\
\bottomrule
\bottomrule

train acc.  & 100.000 $\pm$ 0.000 & 100.000 $\pm$ 0.001 & &    \\
val. acc. & 84.002 $\pm$ 0.322 & 83.468 $\pm$ 0.669 & & \\
test acc.  & 83.478 $\pm$ 0.107 & 83.202 $\pm$ 0.656 & (2, 0.5, 0.5, 0.5) & 175 \\
train loss & 0.000 $\pm$ 0.000 & 0.000 $\pm$ 0.000 & & \\
val. loss & 0.985 $\pm$ 0.028 & 0.918 $\pm$ 0.027 & & \\
\bottomrule

train acc.  & 100.000 $\pm$ 0.001 & 100.000 $\pm$ 0.001 & &    \\
val. acc. & 84.178 $\pm$ 0.238 & \textbf{ \color{blue} 84.854	$\pm$ 0.255 \color{black} } & & \\
test acc.  & 83.784 $\pm$ 0.263 & \textbf{ \color{red} 84.330 $\pm$ 0.232 \color{black} } & (1, 0.5, 0.5, 0.5) & 175 \\
train loss & 0.000 $\pm$ 0.000 & 0.000 $\pm$ 0.000 & & \\
val. loss & 0.915 $\pm$ 0.013 & 0.780 $\pm$ 0.023 & & \\
\bottomrule
\bottomrule

train acc.  & 100.000 $\pm$ 0.001 & 99.998 $\pm$ 0.002 & &    \\
val. acc. & 83.896 $\pm$ 0.353 & 83.382 $\pm$ 0.595 & & \\
test acc.  & 83.454 $\pm$ 0.160 & 83.204 $\pm$ 0.585 & (2, 0.5, 0.5, 0.5) & 200 \\
train loss & 0.000 $\pm$ 0.000 & 0.000 $\pm$ 0.000 & & \\
val. loss & 0.983 $\pm$ 0.030 & 0.934 $\pm$ 0.027 & & \\
\bottomrule

train acc.  & 99.999 $\pm$ 0.001 & 99.999 $\pm$ 0.002 & &    \\
val. acc. & 84.144 $\pm$ 0.233 & 84.728 $\pm$ 0.322 & & \\
test acc.  & 83.746 $\pm$ 0.207 & 84.272 $\pm$ 0.195 & (1, 0.5, 0.5, 0.5) & 200 \\
train loss & 0.000 $\pm$ 0.000 & 0.000 $\pm$ 0.000 & & \\
val. loss & 0.922 $\pm$ 0.024 & 0.786 $\pm$ 0.030 & & \\
\bottomrule

\end{tabular}
\centering
\vspace{+0.5em}
\caption{Grid search for \texttt{VGG-CIFAR10}, \texttt{ResNet-CIFAR10} and \eqref{AdaptiveAcceleratedMomentumMethodShiftedUpdates}. The batch size \texttt{bs} was set to $128$ and the initial value of the learning rate $\eta$ to $1\text{e-}3$, which finally decreases to $1\text{e-}6$ due to the scheduler. Also, we have taken \texttt{n$\_$runs}=$5$.}\label{table_VGG_ResNet_AAMMSU}
\end{table}


\setlength{\tabcolsep}{10pt}
\renewcommand{\arraystretch}{0.75}
\begin{table}[ht!]
\begin{tabular}{>{\small}l || >{\scriptsize}c | >{\scriptsize}c || >{\scriptsize}c   }
\toprule
& \multicolumn{1}{c}{\texttt{VGG-CIFAR10}} 
& \multicolumn{1}{c}{\texttt{ResNet-CIFAR10}}\\
\cmidrule(lr){1-1} \cmidrule(lr){2-2} \cmidrule(lr){3-3} \cmidrule(lr){4-4}
\textbf{Metrics}  & \textbf{Results} & \textbf{Results} & \texttt{epochs} \\
\midrule
\midrule

train acc.  & 99.979 $\pm$ 0.041 & 99.976 $\pm$ 0.005 &    \\
val. acc. & 84.452 $\pm$ 0.375 & 85.766	 $\pm$ 0.294 & \\
test acc.  & 84.046	$\pm$ 0.381 & 85.122 $\pm$ 0.441 & 50 \\
train loss & 0.001 $\pm$ 0.001 & 0.001 $\pm$ 0.000 & \\
val. loss & 0.913 $\pm$ 0.006 & 0.783 $\pm$ 0.016 & \\
\bottomrule
\bottomrule

train acc.  & 100.000 $\pm$ 0.001 & 100.000	$\pm$ 0.001 &    \\
val. acc. & 84.746 $\pm$ 0.234 & 86.328	$\pm$ 0.273 & \\
test acc.  & 84.276	$\pm$ 0.192 & 85.552 $\pm$ 0.273 & 75 \\
train loss & 0.000 $\pm$ 0.000 & 0.000 $\pm$ 0.000 & \\
val. loss & 0.893 $\pm$ 0.021 & 0.760 $\pm$ 0.026 & \\
\bottomrule
\bottomrule

train acc.  & 100.000 $\pm$ 0.000 & 100.000	$\pm$ 0.001 &    \\
val. acc. & 84.694 $\pm$ 0.196 & 86.316	$\pm$ 0.320 & \\
test acc.  & 84.290	$\pm$ 0.194 & 85.680 $\pm$ 0.394 & 100 \\
train loss & 0.000 $\pm$ 0.000 & 0.000 $\pm$ 0.000 & \\
val. loss & 0.908 $\pm$ 0.020 & 0.765 $\pm$ 0.020 & \\
\bottomrule
\bottomrule

train acc.  & 100.000 $\pm$ 0.000 & 100.000 $\pm$ 0.000 &    \\
val. acc. & 84.714 $\pm$ 0.210 & 86.424 $\pm$ 0.290 & \\
test acc.  & 84.270	$\pm$ 0.153 & \textbf{ \color{red} 85.820 $\pm$ 0.349 \color{black} } & 150 \\
train loss & 0.000 $\pm$ 0.000 & 0.000 $\pm$ 0.000 & \\
val. loss & 0.905 $\pm$ 0.021 & 0.763 $\pm$ 0.021 & \\
\bottomrule
\bottomrule

train acc.  & 100.000 $\pm$ 0.001 & 99.999 $\pm$ 0.001 &    \\
val. acc. & \textbf{ \color{blue} 84.784 $\pm$ 0.174 \color{black} } & 86.394 $\pm$ 0.363 & \\
test acc.  & \textbf{ \color{red} 84.344 $\pm$ 0.174 \color{black} } & 85.660 $\pm$ 0.364 & 175 \\
train loss & 0.000 $\pm$ 0.000 & 0.000 $\pm$ 0.000 & \\
val. loss & 0.900 $\pm$ 0.017 & 0.767 $\pm$ 0.025 & \\
\bottomrule
\bottomrule

train acc.  & 100.000 $\pm$ 0.001 & 99.999 $\pm$ 0.001 &    \\
val. acc. & 84.720 $\pm$ 0.200 & \textbf{ \color{blue} 86.450 $\pm$ 0.278 \color{black} } & \\
test acc.  & 84.344	$\pm$ 0.195 & 85.682 $\pm$ 0.333 & 200 \\
train loss & 0.000 $\pm$ 0.000 & 0.000 $\pm$ 0.000 & \\
val. loss & 0.901 $\pm$ 0.019 & 0.764 $\pm$ 0.019 & \\
\bottomrule

\end{tabular}
\centering
\vspace{+0.5em}
\caption{Grid search for \texttt{VGG-CIFAR10}, \texttt{ResNet-CIFAR10} and \emph{AMSGrad}. The batch size \texttt{bs} was set to $128$ and the initial value of the learning rate $\eta$ to $1\text{e-}3$, which finally decreases to $1\text{e-}6$ due to the scheduler. Also, we have taken \texttt{n$\_$runs}=$5$.}\label{table_grid_search_vgg_resnet_AMSGrad}
\end{table}


\setlength{\tabcolsep}{10pt}
\renewcommand{\arraystretch}{0.75}
\begin{table}[hb!]
\centering
\begin{tabular}{>{\small}l || >{\scriptsize}c | >{\scriptsize}c  | >{\scriptsize}c  || >{\scriptsize}p{0.3cm} | >{\scriptsize}p{1.2cm}}
\toprule
 & \multicolumn{1}{c}{\texttt{epochs=15}} & \multicolumn{1}{c}{\texttt{epochs=35}}& \multicolumn{1}{c}{\texttt{epochs=50}}\\
\cmidrule(lr){1-1} \cmidrule(lr){2-2} \cmidrule(lr){3-3} \cmidrule(lr){4-4} \cmidrule(lr){5-5} \cmidrule(lr){6-6} 
\textbf{Metrics} &Results  & Results & Results & \texttt{bs} & Model \\
\midrule
\midrule

train acc.  & 93.074 $\pm$ 0.037 & 93.459 $\pm$ 0.034 & 93.554 $\pm$ 0.064 & & \\
val. acc.   & \textbf{ \color{blue} 92.037 $\pm$ 0.301 \color{black} } & 91.968 $\pm$ 0.204 & 91.848 $\pm$ 0.267 & & \\
test acc.   & \textbf{ \color{red} 92.402 $\pm$ 0.195 \color{black} } & 92.284 $\pm$ 0.160 &  92.322 $\pm$ 0.145 & 128 & \texttt{LR-MNIST} \\
train loss  & 0.248 $\pm$ 0.002 & 0.234 $\pm$ 0.002 & 0.229 $\pm$ 0.002 & & \\
val. loss   & 0.290 $\pm$ 0.008 & 0.300	$\pm$ 0.008 &  0.303 $\pm$ 0.009 & & \\
\bottomrule
\end{tabular}

\smallskip

\begin{tabular}{>{\small}l || >{\scriptsize}c | >{\scriptsize}c  | >{\scriptsize}c || >{\scriptsize}p{0.3cm} | >{\scriptsize}p{1.2cm}}
\toprule
 & \multicolumn{1}{c}{\texttt{epochs=10}} & \multicolumn{1}{c}{\texttt{epochs=17}}& \multicolumn{1}{c}{\texttt{epochs=30}}\\
\cmidrule(lr){1-1} \cmidrule(lr){2-2} \cmidrule(lr){3-3} \cmidrule(lr){4-4} \cmidrule(lr){5-5} \cmidrule(lr){6-6}  
\textbf{Metrics} &Results  & Results & Results & \texttt{bs} & Model       \\
\midrule
\midrule

train acc.  & 77.888 $\pm$ 1.278 & 85.487 $\pm$ 0.958 & 93.043 $\pm$ 0.472 & & \\
val. acc.   & 75.863 $\pm$ 0.742 & 78.067 $\pm$ 0.763  &  \textbf{ \color{blue} 79.383 $\pm$ 0.585 \color{black} } & & \\
test acc.   & 75.697 $\pm$ 0.956 & 77.720 $\pm$ 0.747 & \textbf{ \color{red} 79.077 $\pm$ 0.319 \color{black} }
& 20 & \texttt{CNN-CIFAR10} \\
train loss  & 0.634	$\pm$ 0.033 & 0.414	$\pm$ 0.029 & 0.201 $\pm$ 0.013 & & \\
val. loss   & 0.692	$\pm$ 0.026 & 0.671	$\pm$ 0.020 &  0.738
 $\pm$ 0.031 & &  \\
\bottomrule

\end{tabular}
\centering
\vspace{+0.5em}
\caption{Grid search for \texttt{LR-MNIST}, \texttt{CNN-CIFAR10} and \emph{AMSGrad} with learning rate $\eta = 1\text{e-}3$. We have used \texttt{n$\_$runs}=$5$ for \texttt{LR-MNIST} and \texttt{n$\_$runs}=$3$ for \texttt{CNN-CIFAR10}.}\label{table_models_lr_cnn_AMSGrad}
\end{table}

\end{document}